\documentclass[final,1p,times,number]{elsarticle}

\usepackage{amsmath}
\usepackage{amsthm}
\usepackage{amssymb}

\usepackage{amsfonts}
\usepackage{graphicx}
\usepackage{epstopdf}
\usepackage{color}
\usepackage{bm}
\usepackage{multirow}
\usepackage{natbib}
\usepackage{hyperref}
\usepackage{cleveref}

\usepackage{color}
\usepackage{float}
\usepackage{caption}

\usepackage{geometry}

\newtheorem{remark}{Remark}
\newtheorem{definition}{Definition}[section]
\newtheorem{theorem}{Theorem}[section]

\usepackage[caption=false]{subfig}

\usepackage{ifpdf}
\ifpdf%
\usepackage{pdflscape}
\else
\usepackage{lscape}
\fi

\hypersetup{
	colorlinks,
	linkcolor={red!50!black},
	citecolor={blue!50!black},
	urlcolor={blue!80!black}
}

\def\e{\mbox{\boldmath $e$}}
\def\f{\mbox{\boldmath $f$}}

\def\h{\mbox{\boldmath $h$}}
\def\m{\mbox{\boldmath $m$}}

\def\0{\mbox{\boldmath $0$}}

\biboptions{sort&compress}

\begin{document}

\begin{frontmatter}
\title{Improved Energy Stable Symmetric Gauss-Seidel Projection Method for Micromagnetics Simulations}

\author[XJTLU]{Yingxi Miao\corref{cor1}}
\cortext[cor1]{First author.} 
\ead{Yingxi.Miao22@student.xjtlu.edu.cn}

\author[XJTLU]{Changjian Xie\corref{cor2}}
\cortext[cor2]{Corresponding author.} 
\ead{Changjian.Xie@xjtlu.edu.cn}

\address[XJTLU]{School of Mathematics and Physics, Xi'an-Jiaotong-Liverpool University, Re'ai Rd. 111, Suzhou, 215123, Jiangsu, China.}

\begin{abstract}

The Gauss–Seidel projection method (GSPM) constitutes an efficient and numerically stable numerical framework for micromagnetic simulations of ferromagnetic media. This scheme attains first-order temporal accuracy and second-order spatial accuracy. Fast Fourier transform (FFT) techniques can be incorporated to accelerate both the solution of the arising linear algebraic systems and the evaluation of stray magnetic fields. The conventional GSPM relies on a single-sided Gauss–Seidel iteration, which leverages the latest updated state variables associated with the heat-diffusion subproblem.
In this work, we develop a symmetric Gauss–Seidel projection method (SGSPM) that retains first-order temporal accuracy and second-order spatial consistency. The proposed symmetric variant exhibits superior stability properties relative to the standard GSPM. Specifically, SGSPM adopts a two-pass symmetric Gauss–Seidel iteration, where updated information from the heat-diffusion stage is fully exploited to rigorously guarantee discrete energy stability. We validate the performance of the devised scheme through numerical investigations of magnetization dynamic evolution and magnetic domain-wall propagation. Numerical evidence demonstrates that the improved symmetric scheme delivers enhanced stability for capturing magnetization motion dynamics.
\end{abstract}

\begin{keyword}
{Micromagnetic Simulation\sep Symmetric Gauss–Seidel Projection Scheme\sep Discrete Energy Stability\sep Weak Damping}
\end{keyword}

\end{frontmatter}

\section{Introduction}

Micromagnetics simulation based on the Landau–Lifshitz–Gilbert (LLG) equation is the fundamental numerical tool for investigating magnetization dynamics, including magnetic switching, domain wall motion, spin wave propagation, and thermal magnetization fluctuations in ferromagnetic materials . The LLG equation governs the temporal evolution of magnetization vectors with unit-length constraints, highly nonlinear gyromagnetic torque, and multi-field coupling effects involving exchange field, anisotropy field, demagnetization stray field, and external applied field. Numerically, the stiff, constrained, and strongly coupled characteristics of the LLG system pose severe challenges to traditional time-stepping schemes, including numerical instability, strict time-step restrictions, constraint violation, and excessive computational overhead for large-scale sparse systems.

The Gilbert damping parameter $\alpha$ characterizes the energy dissipation rate in magnetization dynamics. For conventional soft magnetic materials, such as permalloy (NiFe) and typical CoFeB thin films, the room-temperature damping coefficient generally ranges from $\alpha=0.01$ to $\alpha=0.05$. Nevertheless, numerous experimental studies have verified that high-quality magnetic materials can achieve ultra-low Gilbert damping \cite{schoen2016ultra} down to the order of $10^{-3}\sim 10^{-4}$, which is physically feasible and widely adopted in micromagnetic simulations.
Yttrium iron garnet (YIG) \cite{onbasli2014pulsed} is well known as the low-loss magnetic material with the smallest room-temperature damping. Epitaxial YIG thin films can reach an extremely low damping of $\alpha \approx 3.4\times 10^{-4}$, while bulk single-crystal YIG exhibits damping coefficients as low as $\alpha \sim 10^{-4}$. Such ultra-low-damping insulators are standard benchmark materials for investigating domain wall motion and magnon propagation.
For metallic ferromagnets, high-quality epitaxial FeCo alloys can achieve a record-low Gilbert damping of $\alpha \approx 8\times 10^{-4}$ at room temperature, which is nearly two orders of magnitude smaller than that of traditional soft magnetic films. Sputtered FeCo thin films \cite{wei2021ultralow} also stably maintain $\alpha < 0.001$ in practical experimental conditions. In addition, half-Heusler ferromagnetic alloys and ferrimagnetic GdFeCo thin films possess moderate low damping values around $0.001\sim 0.007$, which are frequently employed in ultrafast domain wall dynamic simulations.

Conventional numerical methods for micromagnetic simulations can be categorized into explicit, fully implicit, and semi-implicit frameworks. Explicit methods, such as the forward Euler and Runge–Kutta schemes, feature simple implementation and low per-step cost but suffer from stringent Courant–Friedrichs–Lewy (CFL) stability constraints, requiring extremely small time steps for stiff magnetic systems and resulting in low computational efficiency. Fully implicit methods achieve unconditional stability but necessitate solving large-scale nonlinear coupled systems at each time step, which introduces enormous computational complexity and poor scalability for micron-scale magnetic structures. To balance numerical stability, accuracy, and computational efficiency, semi-implicit fractional-step methods have become the dominant paradigm in modern micromagnetic simulation, among which the Gauss–Seidel Projection Method (GSPM) stands out as one of the most classic and robust algorithms.

The GSPM was first proposed for the LLG-based micromagnetic dynamics, integrating Gauss–Seidel iterative updating, fractional-step decomposition, and manifold projection techniques. The original GSPM was formally established by Wang et al. (2001) \cite{wang2001gauss}, who innovatively combined Gauss–Seidel iterative relaxation with a projection method to solve the constrained LLG equation. Different from traditional unified implicit discretization, the GSPM adopts a physical-field fractional-step strategy, decoupling the stiff gyromagnetic precession term, Gilbert damping term, and unit-length magnetization constraint in the LLG equation. The core innovation lies in replacing the solution of complex nonlinear coupled systems with a series of simple linear constant-coefficient heat equations and low-cost projection corrections. In the original framework, each complete time iteration of GSPM is divided into two key sub-steps. First, a Gauss–Seidel semi-implicit sweep is performed to update the magnetization field, which efficiently stabilizes the precession and damping evolution without solving nonlinear systems. Second, a manifold projection step is implemented to strictly enforce the unit-length constraint of magnetization vectors, eliminating the numerical drift widely existing in unconstrained iterative methods. Compared with Jacobi-type fractional methods, the Gauss–Seidel updating mechanism utilizes real-time updated intermediate variables in each sweep, significantly accelerating error decay and improving iterative convergence efficiency. The original GSPM achieves unconditional stability with respect to magnetic damping parameters, second-order spatial accuracy, and first-order temporal accuracy, forming a lightweight, robust, and easy-to-parallelize numerical framework for micromagnetic simulations. Early studies verified that GSPM effectively solves the core defects of classic micromagnetic algorithms. Explicit schemes are limited by linear stability constraints and fail in long-time dynamic simulations of low-damping magnetic systems. Fully implicit methods require iterative solving of variable-coefficient asymmetric linear systems, leading to huge computational costs for large-scale micron-sized magnetic domains. In contrast, GSPM only solves constant-coefficient linear systems in each time step, with computational complexity equivalent to that of scalar heat equation solving, which greatly reduces hardware resource consumption while maintaining unconditional numerical stability. This unique advantage enables GSPM to realize fully resolved numerical simulations of magnetization switching and domain evolution in microscale magnetic devices, which was difficult for previous mainstream algorithms.

Although the original GSPM exhibits excellent stability and simplicity, it still has inherent limitations in practical micromagnetic applications, including low temporal accuracy, redundant stray field updates, and degraded performance under extreme low/high damping conditions. In the past two decades, one has proposed a series of modified GSPM variants targeting these defects. The original GSPM suffers from accuracy degradation for micromagnetic systems with extremely small damping parameters, which are common in high-frequency spin dynamics and ultrafast magnetization reversal simulations. One of the major limitations of the original GSPM is its first-order temporal accuracy, which restricts the simulation efficiency of long-time dynamic evolution and high-frequency spin wave problems. Recent studies have focused on high-order extension of the GSPM framework. Li et al. (2026) \cite{li2026enhanced} proposed second-order accurate GSPM variants by combining biharmonic operator discretization and second-order backward differentiation formula (BDF2). The high-order GSPM retains the unconditional stability and low-complexity advantages of the original algorithm, upgrades the temporal accuracy to second order, and maintains the computational cost comparable to implicit scalar biharmonic equation solving. This improvement effectively compensates for the low-accuracy defect of traditional GSPM in long-term micromagnetic evolution simulations.

Another major limitation of the original GSPM is its energy stability with weak damping parameters.
In this paper, combined with the symmetric Gauss–Seidel (SGS) iterative theory, we further develop a symmetric Gauss–Seidel projection scheme. The standard single-direction Gauss–Seidel sweep introduces directional numerical bias and weak error accumulation in long-time iteration. The symmetric GSPM adopts forward–backward bidirectional sweeps to eliminate directional iteration deviation, enhances the numerical stability of weak diagonally dominant magnetic field systems, and suppresses residual oscillation in the late iteration stage. This variant improves the robustness of GSPM for ill-conditioned micromagnetic systems with non-strictly diagonal dominant coefficient matrices.

The rest of the paper is organized below. 
\Cref{sec: model} introduces the full governing equation of micromagnetics, then \Cref{sec:method} first revises the original GSPM method and further proposes a symmetric Gauss-Seidel projection method. We give the theoretical and numerical analysis for the stability and convergence in \Cref{sec:theory} for GSPM and SGSPM. \Cref{sec:experiments} presents the simulation of micromagnetic domain wall dynamics to study the domain wall motion. Finally, concluding remarks and perspectives for future work are given in \Cref{sec:conclusions}.

\section{The physical governing equation}
\label{sec: model}

The micromagnetic equation integrates gyromagnetic precession and dissipative relaxation \cite{Landau1935On,Brown1963micromagnetics}. Its nondimensionalized model is given by
\begin{align}\label{c1-large}
{\m}_t =-{\m}\times{\bm h}_{\text{eff}}-\alpha{\m}\times({\m}\times{\bm h}_{\text{eff}}),
\end{align}
where $\partial_t \m$ is the time derivative for the magnetization and ${\bm h}_{\text{eff}}$ is the term from the energy variation, say
\(\bm h_{\text{eff}} = -\delta E[\m]/\delta \m\). In this paper, we take the free energy below,
\begin{align*}
    E[\m]=\epsilon \int_{\Omega} |\nabla \m|^2\;dx+Q\int_{\Omega} (m_2^2+m_3^2)\;dx- \int_{\Omega} \h_s\cdot \m\;dx-\int_{\Omega}\h_e\cdot \m\;dx,
\end{align*}
(here $Q$ and $\epsilon$ are nondimensionalized parameters) and get the effective field below,
\begin{align*}
    {\bm h}_{\text{eff}}=\epsilon \Delta \m-Q(m_2\e_2+m_3\e_3)+\h_s+\h_e,
\end{align*}
where $\epsilon \Delta \m$ is the exchange field, $-Q(m_2\e_2+m_3\e_3)$ is anisotropy field, $\h_s$ is the stray field and $\h_e$ is the applied external field.

This equation is subject to the homogeneous Neumann boundary condition
\begin{equation}\label{boundary-large}
\frac{\partial{\m}}{\partial {\bm \nu}}\Big|_{\partial \Omega}=0,
\end{equation}
where \(\Omega \subset \mathbb{R}^d\) (\(d=1,2,3\)) denotes the bounded spatial domain occupied by the ferromagnetic material, and \(\bm \nu\) represents the unit outward normal vector on the domain boundary \(\partial \Omega\). It is noteworthy that this boundary condition is physically consistent for isolated ferromagnetic systems, as it inherently ensures the absence of magnetic surface charge—an essential prerequisite for accurately modeling unperturbed magnetic dynamics.

To achieve a comprehensive understanding of the LLG equation, it is essential to elucidate the physical essence of its key components. The magnetization field \(\m: \Omega \to \mathbb{R}^3\) is a three-dimensional vector field subject to the pointwise constraint \(|\m|=1\), which indicates a fundamental characteristic derived from the quantum mechanical alignment of electron spins in ferromagnetic materials. With respect to the right-hand side of \eqref{c1-large}, the first term describes the gyromagnetic precession effect, whereby magnetic moments undergo precessional motion around the effective magnetic field \(\bm h_{\text{eff}}\). The second term represents dissipative relaxation, where the parameter \(\alpha > 0\) stands for the dimensionless Gilbert damping coefficient, which quantifies the rate of energy transfer from the magnetic subsystem to the lattice structure. Consequently, the model \eqref{c1-large} inherently satisfies the norm-preserving constraint.


\section{Proposed method}\label{sec:method}



In the light of the linear algebra,
we consider the nonsingular linear system
\begin{equation}
A\boldsymbol{x} = \boldsymbol{b},\quad A\in\mathbb{R}^{n\times n},\quad \boldsymbol{x},\boldsymbol{b}\in\mathbb{R}^n,
\end{equation}
and we decompose the coefficient matrix $A$ additively as
\begin{equation}
A = D - L - U,
\end{equation}
where $D$ denotes the diagonal matrix consisting of diagonal entries of $A$, invertible (no zero pivots for iteration), $L$ be strictly lower triangular matrix, $L_{ij}=-A_{ij}$ for $i>j$, zero otherwise, $U$ be strictly upper triangular matrix, $U_{ij}=-A_{ij}$ for $i<j$, zero otherwise.

This splitting forms the foundation for constructing Gauss-Seidel and symmetric Gauss-Seidel iterative schemes. For Gauss-Seidel (GS) iteration, 
let $\boldsymbol{x}^{(k)} = \big(x_1^{(k)},x_2^{(k)},\dots,x_n^{(k)}\big)^T$ denote the iterate at step $k$. The core idea of GS is \textit{instant update}: when computing $x_i^{(k+1)}$, all preceding components $x_1^{(k+1)},\dots,x_{i-1}^{(k+1)}$ already updated in the current sweep are utilized, while trailing components retain values from the previous iteration $\boldsymbol{x}^{(k)}$.

Rewrite the $i$-th scalar equation of $A\boldsymbol{x}=\boldsymbol{b}$:
\[
a_{ii}x_i = b_i - \sum_{j=1}^{i-1}a_{ij}x_j - \sum_{j=i+1}^n a_{ij}x_j.
\]
Divide by $a_{ii}\neq 0$ and apply the instant-update rule to obtain the component-wise GS update:
\begin{equation}\label{eq:gs-component}
x_i^{(k+1)} = \frac{1}{a_{ii}}\left( b_i - \sum_{j=1}^{i-1}a_{ij}x_j^{(k+1)} - \sum_{j=i+1}^n a_{ij}x_j^{(k)} \right),\quad i=1,2,\dots,n.
\end{equation}
In matrix form of GS iteration, 
rearrange $A=D-L-U$ into
\[
(D-L)\boldsymbol{x}^{(k+1)} - U\boldsymbol{x}^{(k)} = \boldsymbol{b}.
\]
Since $D-L$ is invertible lower triangular, left-multiply by $(D-L)^{-1}$:
\begin{equation}\label{eq:gs-matrix}
\boldsymbol{x}^{(k+1)} = (D-L)^{-1}U \boldsymbol{x}^{(k)} + (D-L)^{-1}\boldsymbol{b}.
\end{equation}
Define the \textit{GS iteration matrix}:
\begin{equation}
G_{\text{GS}} = (D-L)^{-1}U, \quad \boldsymbol{f}_{\text{GS}} = (D-L)^{-1}\boldsymbol{b}.
\end{equation}
The unified linear iteration format reads
\[
\boldsymbol{x}^{(k+1)} = G_{\text{GS}} \boldsymbol{x}^{(k)} + \boldsymbol{f}_{\text{GS}}.
\]

Standard GS employs a unidirectional forward sweep ($1\to n$), which introduces directional bias. SGS eliminates this asymmetry by combining a forward GS sweep followed by a backward GS sweep ($n\to 1$) in one full iteration cycle. Let $\bar{\boldsymbol{x}}^{(k+1)}$ be the intermediate solution after forward GS.

\begin{itemize}
    \item Step 1: Forward GS sweep ($1\to n$):
\begin{align}
\bar{x}_i^{(k+1)} &= \frac{1}{a_{ii}}\left( b_i - \sum_{j=1}^{i-1}a_{ij}\bar{x}_j^{(k+1)} - \sum_{j=i+1}^n a_{ij}x_j^{(k)} \right),\\
\bar{\boldsymbol{x}}^{(k+1)} &= (D-L)^{-1}U \boldsymbol{x}^{(k)} + (D-L)^{-1}\boldsymbol{b}.
\end{align}
\item Step 2: Backward GS sweep ($n\to 1$):
Now update from the last component to the first, using intermediate $\bar{\boldsymbol{x}}^{(k+1)}$ as input. The matrix form of backward GS is
\[
\boldsymbol{x}^{(k+1)} = (D-U)^{-1}L \bar{\boldsymbol{x}}^{(k+1)} + (D-U)^{-1}\boldsymbol{b}.
\]
\item Global SGS Iteration Matrix:
Substitute the forward sweep expression into the backward update to eliminate $\bar{\boldsymbol{x}}^{(k+1)}$:
\begin{align*}
\boldsymbol{x}^{(k+1)} &= (D-U)^{-1}L \big[(D-L)^{-1}U \boldsymbol{x}^{(k)} + (D-L)^{-1}\boldsymbol{b}\big] + (D-U)^{-1}\boldsymbol{b}\\
&= \underbrace{(D-U)^{-1}L(D-L)^{-1}U}_{G_{\text{SGS}}} \boldsymbol{x}^{(k)} + \underbrace{\big[(D-U)^{-1}L(D-L)^{-1} + (D-U)^{-1}\big]\boldsymbol{b}}_{\boldsymbol{f}_{\text{SGS}}}.
\end{align*}
The \textit{SGS iteration matrix} is defined as
\begin{equation}\label{eq:sgs-matrix}
G_{\text{SGS}} = (D-U)^{-1}L(D-L)^{-1}U,
\end{equation}
with unified iteration form
\[
\boldsymbol{x}^{(k+1)} = G_{\text{SGS}} \boldsymbol{x}^{(k)} + \boldsymbol{f}_{\text{SGS}}.
\]
\end{itemize}

\begin{definition}
A linear iteration $\boldsymbol{x}^{(k+1)}=G\boldsymbol{x}^{(k)}+\boldsymbol{f}$ converges if and only if the spectral radius of iteration matrix satisfies $\rho(G)<1$. A sufficient condition is $\|G\|<1$ for any consistent matrix norm.
\end{definition}
Let $\boldsymbol{x}^*$ denote the exact solution satisfying $\boldsymbol{x}^*=G\boldsymbol{x}^*+\boldsymbol{f}$. The error vector is defined as $\boldsymbol{e}^{(k)} = \boldsymbol{x}^{(k)} - \boldsymbol{x}^*$, yielding the error propagation law:
\begin{equation}\label{eq:error-recur}
\boldsymbol{e}^{(k+1)} = G \boldsymbol{e}^{(k)},\quad \boldsymbol{e}^{(k)} = G^k \boldsymbol{e}^{(0)}.
\end{equation}

\begin{remark}[Convergence of Gauss-Seidel Iteration]
    \begin{theorem}\label{thm:gs-diagdom}
If $A$ is strictly row diagonally dominant, i.e., $|a_{ii}|>\sum_{j\neq i}|a_{ij}|$ for all $i$, then $\rho(G_{\text{GS}})<1$, and GS iteration converges.
\end{theorem}
\end{remark}

\begin{remark}[Convergence of Symmetric Gauss-Seidel Iteration]
\begin{theorem}\label{thm:sgs-diagdom}
If $A$ is strictly diagonally dominant, SGS iteration converges, with $\rho(G_{\text{SGS}})\le \rho(G_{\text{GS}})<1$.
\end{theorem}
\end{remark}

\subsection{Gauss-Seidel Projection Method}





The Gauss-Seidel projection method (GSPM) is given as the following three steps:

\paragraph{Step 1. Implicit Gauss-Seidel:}
\begin{align}
g_i^n &= (I - \epsilon \Delta t \Delta_h)^{-1} (m_i^n + \Delta t f_i^n), \\
g_i^* &= (I - \epsilon \Delta t \Delta_h)^{-1} (m_i^* + \Delta t f_i^*), \quad i=1,2,3
\end{align}
\begin{equation}
\begin{pmatrix}
m_1^* \\
m_2^* \\
m_3^*
\end{pmatrix}
=
\begin{pmatrix}
m_1^n + (g_2^n m_3^n - g_3^n m_2^n) \\
m_2^n + (g_3^n m_1^* - g_1^* m_3^n) \\
m_3^n + (g_1^* m_2^* - g_2^* m_1^*)
\end{pmatrix}.
\label{eq:step1_update}
\end{equation}

\paragraph{Step 2. Heat flow without constraints:}
\begin{equation}
\boldsymbol{f}^* = -Q(m_2^* \boldsymbol{e}_2 + m_3^* \boldsymbol{e}_3) + \boldsymbol{h}_s^* + \boldsymbol{h}_e,
\label{eq:heat_flow_f}
\end{equation}
\begin{equation}
\begin{pmatrix}
m_1^{**} \\
m_2^{**} \\
m_3^{**}
\end{pmatrix}
=
\begin{pmatrix}
m_1^* + \alpha \Delta t (\epsilon \Delta_h m_1^{**} + f_1^*) \\
m_2^* + \alpha \Delta t (\epsilon \Delta_h m_2^{**} + f_2^*) \\
m_3^* + \alpha \Delta t (\epsilon \Delta_h m_3^{**} + f_3^*)
\end{pmatrix}.
\label{eq:heat_flow_update}
\end{equation}

\paragraph{Step 3. Projection onto $S^2$:}
\begin{equation}
\begin{pmatrix}
m_1^{n+1} \\
m_2^{n+1} \\
m_3^{n+1}
\end{pmatrix}
=
\frac{1}{|m^{**}|}
\begin{pmatrix}
m_1^{**} \\
m_2^{**} \\
m_3^{**}
\end{pmatrix},
\label{eq:projection}
\end{equation}
where $\boldsymbol{m}^*$ denotes the intermediate values of $\boldsymbol{m}$. The stray field $\boldsymbol{h}_s$ is computed using the intermediate values $\boldsymbol{m}^*$ in \eqref{eq:step1_update} and \eqref{eq:heat_flow_f}.

\begin{figure}[htbp]
    \centering
    \subfloat[GSPM in 1D, $\alpha=0$]{\includegraphics[width=0.25\linewidth]{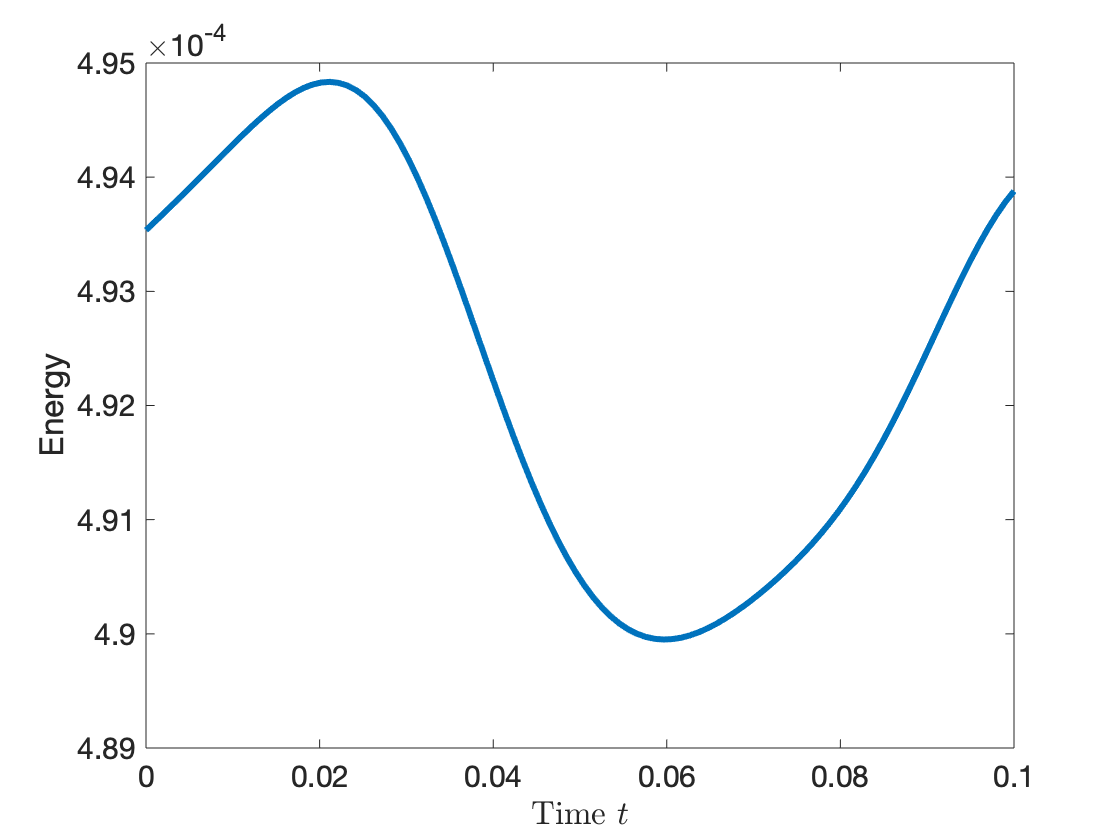}}
    \subfloat[GSPM in 1D, $\alpha=0.005$]{\includegraphics[width=0.25\linewidth]{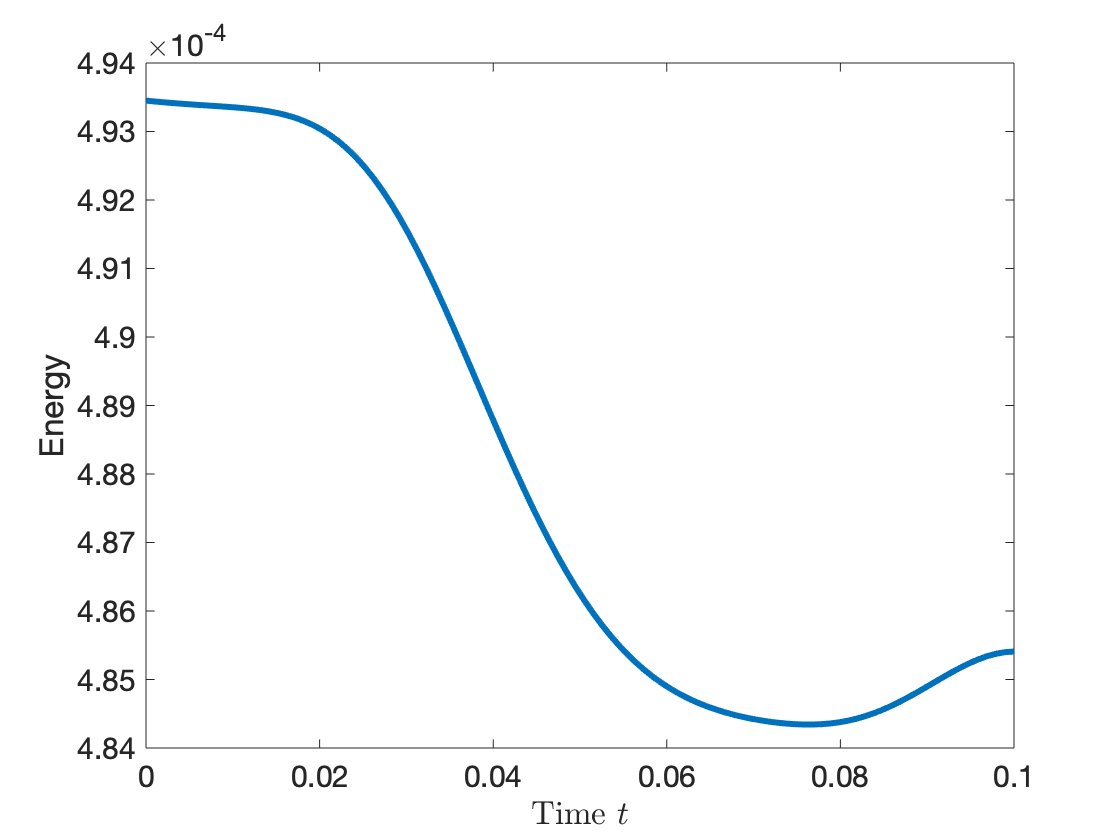}}
    \subfloat[GSPM in 1D, $\alpha=0.01$]{\includegraphics[width=0.25\linewidth]{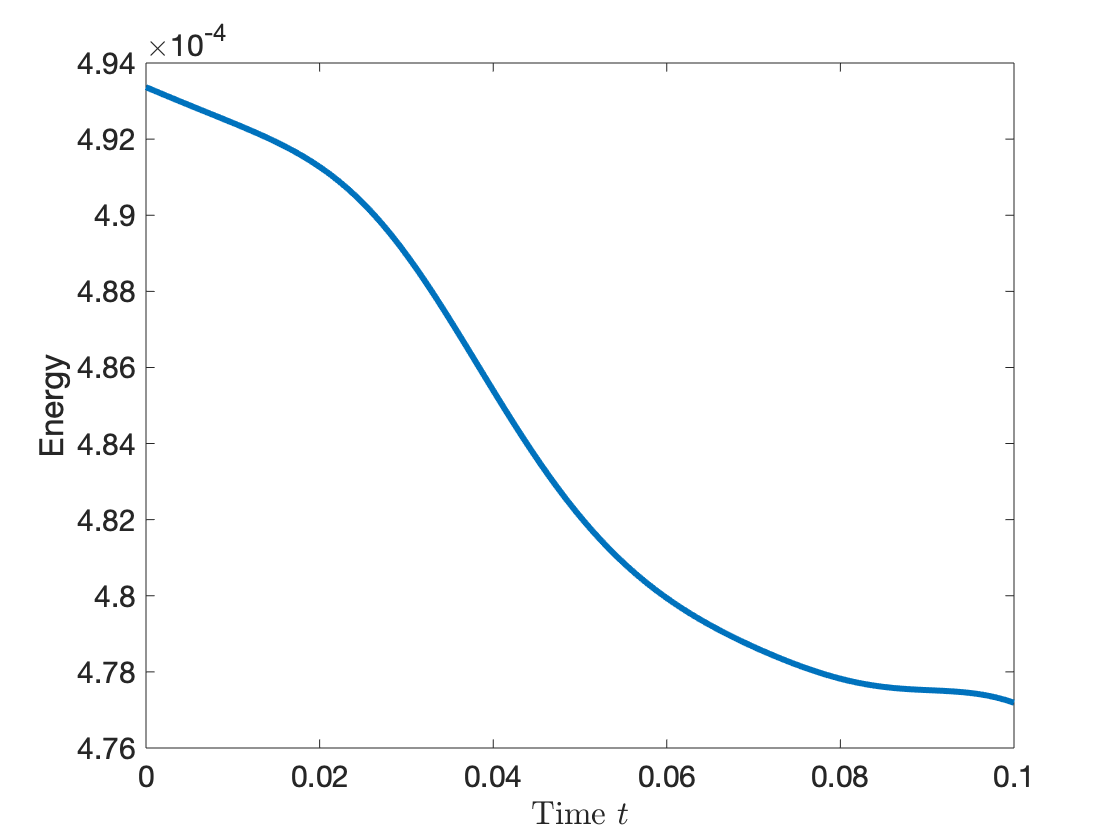}}
    \subfloat[GSPM in 1D, $\alpha=0.05$]{\includegraphics[width=0.25\linewidth]{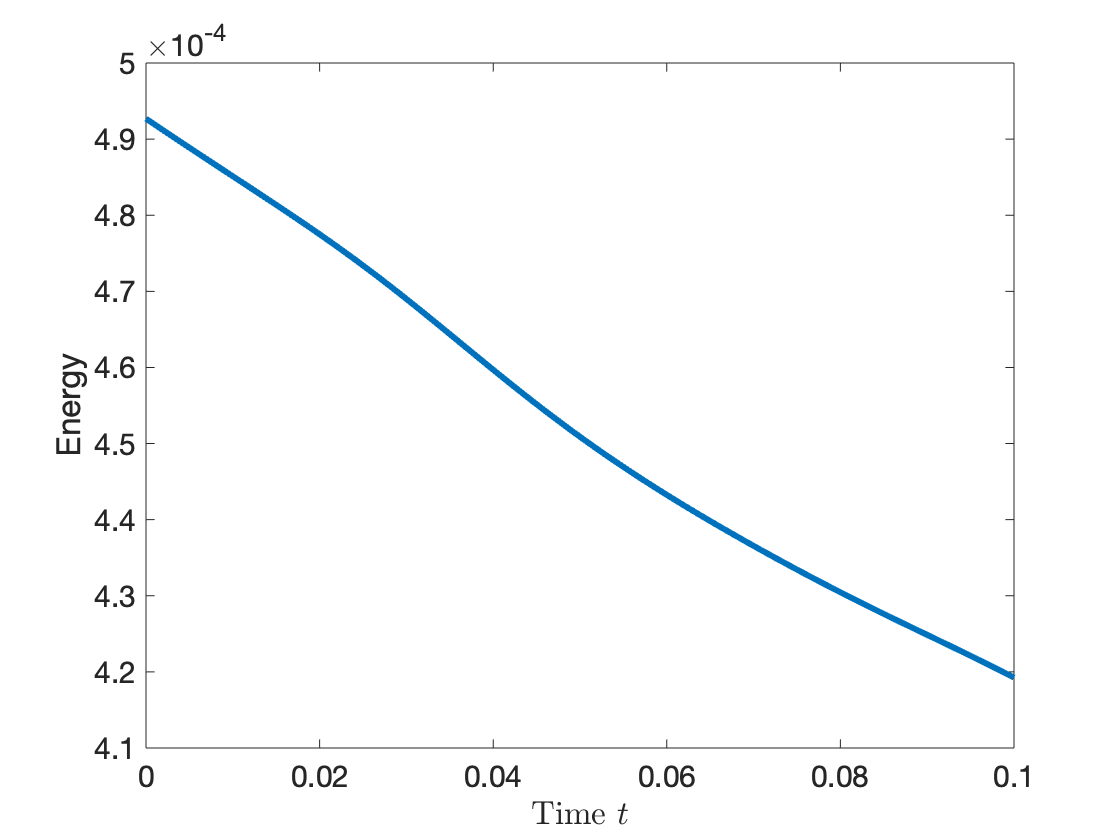}}
    \hspace{0.1in}
    \subfloat[GSPM in 3D, $\alpha=0$]{\includegraphics[width=0.25\linewidth]{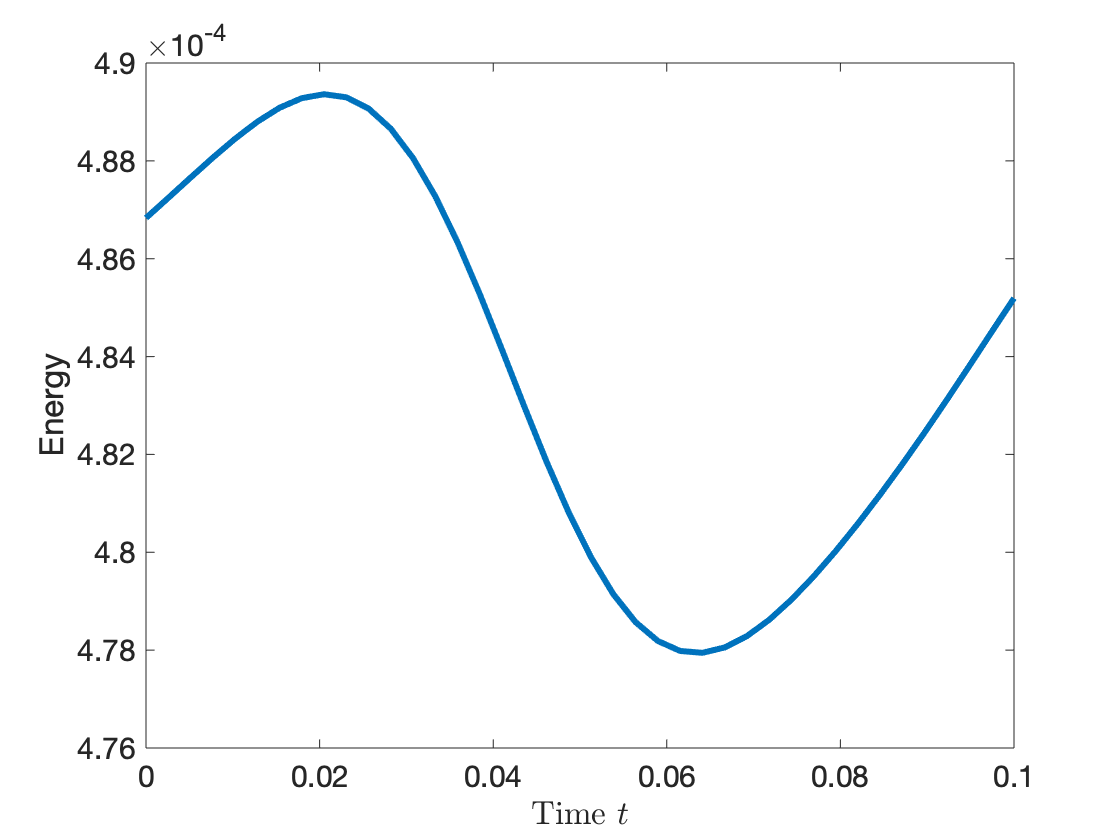}}
    \subfloat[GSPM in 3D, $\alpha=0.005$]{\includegraphics[width=0.25\linewidth]{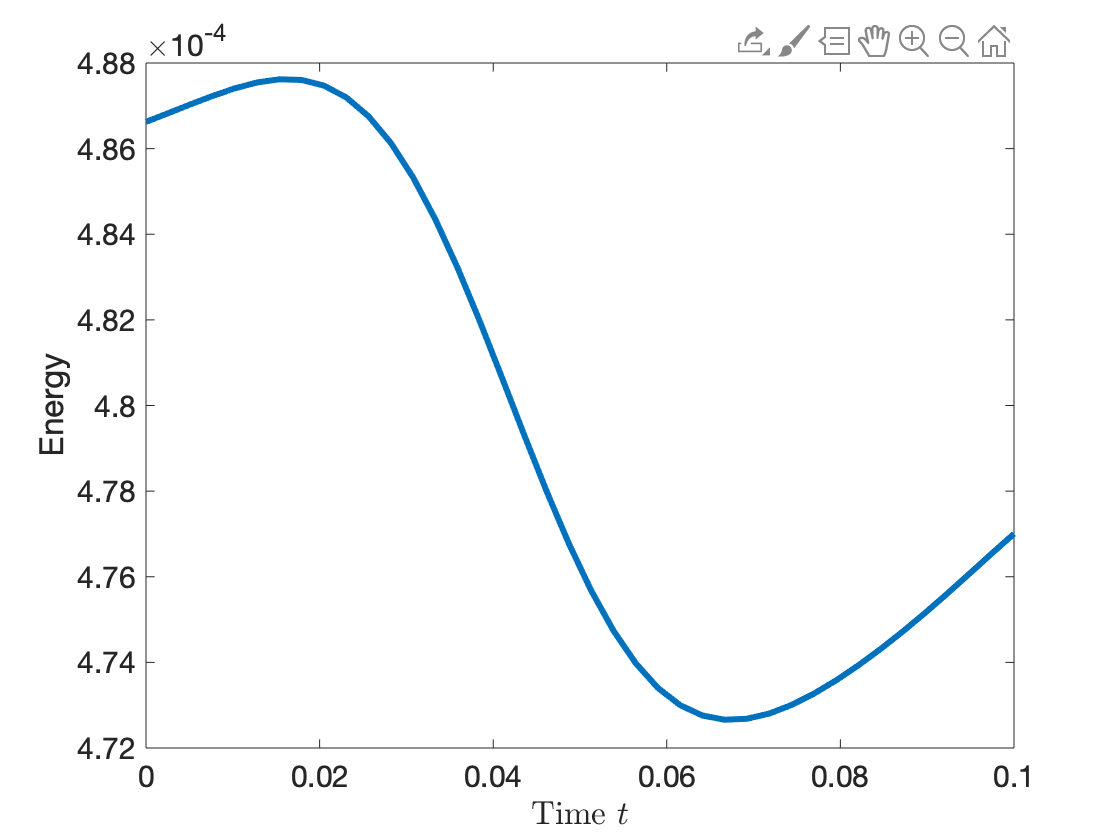}}
    \subfloat[GSPM in 3D, $\alpha=0.01$]{\includegraphics[width=0.25\linewidth]{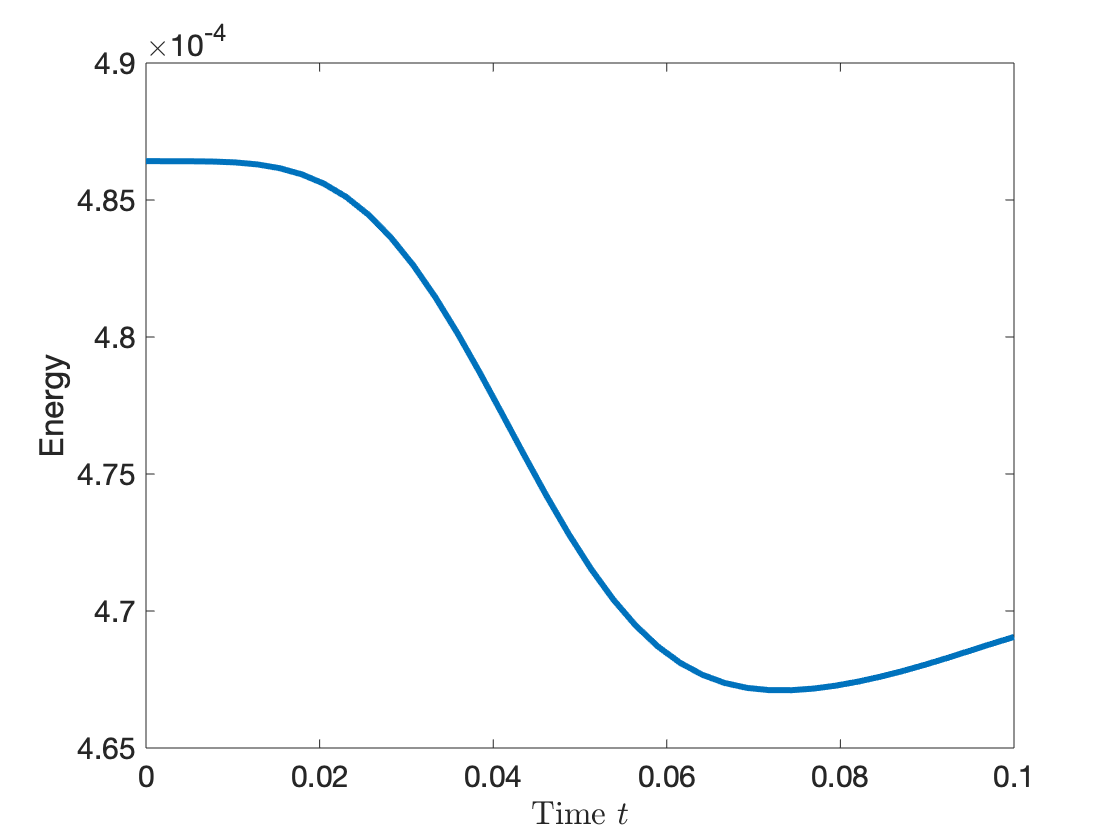}}
    \subfloat[GSPM in 3D, $\alpha=0.05$]{\includegraphics[width=0.25\linewidth]{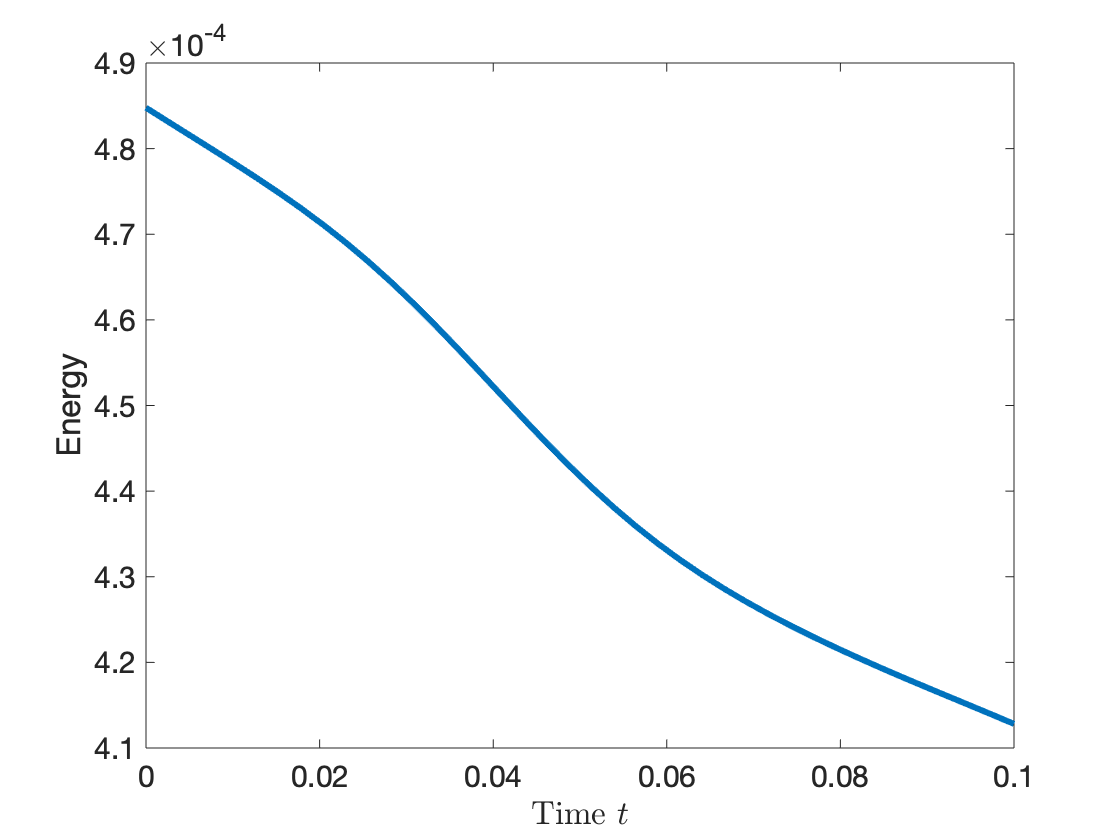}}
    \caption{The energy evolution between the GSPM. The initial condition is set to be $\m_0=\left(\cos(\cos(\pi x))\sin (0.01), \sin(\cos(\pi x))\sin (0.01), \cos (0.01)\right)^T$ for both 1D and 3D. The energy evolution with different damping parameters. From left to right panel with $\alpha=0,0.005,0.01,0.05$. The top row for 1D; The bottom for 3D.}
    \label{fig:energy-1}
\end{figure}

\begin{figure}[htbp]
    \centering
    \subfloat[GSPM in 1D]{\includegraphics[width=0.4\linewidth]{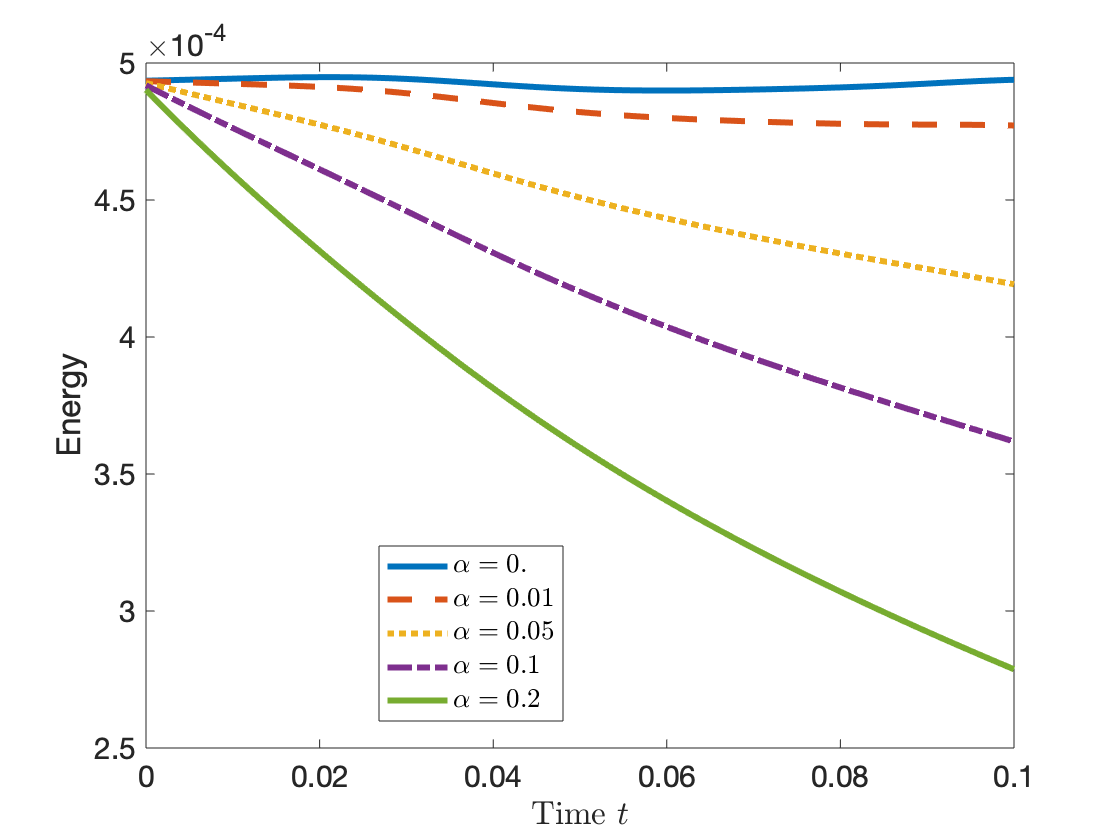}}
    \subfloat[GSPM in 3D]{\includegraphics[width=0.4\linewidth]{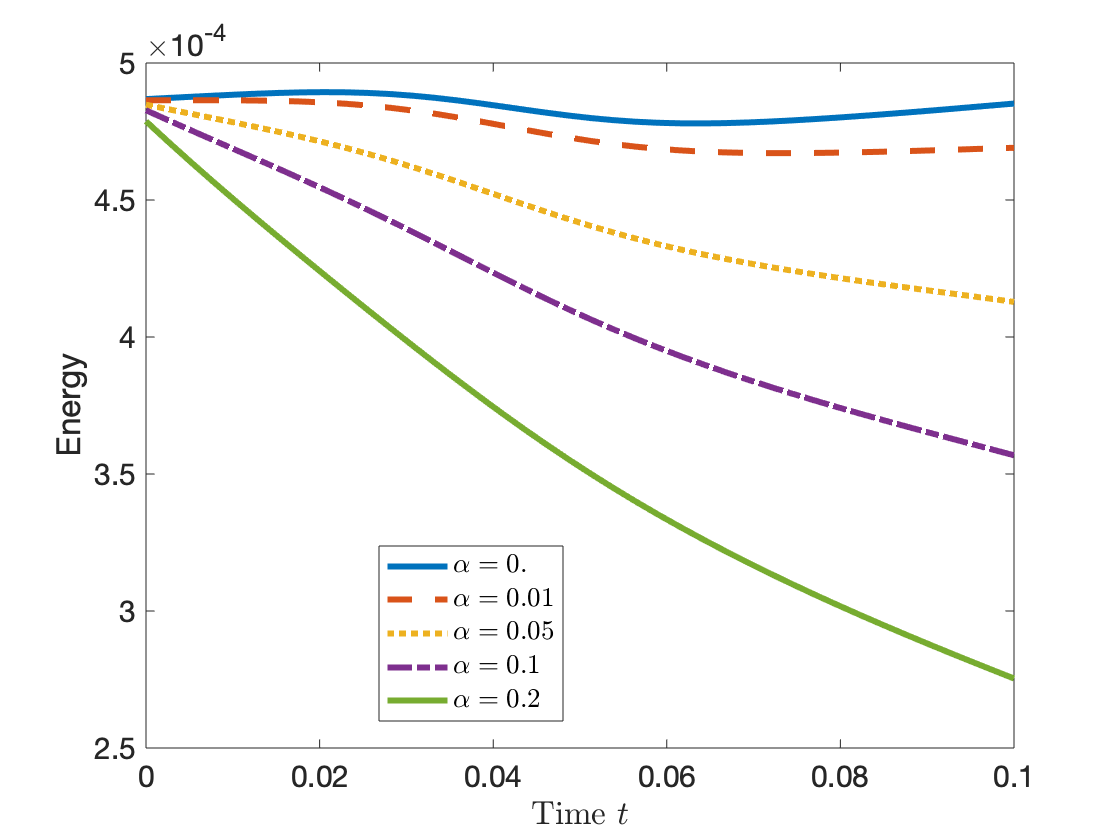}}
    \caption{The energy evolution of the GSPM. The initial condition is set to be $\m_0=\left(\cos(\cos(\pi x))\sin (0.01), \sin(\cos(\pi x))\sin (0.01), \cos (0.01)\right)^T$ for both 1D and 3D. The energy evolution with different damping parameters with $\alpha=0,0.01,0.05,0.1,0.2$.}
    \label{fig:energy-2}
\end{figure}

If we take the simple case with the exchange constant $\epsilon=1$, other fields $\f=0$, we take the initial condition as 
$\m_0=\left(\cos(\cos(\pi x))\sin (0.01), \sin(\cos(\pi x))\sin (0.01), \cos (0.01)\right)^T$. The results are shown in \Cref{fig:energy-1} and \Cref{fig:energy-2} which give the energy evolution of GSPM in 1D (top row) and 3D (bottom row) under four damping coefficients \(\alpha=0,0.005,0.01,0.05\) with the same initial magnetization. For undamped \(\alpha=0\), energy rises to a peak, drops to a minimum and bounces back. Larger \(\alpha\) weakens the early energy peak and suppresses the later rebound, leading to monotonic energy decay when \(\alpha=0.05\). The 1D and 3D energy curves share nearly identical qualitative trends with minor quantitative differences. 

\subsection{Proposed methods}

To improve its energy stability for GSPM with $\alpha=0$, we propose a symmetric Gauss-Seidel projection method (SGSPM), which is given as the following three steps:

\paragraph{Step 1. Implicit Symmetric Gauss-Seidel:}
\begin{align}
g_i^n &= (I - \epsilon \Delta t \Delta_h)^{-1} (m_i^n + \Delta t f_i^n), \\
g_i^* &= (I - \epsilon \Delta t \Delta_h)^{-1} (m_i^* + \Delta t f_i^*), \quad i=1,2,3
\end{align}
\begin{equation}
\begin{pmatrix}
m_1^* \\
m_2^* \\
m_3^*
\end{pmatrix}
=
\begin{pmatrix}
m_1^n + (g_2^n m_3^n - g_3^n m_2^n) \\
m_2^n + (g_3^n m_1^* - g_1^* m_3^n) \\
m_3^n + (g_1^* m_2^* - g_2^* m_1^*)
\end{pmatrix}.
\end{equation}
\begin{equation}
\begin{pmatrix}
m_3^{**}\\
m_2^{**} \\
m_1^{**} 
\end{pmatrix}
=
\begin{pmatrix}
m_3^n + (g_1^{*} m_2^{*} - g_2^{*} m_1^{*})\\
m_2^n + (g_3^{**} m_1^{*} - g_1^{*} m_3^{**}) \\
m_1^n + (g_2^{**} m_3^{**} - g_3^{**} m_2^{**}) 
\end{pmatrix}.
\end{equation}

\paragraph{Step 2. Heat flow without constraints:}
\begin{equation}
\boldsymbol{f}^{**} = -Q(m_2^{**} \boldsymbol{e}_2 + m_3^{**} \boldsymbol{e}_3) + \boldsymbol{h}_s^{**} + \boldsymbol{h}_e,
\end{equation}
\begin{equation}
\begin{pmatrix}
m_1^{***} \\
m_2^{***} \\
m_3^{***}
\end{pmatrix}
=
\begin{pmatrix}
m_1^{**} + \alpha \Delta t (\epsilon \Delta_h m_1^{***} + f_1^{**}) \\
m_2^{**} + \alpha \Delta t (\epsilon \Delta_h m_2^{***} + f_2^{**}) \\
m_3^{**} + \alpha \Delta t (\epsilon \Delta_h m_3^{***} + f_3^{**})
\end{pmatrix}.
\end{equation}

\paragraph{Step 3. Projection onto $S^2$:}
\begin{equation}
\begin{pmatrix}
m_1^{n+1} \\
m_2^{n+1} \\
m_3^{n+1}
\end{pmatrix}
=
\frac{1}{|m^{***}|}
\begin{pmatrix}
m_1^{***} \\
m_2^{***} \\
m_3^{***}
\end{pmatrix}.
\end{equation}

\begin{figure}[htbp]
    \centering
    \subfloat[SGSPM in 1D, $\alpha=0$]{\includegraphics[width=0.25\linewidth]{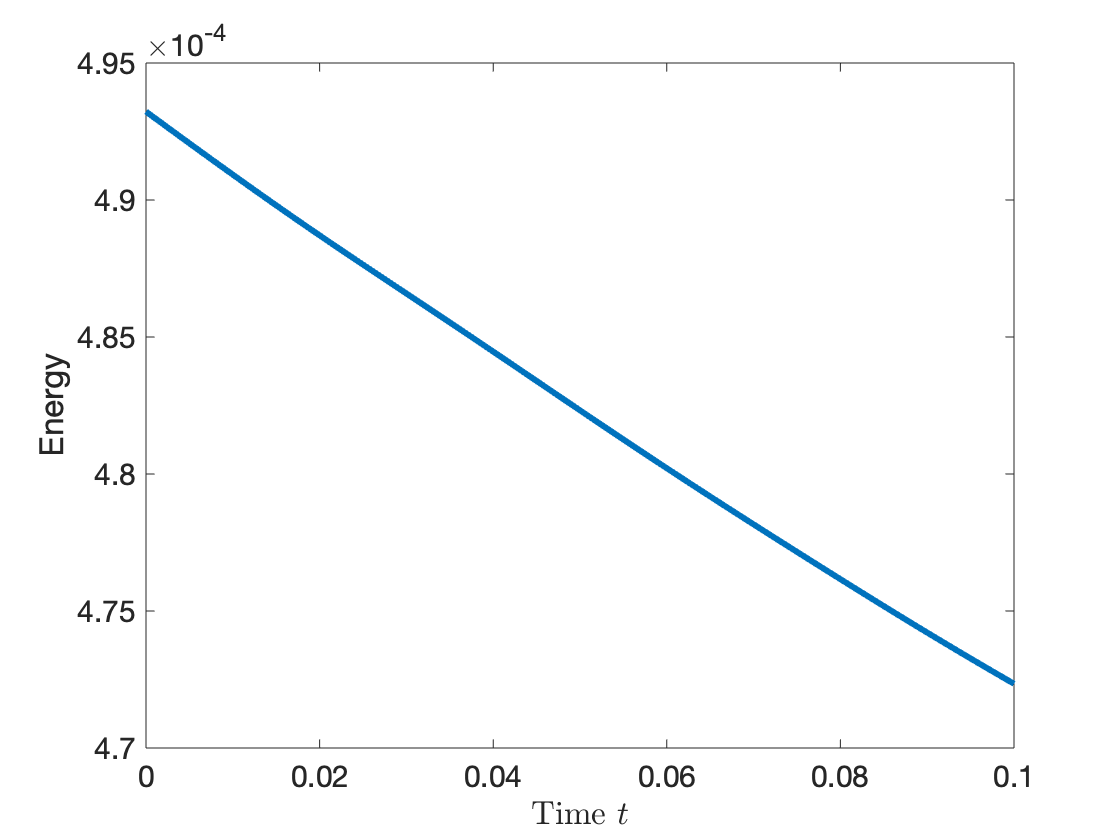}}
    \subfloat[SGSPM in 1D, $\alpha=0.005$]{\includegraphics[width=0.25\linewidth]{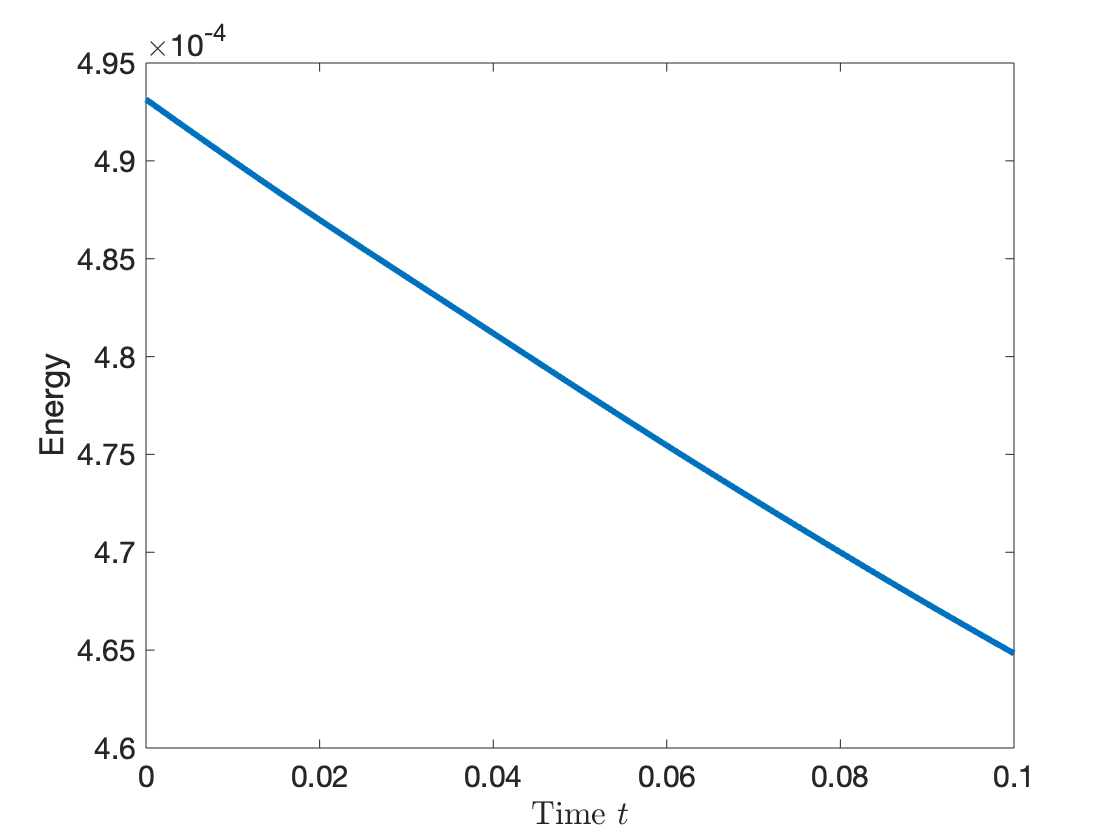}}
    \subfloat[SGSPM in 1D, $\alpha=0.01$]{\includegraphics[width=0.25\linewidth]{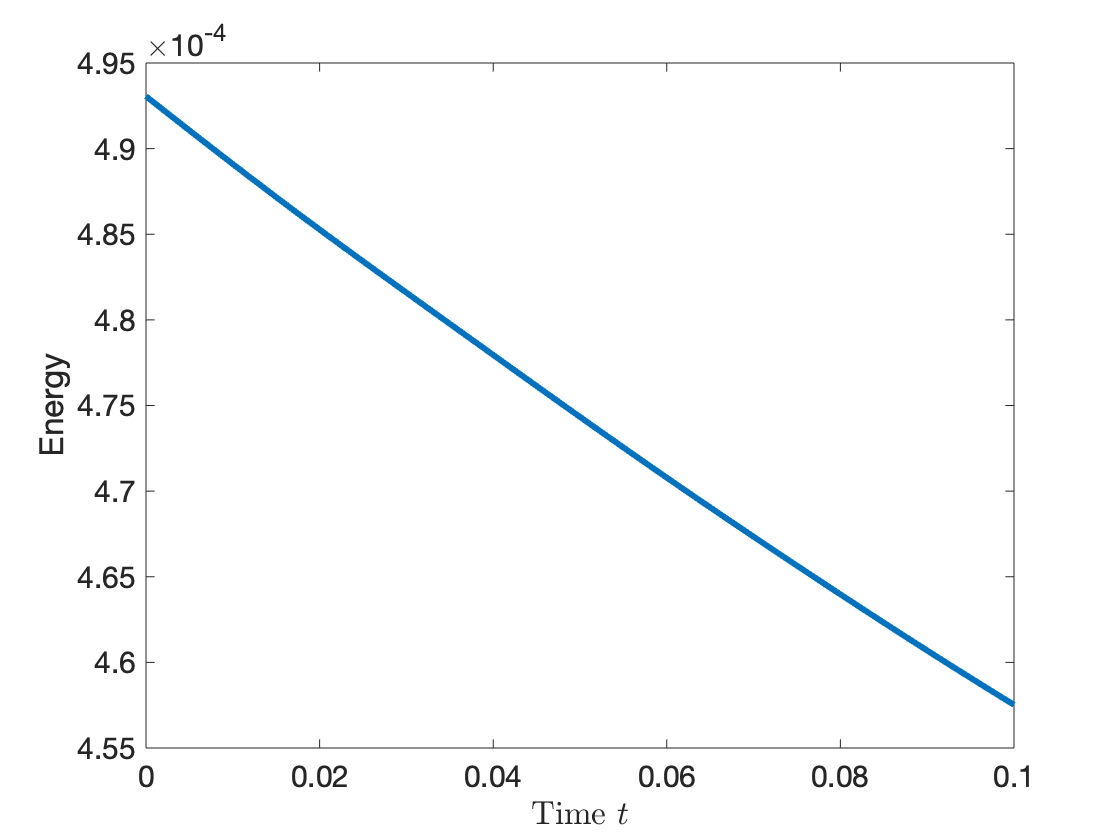}}
    \subfloat[SGSPM in 1D, $\alpha=0.05$]{\includegraphics[width=0.25\linewidth]{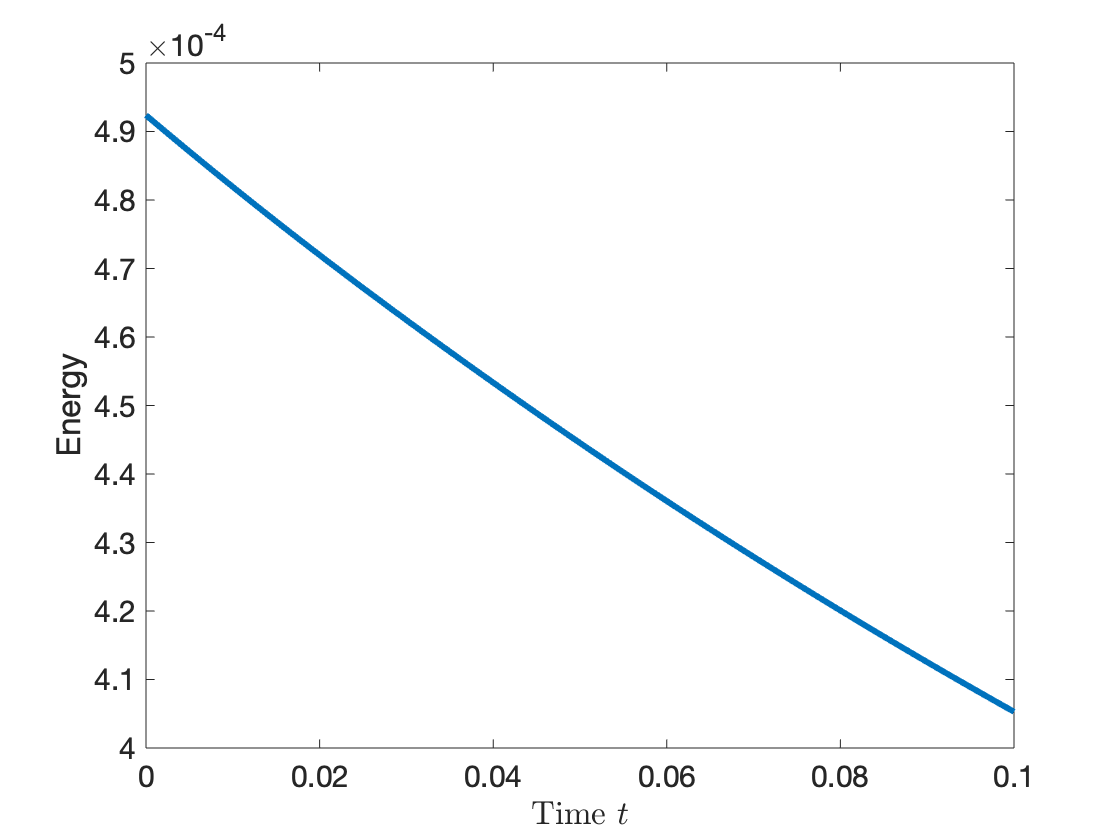}}
    \hspace{0.1in}
    \subfloat[SGSPM in 3D, $\alpha=0$]{\includegraphics[width=0.25\linewidth]{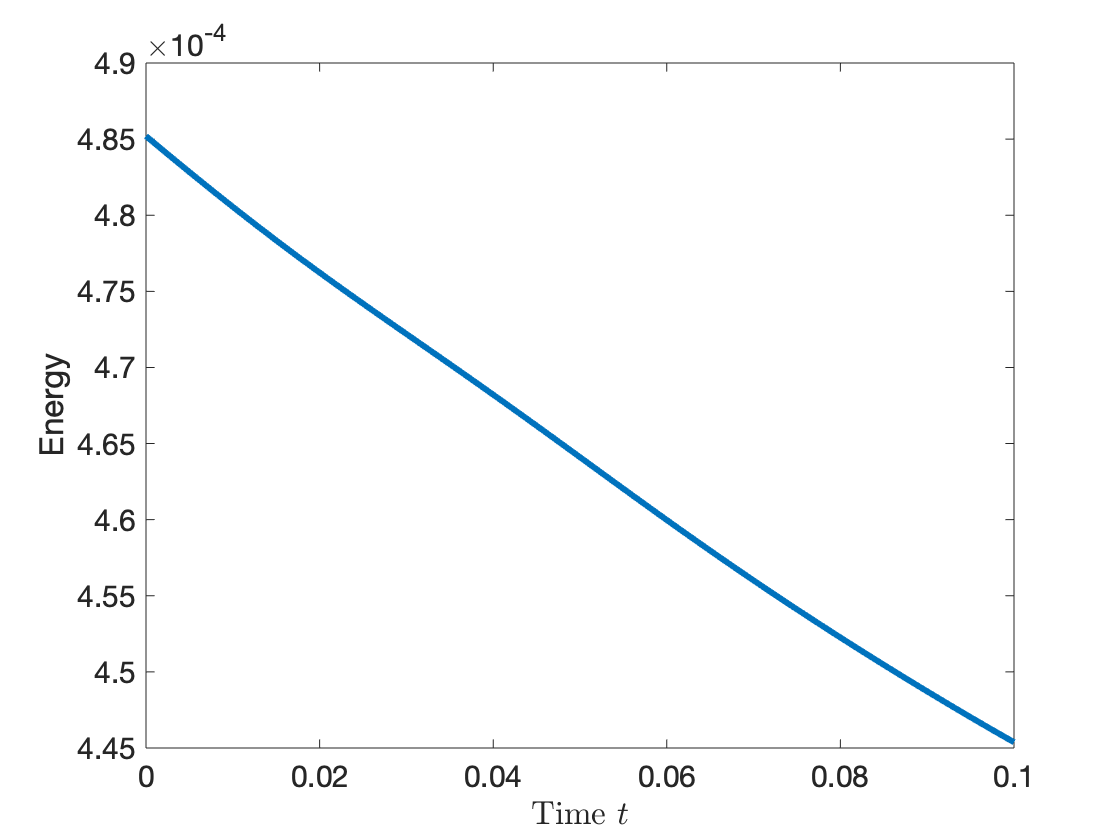}}
    \subfloat[SGSPM in 3D, $\alpha=0.005$]{\includegraphics[width=0.25\linewidth]{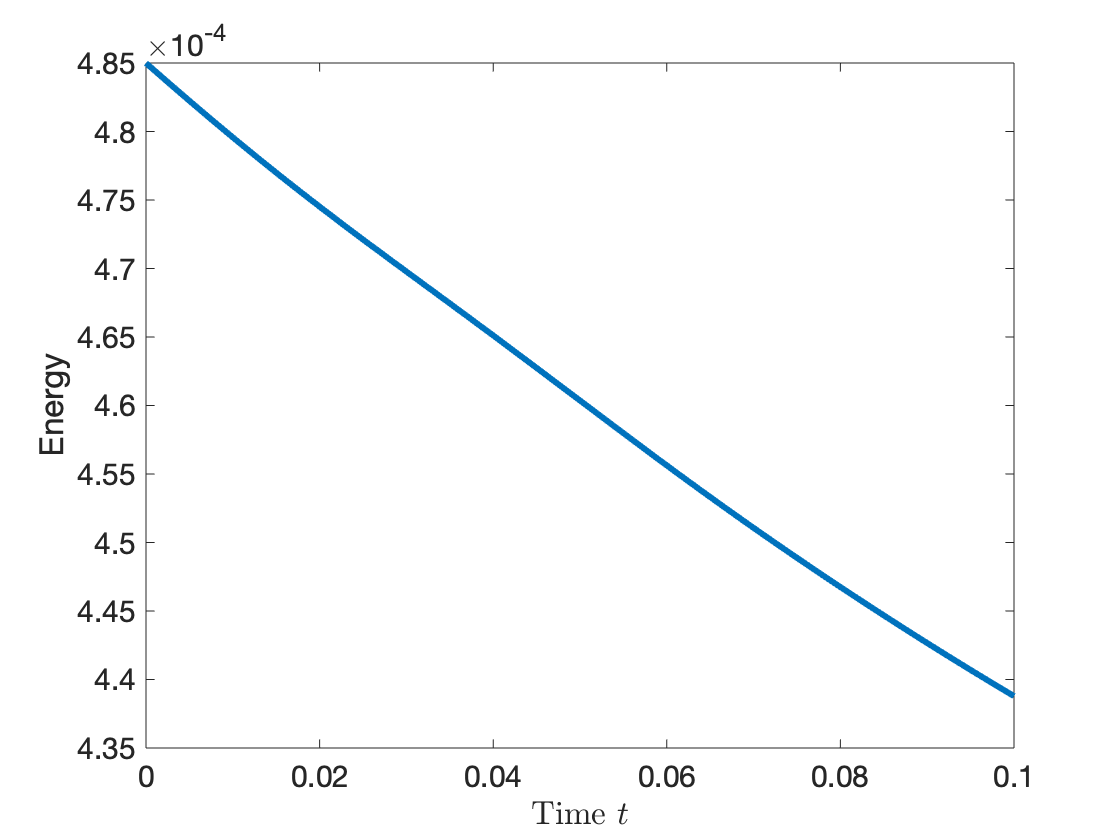}}
    \subfloat[SGSPM in 3D, $\alpha=0.01$]{\includegraphics[width=0.25\linewidth]{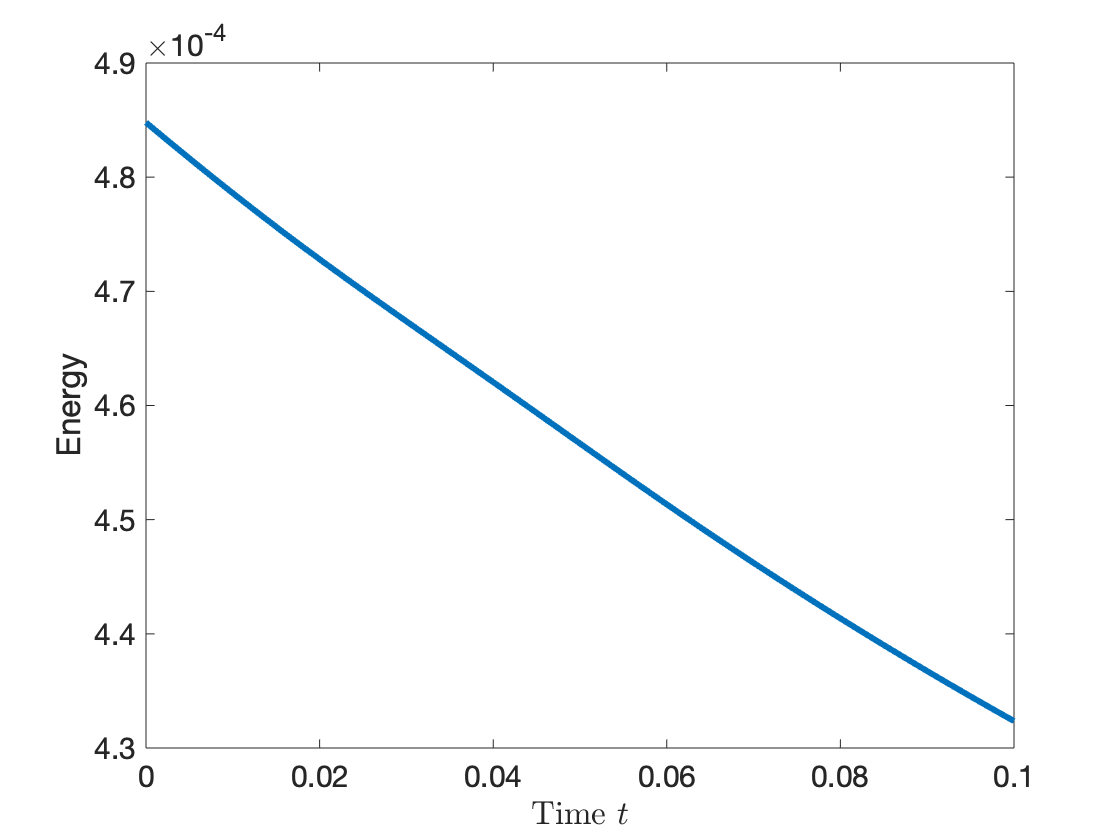}}
    \subfloat[SGSPM in 3D, $\alpha=0.05$]{\includegraphics[width=0.25\linewidth]{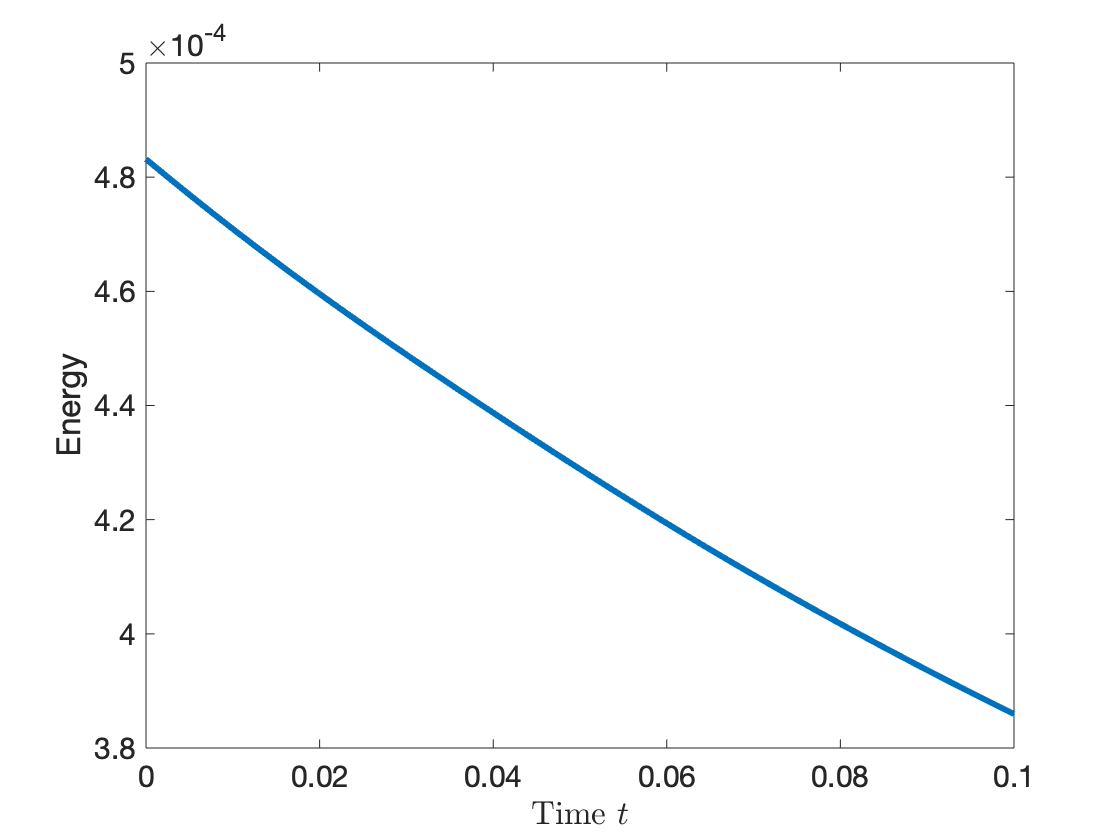}}
    \caption{The energy evolution between the proposed method (SGSPM). The initial condition is set to be $\m_0=\left(\cos(\cos(\pi x))\sin (0.01), \sin(\cos(\pi x))\sin (0.01), \cos (0.01)\right)^T$ for both 1D and 3D. The energy evolution with different damping parameters. From left to right panel with $\alpha=0,0.005,0.01,0.05$. The top row for 1D; The bottom for 3D.}
    \label{fig:energy-3}
\end{figure}

\begin{figure}[htbp]
    \centering
    \subfloat[SGSPM in 1D]{\includegraphics[width=0.4\linewidth]{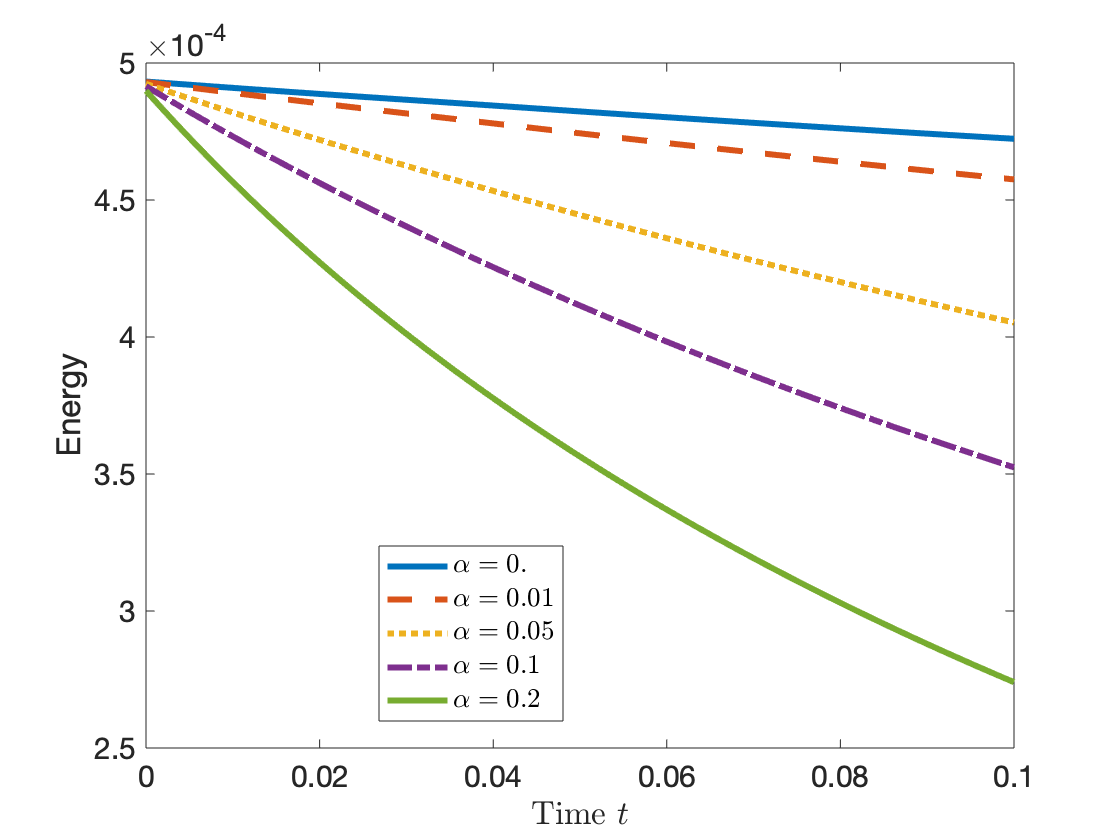}}
    \subfloat[SGSPM in 3D]{\includegraphics[width=0.4\linewidth]{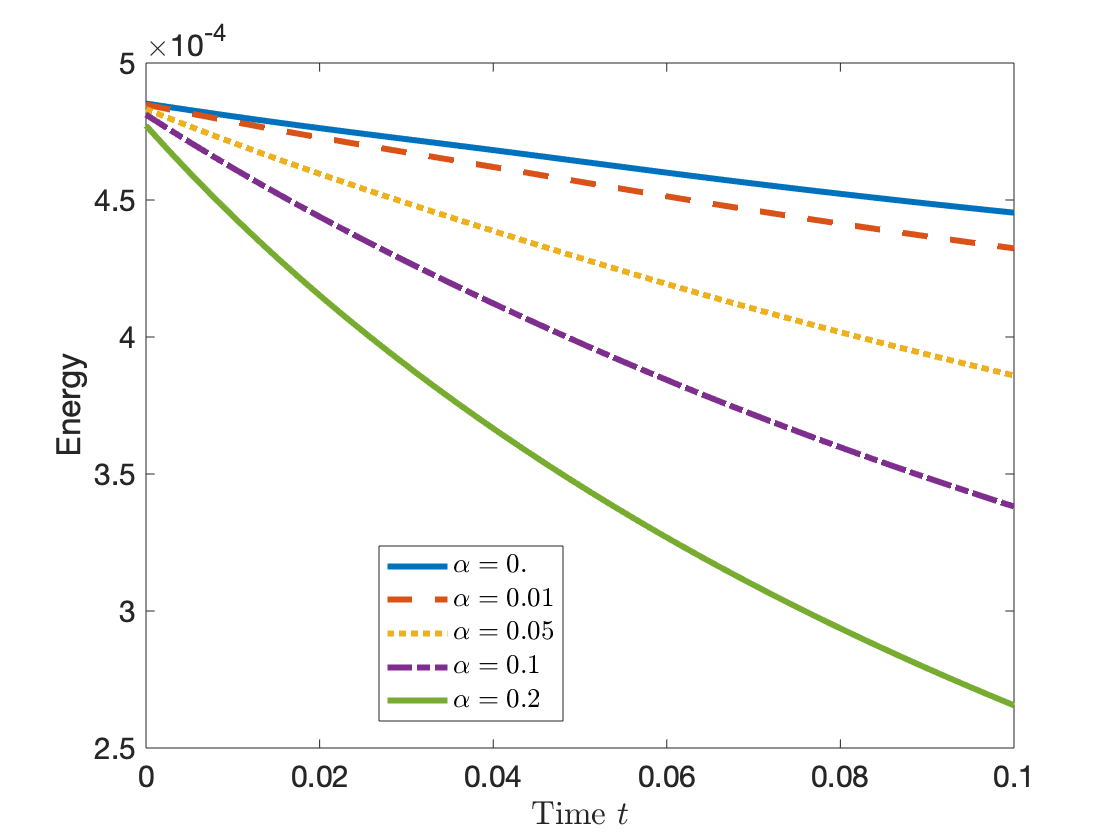}}
    \caption{The energy evolution of the proposed method (SGSPM). The initial condition is set to be $\m_0=\left(\cos(\cos(\pi x))\sin (0.01), \sin(\cos(\pi x))\sin (0.01), \cos (0.01)\right)^T$ for both 1D and 3D. The energy evolution with different damping parameters with $\alpha=0,0.01,0.05,0.1,0.2$.}
    \label{fig:energy-4}
\end{figure}

If we take the simple case with the exchange constant $\epsilon=1$, other fields $\f=0$, we take the initial condition as 
$\m_0=\left(\cos(\cos(\pi x))\sin (0.01), \sin(\cos(\pi x))\sin (0.01), \cos (0.01)\right)^T$. The results for SGSPM are shown in \Cref{fig:energy-3} and \Cref{fig:energy-4}, which display energy curves of the proposed SGSPM for 1D (top) and 3D (bottom) cases under damping coefficients \(\alpha=0,0.005,0.01,0.05\) with identical initial magnetization. Unlike GSPM, SGSPM yields strictly monotonic energy decay over time for all \(\alpha\) values, with faster energy dissipation observed as \(\alpha\) increases. The 1D and 3D configurations show consistent qualitative energy evolution trends.

\section{Theoretical analysis}\label{sec:theory}


For simplicity, we take the exchange constant $\epsilon=1$, other fields $\f=0$, we analyze the step 1, for our model $\partial_t \m=-\m\times \Delta \m$, the GS method is given by

\begin{align*}
    \frac{m_1^{*}-m_1^n}{\Delta t}&=-m_2^n\cdot \Delta_h g_3^n+m_3^n\cdot\Delta_h g_2^n\\
    \frac{m_2^{*}-m_2^n}{\Delta t}&=-m_3^n\cdot \Delta_h g_1^{*}+m_1^*\cdot\Delta_h g_3^n\\
    \frac{m_3^{*}-m_3^n}{\Delta t}&=-m_1^* \cdot\Delta_h g_2^*+m_2^*\cdot\Delta_h g_1^*.
\end{align*}
The SGS method is given by
forwardly
\begin{align*}
    \frac{m_1^{*}-m_1^n}{\Delta t}&=-m_2^n\cdot \Delta_h g_3^n+m_3^n\cdot\Delta_h g_2^n\\
    \frac{m_2^{*}-m_2^n}{\Delta t}&=-m_3^n\cdot \Delta_h g_1^{*}+m_1^*\cdot\Delta_h g_3^n\\
    \frac{m_3^{*}-m_3^n}{\Delta t}&=-m_1^* \cdot\Delta_h g_2^*+m_2^*\cdot\Delta_h g_1^*,
\end{align*}
and backwardly
\begin{align*}
\frac{m_3^{**}-m_3^n}{\Delta t}&=-m_1^* \cdot\Delta_h g_2^*+m_2^*\cdot\Delta_h g_1^*\\
 \frac{m_2^{**}-m_2^n}{\Delta t}&=-m_3^{**}\cdot \Delta_h g_1^{*}+m_1^*\cdot\Delta_h g_3^{**}\\
    \frac{m_1^{**}-m_1^n}{\Delta t}&=-m_2^{**}\cdot \Delta_h g_3^{**}+m_3^{**}\cdot\Delta_h g_2^{**},
\end{align*}
where $g_i^s=(I-\Delta t\Delta_h)^{-1} m_i^s$, $s=n,*,**$ which is from the heat diffusion to obtain a stabilized $m_i$.

The GS method in a matrix form is given by
\begin{align}\label{GS-matrix}
\begin{pmatrix}
I & 0 & 0 \\
-\Delta t \Delta_h g_3^n & I & 0 \\
\Delta t \Delta_h g_2^* & -\Delta t \Delta_h g_1^* & I
\end{pmatrix}
\begin{pmatrix}m_1^* \\ m_2^* \\ m_3^*\end{pmatrix}
=\begin{pmatrix}
0 & -\Delta t \Delta_h g_3^n & \Delta t \Delta_h g_2^n \\
0 & 0 & -\Delta t \Delta_h g_1^*  \\
0& 0 & 0
\end{pmatrix}
\begin{pmatrix}m_1^n \\ m_2^n \\ m_3^n\end{pmatrix}+\begin{pmatrix}m_1^n \\ m_2^n \\ m_3^n\end{pmatrix}.
\end{align}

Originally, we solve to $A\m^{*}=\m^n$
where
\begin{align}\label{eq-M}
    A=\begin{pmatrix}
I & \Delta t \Delta_h g_3(m_3) & -\Delta t \Delta_h g_2(m_2) \\
-\Delta t \Delta_h g_3(m_3) & I & \Delta t \Delta_h g_1(m_1) \\
\Delta t \Delta_h g_2(m_2) & -\Delta t \Delta_h g_1(m_1) & I
\end{pmatrix}.
\end{align}
Gauss-Seidel iteration gives
\begin{align*}
    A=\begin{pmatrix}
I & \Delta t \Delta_h g_3^n & -\Delta t \Delta_h g_2^n \\
-\Delta t \Delta_h g_3^n & I & \Delta t \Delta_h g_1^{*} \\
\Delta t \Delta_h g_2^* & -\Delta t \Delta_h g_1^* & I
\end{pmatrix},
\end{align*}
and
\begin{align*}
    D-L=\begin{pmatrix}
I & 0 & 0 \\
-\Delta t \Delta_h g_3^n & I & 0 \\
\Delta t \Delta_h g_2^* & -\Delta t \Delta_h g_1^* & I
\end{pmatrix},\quad U=\begin{pmatrix}
0 & -\Delta t \Delta_h g_3^n & \Delta t \Delta_h g_2^n \\
0 & 0 & -\Delta t \Delta_h g_1^*  \\
0& 0 & 0
\end{pmatrix},
\end{align*}
where $D$ is the diagonal matrix and $L$ is the lower trigonametric matrix, $U$ is the upper trigonametric matrix. Therfore, the GS iteration scheme is
\begin{align*}
    \m^{*}=(D-L)^{-1}U\m^n+(D-L)^{-1}\m^n.
\end{align*}
The iteration matrix for Gauss-Seidel is given by
\begin{align*}
    G_{GS}=(D-L)^{-1}U.
\end{align*}

The SGS method in a matrix form is given by
\begin{align}\label{SGS-matrix}
\begin{aligned}
    \begin{pmatrix}
I & 0 & 0 \\
-\Delta t \Delta_h g_3^n & I & 0 \\
\Delta t \Delta_h g_2^* & -\Delta t \Delta_h g_1^* & I
\end{pmatrix}
\begin{pmatrix}m_1^* \\ m_2^* \\ m_3^*\end{pmatrix}
=\begin{pmatrix}
0 & -\Delta t \Delta_h g_3^n & \Delta t \Delta_h g_2^n \\
0 & 0 & -\Delta t \Delta_h g_1^*  \\
0& 0 & 0
\end{pmatrix}
\begin{pmatrix}m_1^n \\ m_2^n \\ m_3^n\end{pmatrix}+\begin{pmatrix}m_1^n \\ m_2^n \\ m_3^n\end{pmatrix}\\
\begin{pmatrix}
I & -\Delta t \Delta_h g_3^{**} & \Delta t \Delta_h g_2^{**} \\
0 & I & -\Delta t \Delta_h g_1^{*} \\
0 & 0 & I
\end{pmatrix}
\begin{pmatrix}m_1^{**} \\ m_2^{**} \\ m_3^{**}\end{pmatrix}
=\begin{pmatrix}
0 & 0 & 0 \\
\Delta t \Delta_h g_3^{**} & 0 & 0  \\
-\Delta t \Delta_h g_2^{*} & \Delta t \Delta_h g_1^{*}  & 0
\end{pmatrix}
\begin{pmatrix}m_1^* \\ m_2^* \\ m_3^*\end{pmatrix}+\begin{pmatrix}m_1^n \\ m_2^n \\ m_3^n\end{pmatrix}.
\end{aligned}
\end{align}

The forward GS of SGS gives
\begin{align*}
    \m^{*}=(D-L)^{-1}U\m^n+(D-L)^{-1}\m^n,
\end{align*}
where
\begin{align*}
    D-L=\begin{pmatrix}
I & 0 & 0 \\
-\Delta t \Delta_h g_3^n & I & 0 \\
\Delta t \Delta_h g_2^* & -\Delta t \Delta_h g_1^* & I
\end{pmatrix},\quad U=\begin{pmatrix}
0 & -\Delta t \Delta_h g_3^n & \Delta t \Delta_h g_2^n \\
0 & 0 & -\Delta t \Delta_h g_1^*  \\
0& 0 & 0
\end{pmatrix}.
\end{align*}
The backward GS of SGS gives
\begin{align*}
    \m^{**}=(D-U)^{-1}L\m^{*}+(D-U)^{-1}\m^n,
\end{align*}
where
\begin{align*}
    D-U=\begin{pmatrix}
I & -\Delta t \Delta_h g_3^{**} & \Delta t \Delta_h g_2^{**} \\
0 & I & -\Delta t \Delta_h g_1^{*} \\
0 & 0 & I
\end{pmatrix},\quad L=\begin{pmatrix}
0 & 0 & 0 \\
\Delta t \Delta_h g_3^{**} & 0 & 0  \\
-\Delta t \Delta_h g_2^{*} & \Delta t \Delta_h g_1^{*}  & 0
\end{pmatrix}.
\end{align*}
Therefore, for SGS iteration, we have
\begin{align*}
     \m^{**}&=(D-U)^{-1}L[(D-L)^{-1}U\m^n+(D-L)^{-1}\m^n]+(D-U)^{-1}\m^n\\
     &=(D-U)^{-1}L(D-L)^{-1}U\m^n+(D-U)^{-1}(L(D-L)^{-1}+I)\m^n.
\end{align*}
The SGS iteration matrix is
\begin{align*}
    G_{SGS}=(D-U)^{-1}L(D-L)^{-1}U.
\end{align*}
For the iteration matrix for GS and SGS method, we have
\begin{align*}
 G_{GS}&=(D-L)^{-1}U\\
     G_{SGS}&=(D-U)^{-1}L(D-L)^{-1}U.
\end{align*}
 The spectral radius of the matrix $A$ is defined by
\begin{align*}
    \rho(A)=\max\{|\lambda|,\; \lambda \text{ is the eigenvalue of } A\}.
\end{align*}
Numerically, we compute the spectral radius for $\rho(G_{GS})$ and $\rho(G_{SGS})$.

\begin{figure}[htbp]
    \centering
    \subfloat[$\alpha=0.001$, $N_x=200$, $Nt=10$]{\includegraphics[width=0.35\linewidth]{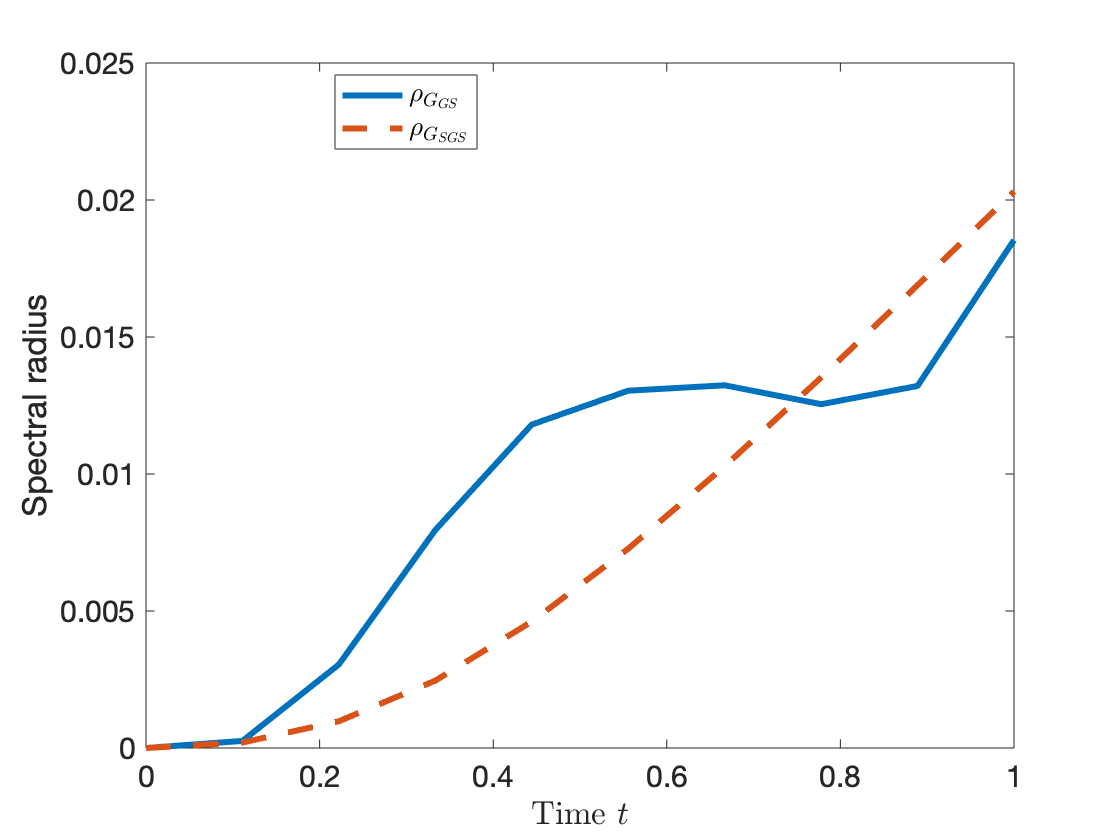}}
     \subfloat[$\alpha=0.001$, $N_x=200$, $Nt=100$]{\includegraphics[width=0.35\linewidth]{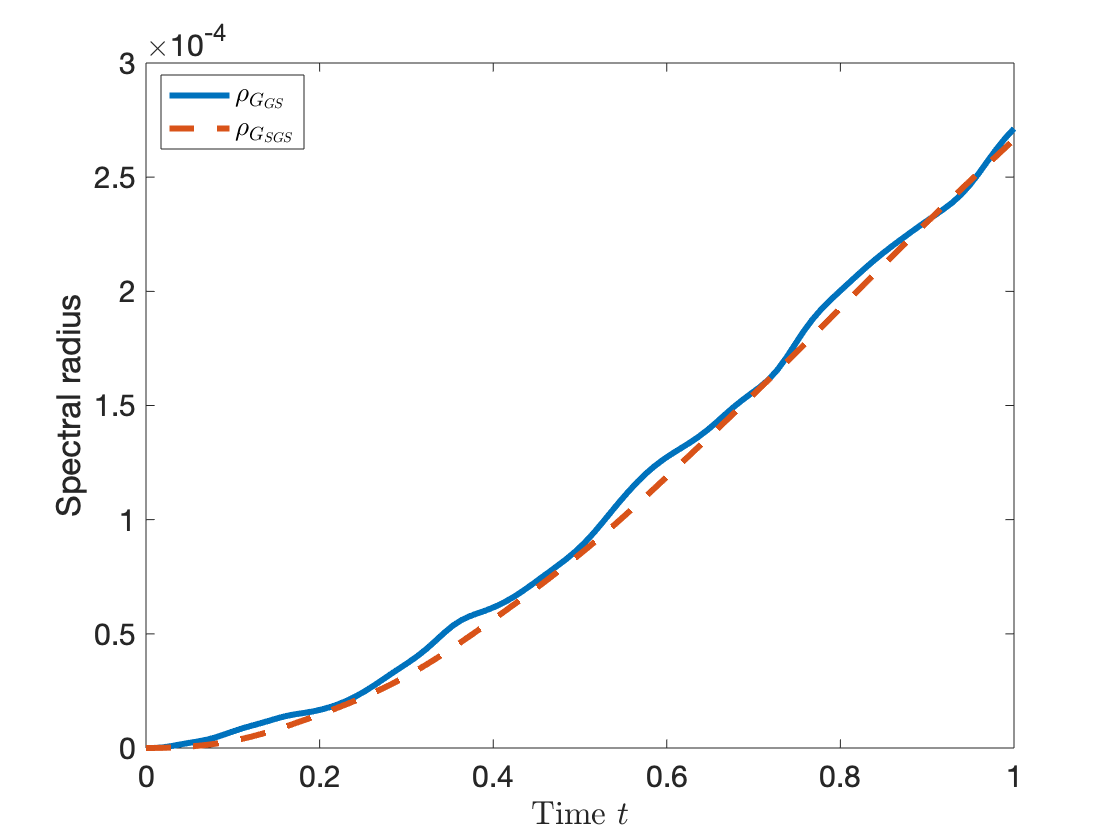}}
      \subfloat[$\alpha=0.001$, $N_x=200$, $Nt=1000$]{\includegraphics[width=0.35\linewidth]{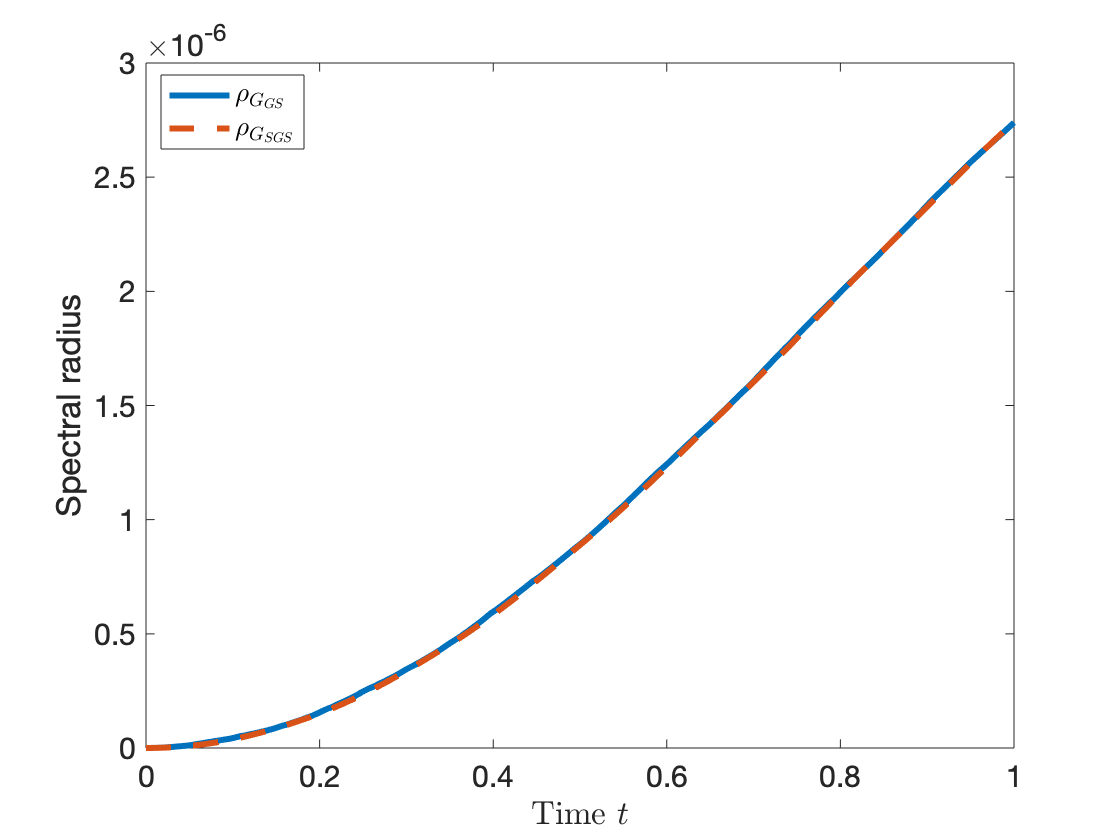}}
      \hspace{0.1in}
      \subfloat[$\alpha=0.001$, $N_x=20$, $Nt=1000$]{\includegraphics[width=0.35\linewidth]{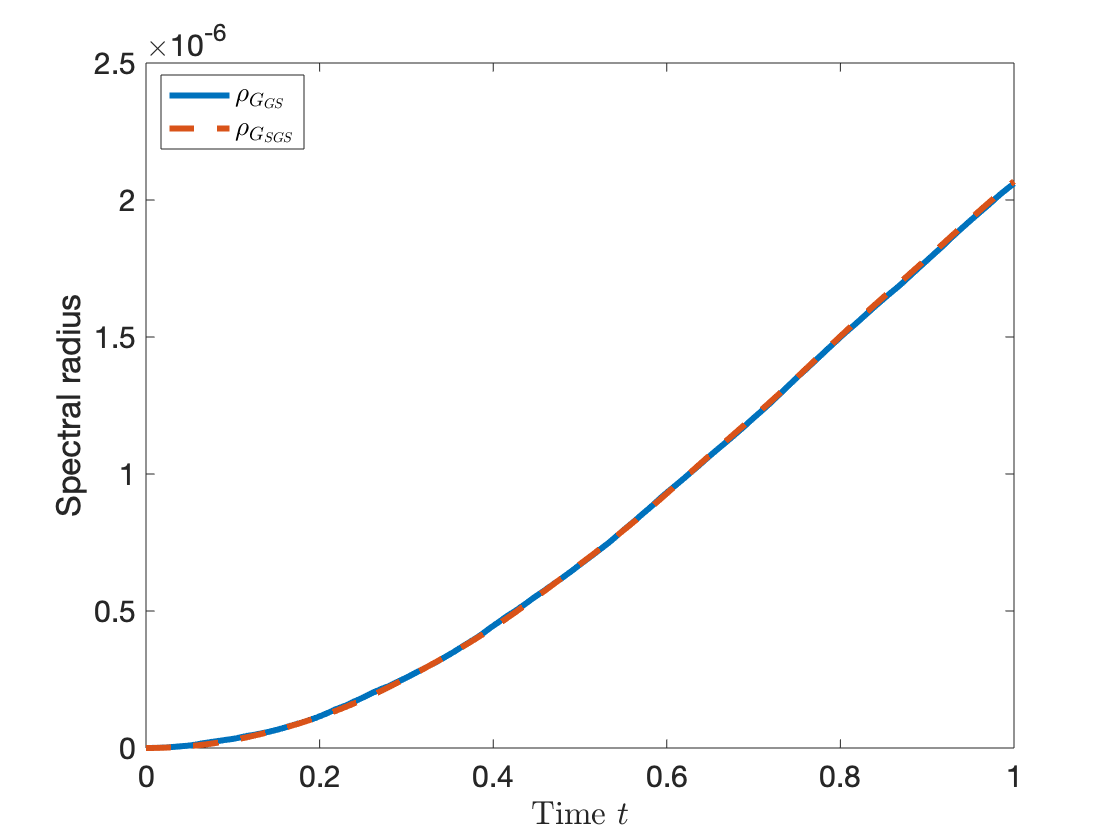}}
      \subfloat[$\alpha=0.001$, $N_x=400$, $Nt=100$]{\includegraphics[width=0.35\linewidth]{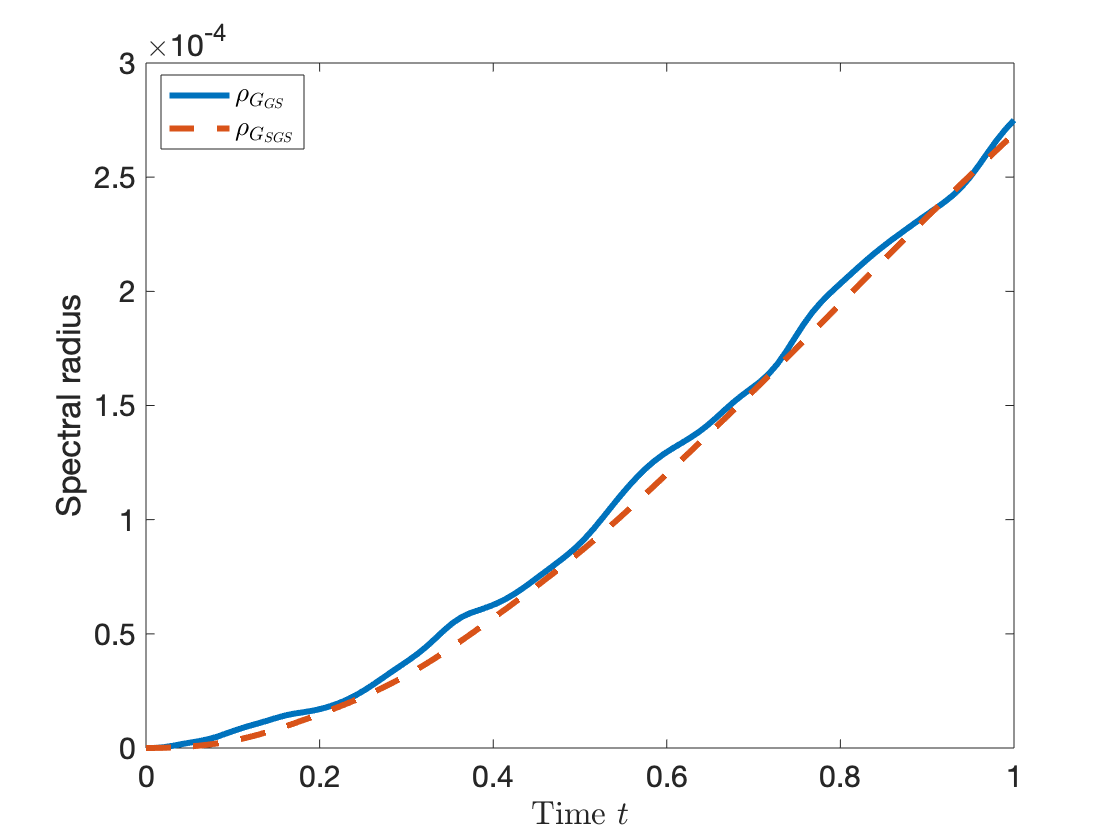}}
      \subfloat[$\alpha=0$, $N_x=1000$, $Nt=10$]{\includegraphics[width=0.35\linewidth]{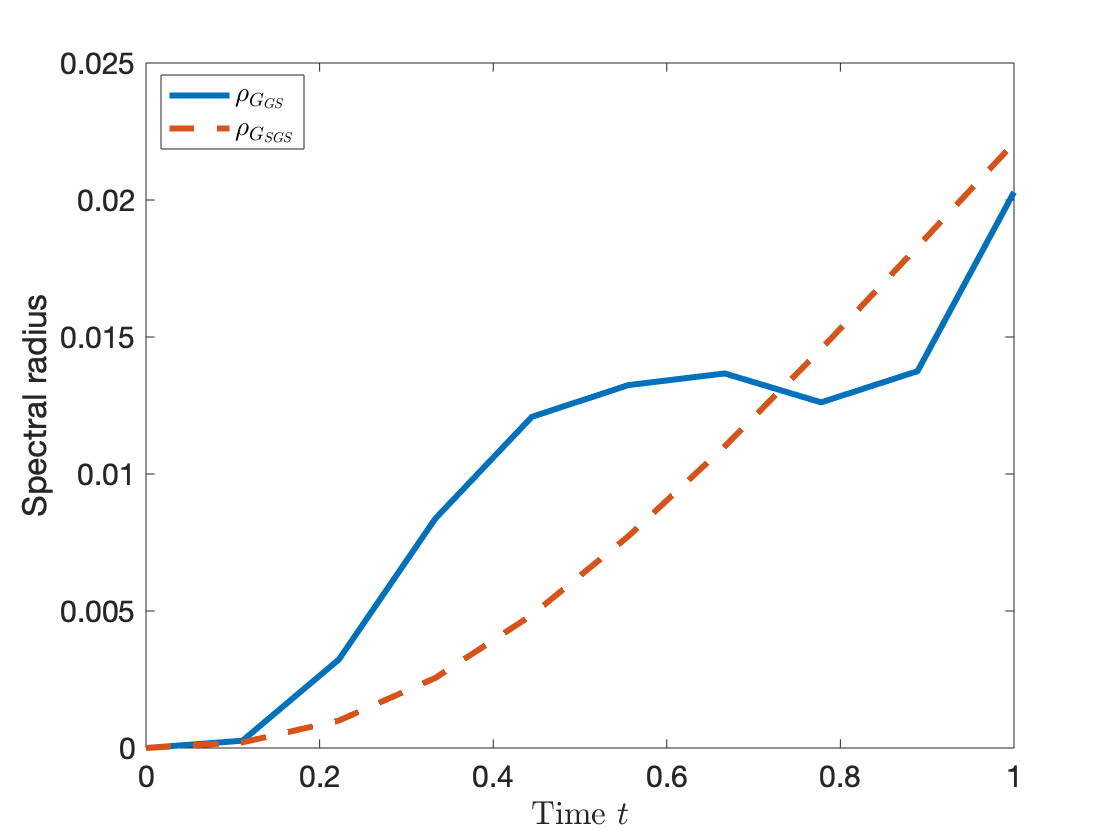}}
      \hspace{0.1in}
       \subfloat[$\alpha=0$, $N_x=1000$, $Nt=100$]{\includegraphics[width=0.35\linewidth]{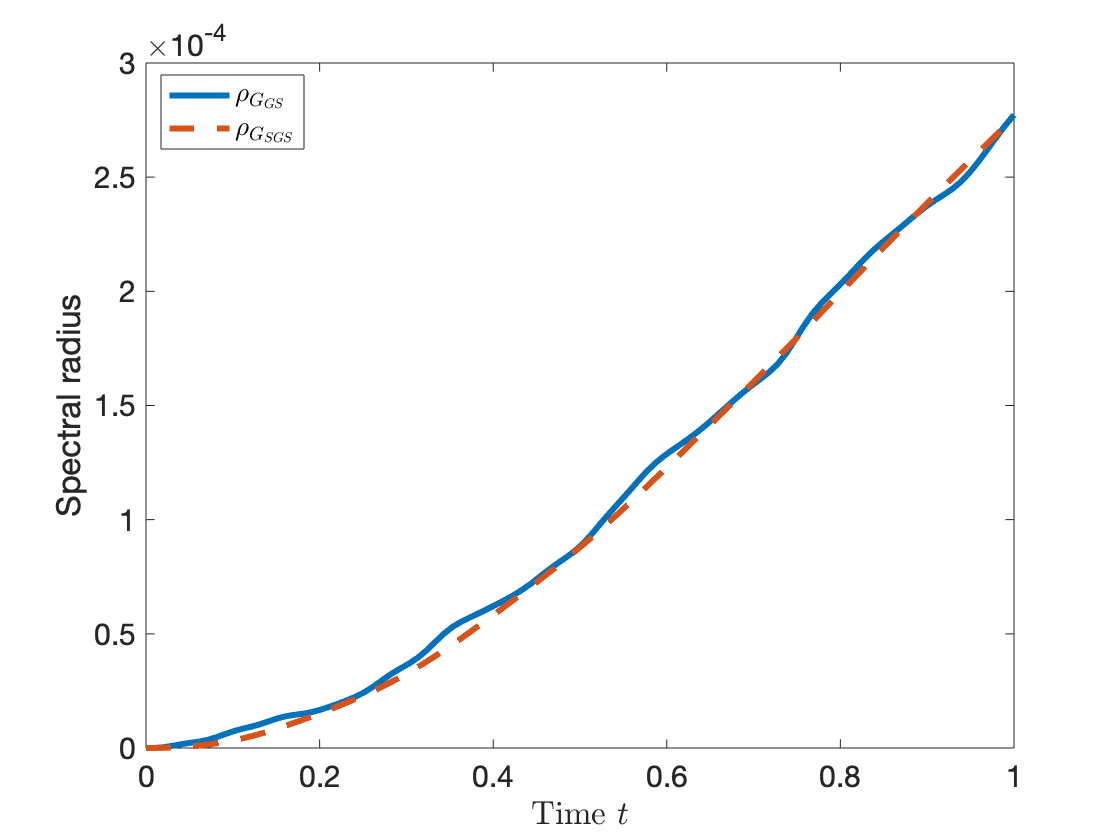}}
      \subfloat[$\alpha=0$, $N_x=20$, $Nt=1000$]{\includegraphics[width=0.35\linewidth]{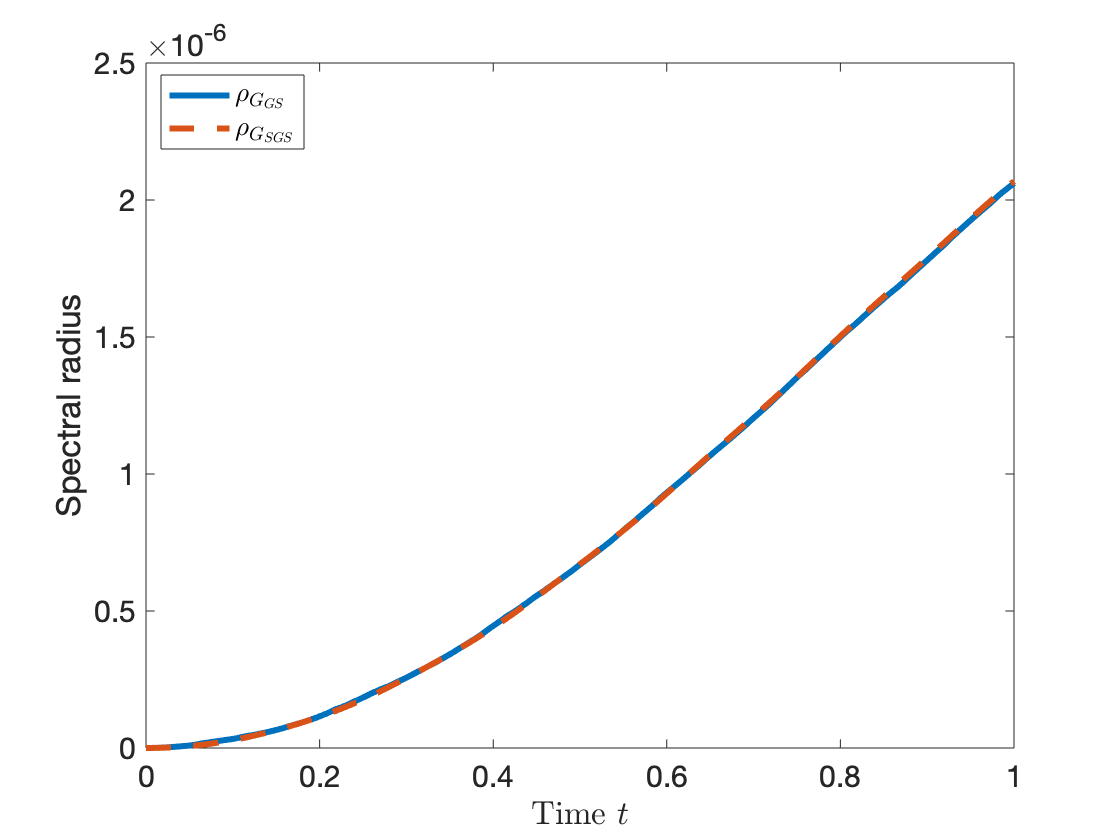}}
    \caption{The spectral radius in 1D for GSPM and SGSPM up to the final time $T=1$ given the initial condition $\m_0=[\cos(x^2(1-x)^2)\sin(0),\sin(x^2(1-x)^2)\sin(0),\cos(0)]^T$ with non-zero force term.}
    \label{fig:rho-1}
\end{figure}

The spectral radius at each time step is presented in \Cref{fig:rho-1} for GS and SGS given the initial condition $\m_0=[\cos(\cos(\pi x))\sin(0.01),\sin(\cos(\pi x))\sin(0.01),\cos(0.01)]^T$ in 1D. The results of the case in 3D is similar. The result given another initial condition $\m_0=[\cos(\cos(\pi x))\sin(0.01),\sin(\cos(\pi x))\sin(0.01),\cos(0.01)]^T$ in 1D is presented in \Cref{fig:rho-2}. The SGS spectral radius curve is more smooth than that of GS method when the time step size is large. The spectral radius of GSPM and SGSPM at the final time $T=1$ is presented in \Cref{tab:GS-SGS-rad-1D} for 1D case and \Cref{tab:GS-SGS-rad-3D} for 3D case. We observe that $\rho_{G_{SGS}}$ is smaller than $\rho_{G_{GS}}$. 

\begin{table}[htbp]
    \centering
     \caption{The spectral radius at the final time $T=1$ for $\rho(G_{GS})$ and $\rho(G_{SGS})$ with $\alpha=0.001$ in 1D.}
    \label{tab:GS-SGS-rad-1D}
    \begin{tabular}{c|c|c|c}
    \hline
    $\Delta x$     & $\Delta t$ & $\rho(G_{GS})$ & $\rho(G_{SGS})$ \\
     \hline
        $5D-4$ & $2D-2$ & $5.896039e-06$ & $3.239712e-06$ \\
        \hline
        $5D-4$ & $1.25D-2$ & $3.641269e-06$  & $1.683222e-06$ \\
        \hline
         $4D-2$ & $1D-6$ & $7.798302e-14$ & $9.728731e-15$ \\
         \hline
         $2.5D-3$ & $1D-6$ & $7.786657e-14$& $9.794043e-15$\\
    \hline
    \end{tabular}
\end{table}

\begin{table}[htbp]
    \centering
     \caption{The spectral radius at the final time $T=1$ for $\rho(G_{GS})$ and $\rho(G_{SGS})$ with $\alpha=0.001$ in 3D.}
    \label{tab:GS-SGS-rad-3D}
    \begin{tabular}{c|c|c|c}
    \hline
    $\Delta x=\Delta y=\Delta z$     & $\Delta t$ & $\rho(G_{GS})$ & $\rho(G_{SGS})$ \\
\hline

$1/2$  & $1\times10^{-4}$   & $4.872153\times10^{-11}$ & $2.696574\times10^{-11}$ \\

$1/8$  & $1\times10^{-4}$   & $5.315656\times10^{-10}$ & $1.106243\times10^{-10}$ \\

$1/16$ & $1\times10^{-4}$   & $5.550360\times10^{-10}$ & $1.203321\times10^{-10}$ \\

$1/20$ & $2.5\times10^{-3}$ & $8.492480\times10^{-7}$  & $2.052615\times10^{-7}$ \\

\hline
    \end{tabular}
\end{table}

\begin{figure}[htbp]
    \centering
    \subfloat[$\alpha=0$, $N_x=200$, $Nt=10$]{\includegraphics[width=0.35\linewidth]{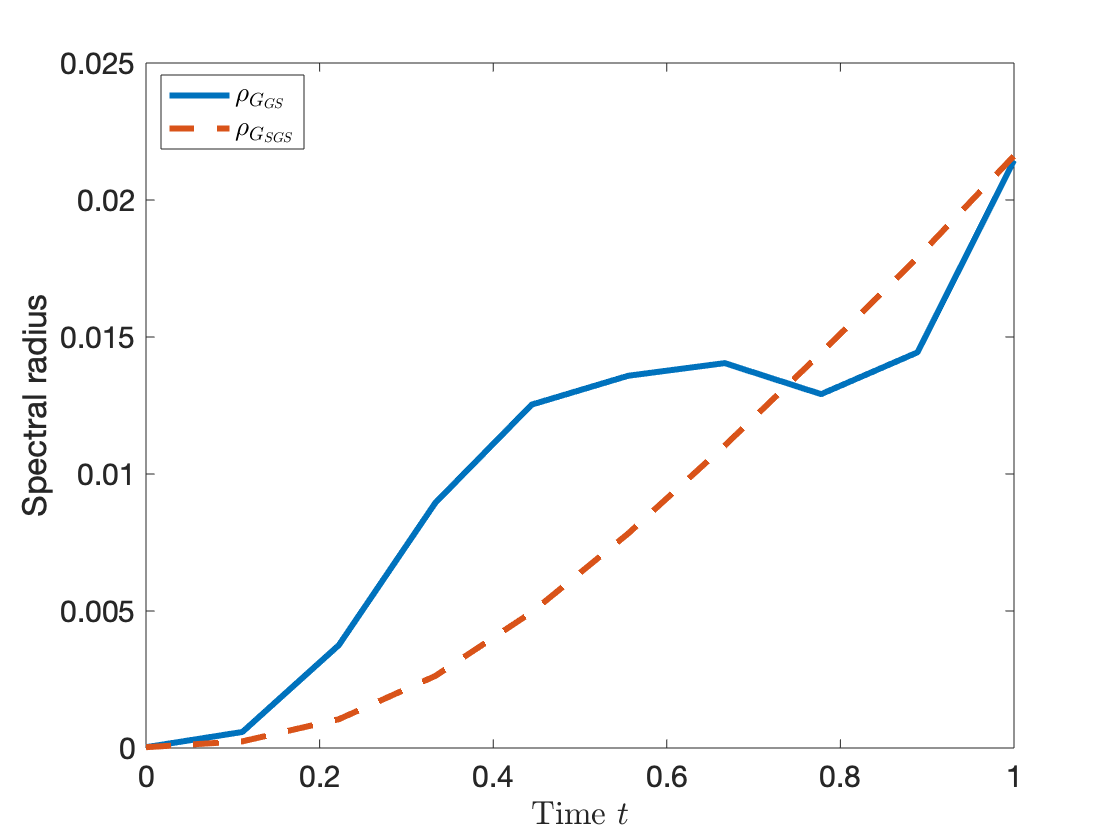}}
     \subfloat[$\alpha=0$, $N_x=200$, $Nt=100$]{\includegraphics[width=0.35\linewidth]{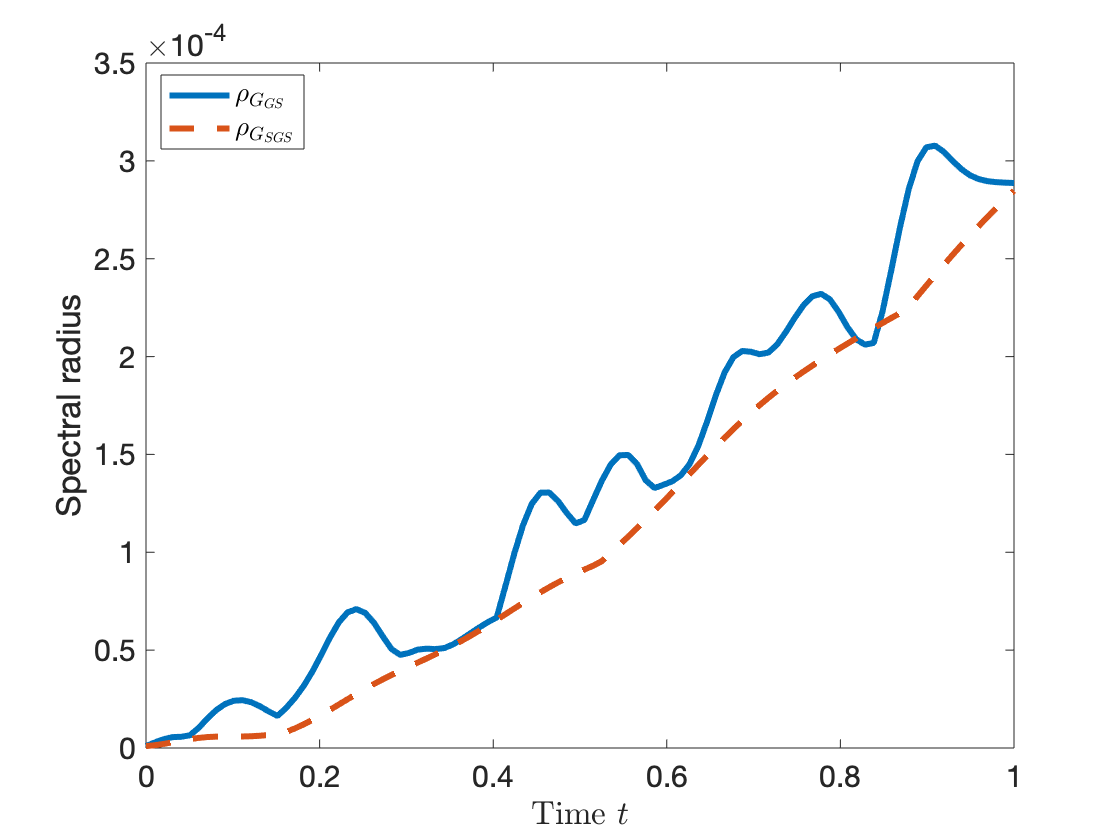}}
      \subfloat[$\alpha=0$, $N_x=1000$, $Nt=10$]{\includegraphics[width=0.35\linewidth]{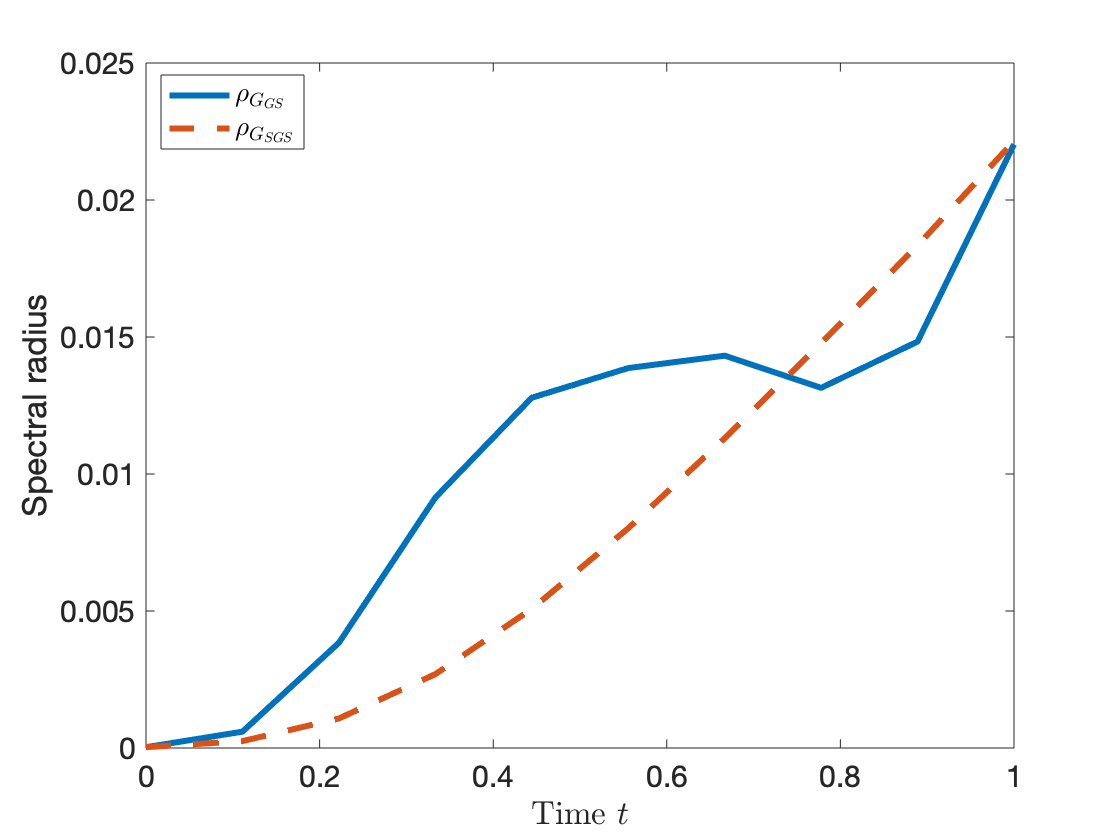}}
      \hspace{0.1in}
      \subfloat[$\alpha=0$, $N_x=200$, $Nt=10$]{\includegraphics[width=0.35\linewidth]{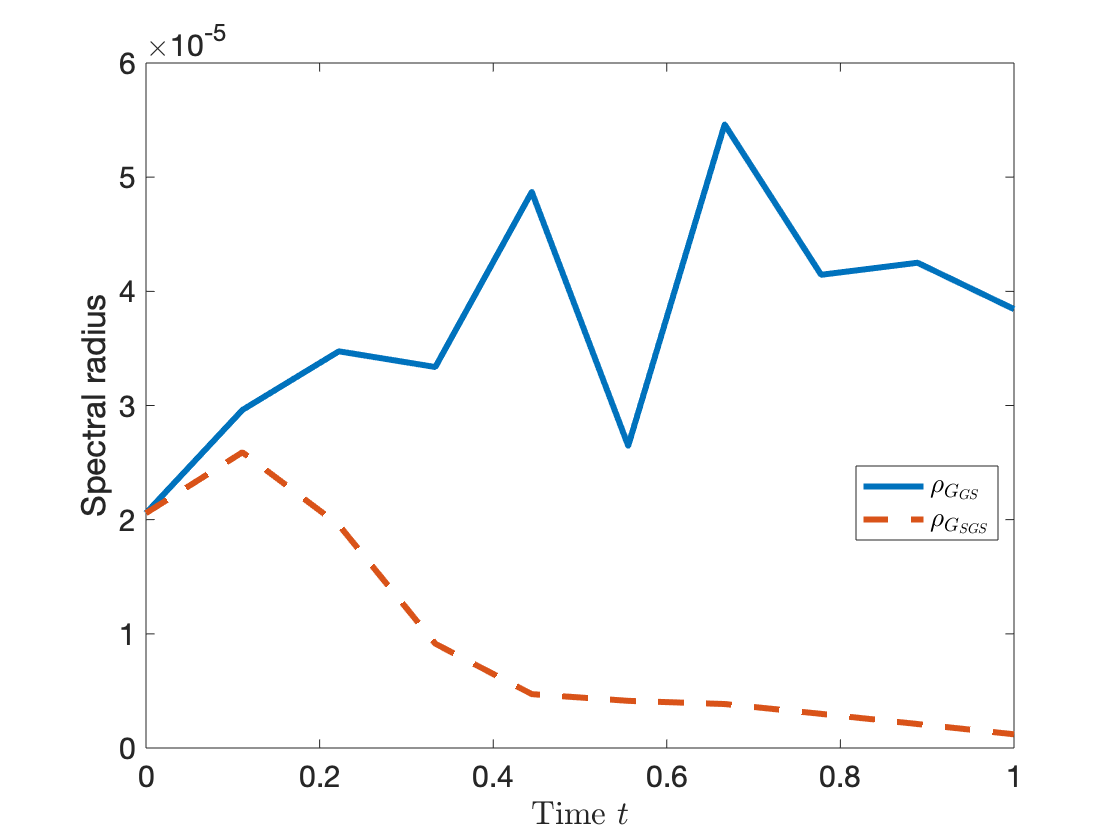}}
     \subfloat[$\alpha=0$, $N_x=200$, $Nt=100$]{\includegraphics[width=0.35\linewidth]{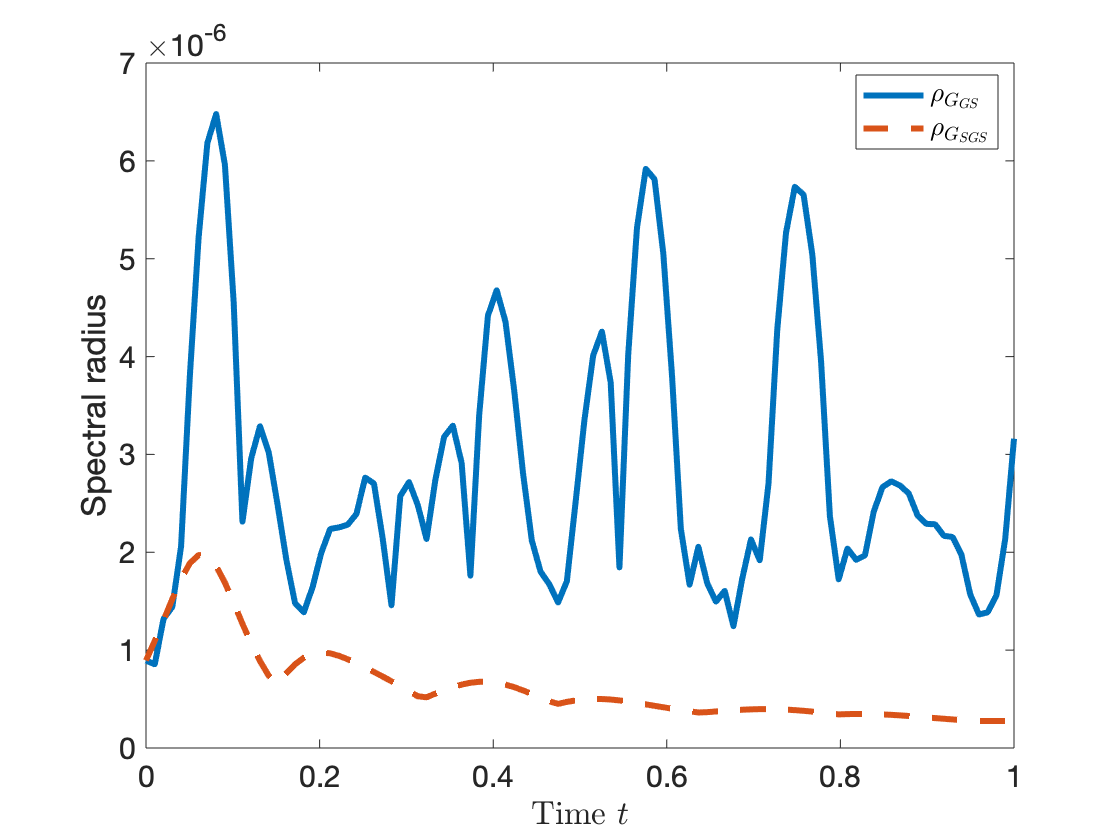}}
      \subfloat[$\alpha=0$, $N_x=1000$, $Nt=10$]{\includegraphics[width=0.35\linewidth]{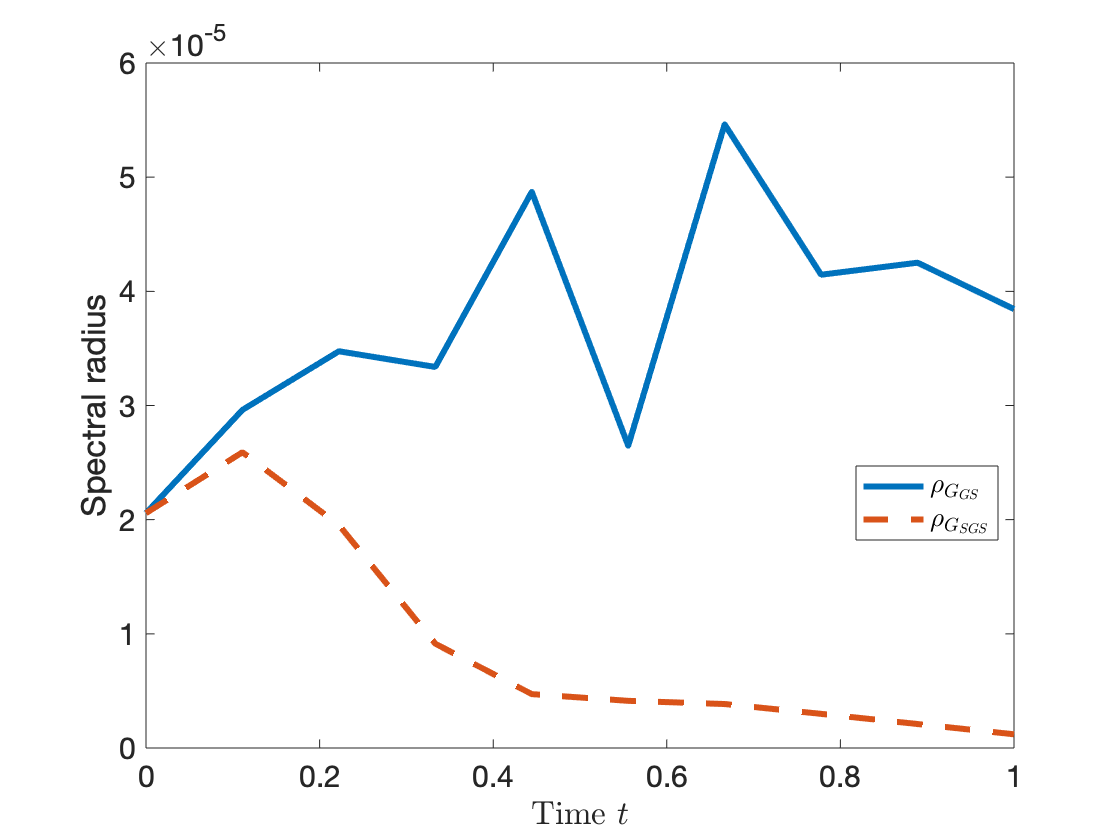}}
    \caption{The spectral radius in 1D for GSPM and SGSPM up to the final time $T=1$ given the initial condition $\m_0=[\cos(\cos(\pi x))\sin(0.01),\sin(\cos(\pi x))\sin(0.01),\cos(0.01)]^T$ with non-zero force term in top rows and zero force term in bottom rows.}
    \label{fig:rho-2}
\end{figure}

\section{Numerical experiments}
\label{sec:experiments}

\subsection{Accuracy and efficiency for SGSPM}

For the model problem, analytical exact solutions are derived for both one-dimensional (1D) and three-dimensional (3D) scenarios to serve as benchmarks for error quantification.

For the 1D case, the exact magnetization solution $\boldsymbol{m}_e$ is:
\[
\boldsymbol{m}_e = \big(\cos\big(x^2(1-x)^2\big)\sin t,\, \sin\big(x^2(1-x)^2\big)\sin t,\, \cos t\big)^T,
\]
while the corresponding 3D exact solution is:
\[
\boldsymbol{m}_e = \big(\cos(XYZ) \sin t,\, \sin(XYZ) \sin t,\, \cos t\big)^T,
\]
where $X = x^2(1-x)^2$, $Y = y^2(1-y)^2$ and $Z = z^2(1-z)^2$.

These exact solutions satisfy the governing equation \eqref{c1-large} with $\h_{\textbf{eff}}=\Delta \m$ when the forcing term is defined as
\[
\boldsymbol{f}_e = \partial_t \boldsymbol{m}_e + \boldsymbol{m}_e \times \Delta \boldsymbol{m}_e + \alpha \boldsymbol{m}_e \times (\boldsymbol{m}_e \times \Delta \boldsymbol{m}_e).
\]
They also comply with the homogeneous Neumann boundary condition, ensuring consistency with simulation constraints. The ferromagetic body $\Omega=[0,1]$ in 1D and $\Omega=[0,1]^3$ in 3D. The final time $T=1$. The temporal accuracy in 1D and 3D for SGSPM and GSPM are presented in \Cref{tab:temp_acc_1d} and \Cref{tab:temp_acc_3d}, respectively. The spatial accuracy with in 1D and 3D for SGSPM and GSPM are presented in \Cref{tab:spatial_acc_1d} and \Cref{tab:cpu}, respectively. These results indicate that the GSPM and SGSPM are both the first-order accurate in time and the second-order accurate in space in both 1D and 3D settings.

\begin{table}[htbp]
  \centering
  \caption{Temporal accuracy in 1D for GSPM, and SGSPM when $\Delta x = 5D-4$ and $\alpha = 1D-3$ up to the final time $T=1$.}
  \begin{tabular}{|c|c|c|c|c|}
    \hline
    $\Delta t$&SGSPM & wall time (s) & GSPM & wall time (s) \\
    &$\|m_h - m_e\|_\infty$ & & $\|m_h - m_e\|_\infty$ &  \\
    \hline
    2.0D-2 &0.031993635118710 & 0.186551 & 0.032024319826967 & 0.105561  \\
    1.0D-2 &0.016258177977747 & 0.204056 & 0.016289023892517 & 0.144290 \\
    5.0D-3 &0.008198757469279 & 0.307139 & 0.007989189016688 & 0.293836 \\
    2.5D-3 &0.004140400463721 &0.622396  & 0.004192511571719 & 0.450300  \\
    1.25D-3 &0.002088708253334 & 1.590927 & 0.002066422804009 & 1.060299 \\
    \hline
    order &0.9848 & -- & 0.9866 & -- \\
    \hline
  \end{tabular}
  \label{tab:temp_acc_1d}
\end{table}

\begin{table}[htbp]
  \centering
  \caption{Spatial accuracy in 1D for SGSPM and GSPM when $\Delta t = 1D-6$ and $\alpha = 1D-3$ up to the final time $T=1$.}
  \begin{tabular}{|c|c|c|c|c|}
    \hline
    $\Delta x$& SGSPM & wall time (s) & GSPM & wall time (s) \\
    &$\|m_h - m_e\|_\infty$ & & $\|m_h - m_e\|_\infty$ &  \\
    \hline
    4.0D-2 &6.211240439639720e-04 & 145.219485 & 6.211187739000312e-04 &99.787014  \\
    2.0D-2 &1.553598118020361e-04 & 173.770068 & 1.553560782340258e-04 & 128.938661 \\
    1.0D-2 &3.898838676885674e-05 & 200.049600 & 3.898639391297642e-05 & 158.970558 \\
    5.0D-3 &9.896762845240659e-06 & 269.238673 & 9.896118137620036e-06 & 195.025552 \\
    2.5D-3 &2.624275358464878e-06 &423.971576  & 2.628912127655703e-06 & 335.403948 \\
    \hline
    order & 1.9746&-- & 1.9741 & -- \\
    \hline
  \end{tabular}
  \label{tab:spatial_acc_1d}
\end{table}

\begin{table}[htbp]
  \centering
  \caption{Temporal accuracy in 3D for GSPM, and SGSPM when $\Delta x=\Delta y=\Delta z = 1/16$ and $\alpha = 1D-3$ up to the final time $T=1$.}
  \begin{tabular}{|c|c|c|c|c|}
    \hline
    $\Delta t$&SGSPM & wall time (s) & GSPM & wall time (s) \\
    &$\|m_h - m_e\|_\infty$ & & $\|m_h - m_e\|_\infty$ &  \\
    \hline
    $1/32$ &2.877587588594097e-04 & 3.938365 & 3.002439695850559e-04  & 4.502316 \\
    $1/64$ &1.473122547213623e-04 & 6.640514 & 1.521432170386424e-04 & 6.102772 \\
    $1/128$ &7.501767594611609e-05 & 13.782621 & 7.736830001073121e-05 & 10.351400 \\
    $1/256$ &3.832297611698108e-05 &25.012785  & 3.817108461035950e-05 & 18.854986 \\
    $1/512$ &2.005449155901271e-05 & 47.747254 & 2.029520809002309e-05 & 39.447516 \\
    \hline
    order & 0.9628& -- & 0.9769 & -- \\
    \hline
  \end{tabular}
  \label{tab:temp_acc_3d}
\end{table}

The comparison of the GSPM and SGSPM is presented in \Cref{tab:cpu} which indicates that the GSPM method is slightly faster than SGSPM method. Both methods are fast and efficient since the FFT has been used for the Poisson solver. The difference between two methods for the cost of computation are the number of solving the heat equation.

\begin{table}[htbp]
  \centering
  \caption{Spatial accuracy in 3D for SGSPM and GSPM when $\Delta t = 1D-4$ and $\alpha = 1D-3$ up to the final time $T=1$.}
  \begin{tabular}{|c|c|c|c|c|}
    \hline
    $\Delta x=\Delta y=\Delta z$& SGSPM & wall time (s) & GSPM & wall time (s) \\
    &$\|m_h - m_e\|_\infty$ & & $\|m_h - m_e\|_\infty$ &  \\
    \hline
    $1/2$ &3.585582209200133e-04 &2.731247  & 3.585582209190141e-04 &2.389967  \\
    $1/4$ &8.123800347159538e-05 & 14.414304 & 8.124627756622704e-05 &14.075774  \\
    $1/8$ & 2.018671673698069e-05&  148.712626& 2.018675673998160e-05 & 124.770155 \\
    $1/16$ &5.040436690473271e-06 &957.476236  & 5.051941406541971e-06 & 752.492960 \\
    \hline
    order & 2.0466& --&2.0457 &--  \\
    \hline
  \end{tabular}
  \label{tab:cpu}
\end{table}

\begin{figure}[htbp]
    \centering
    \subfloat[1D, time]{\includegraphics[width=0.35\linewidth]{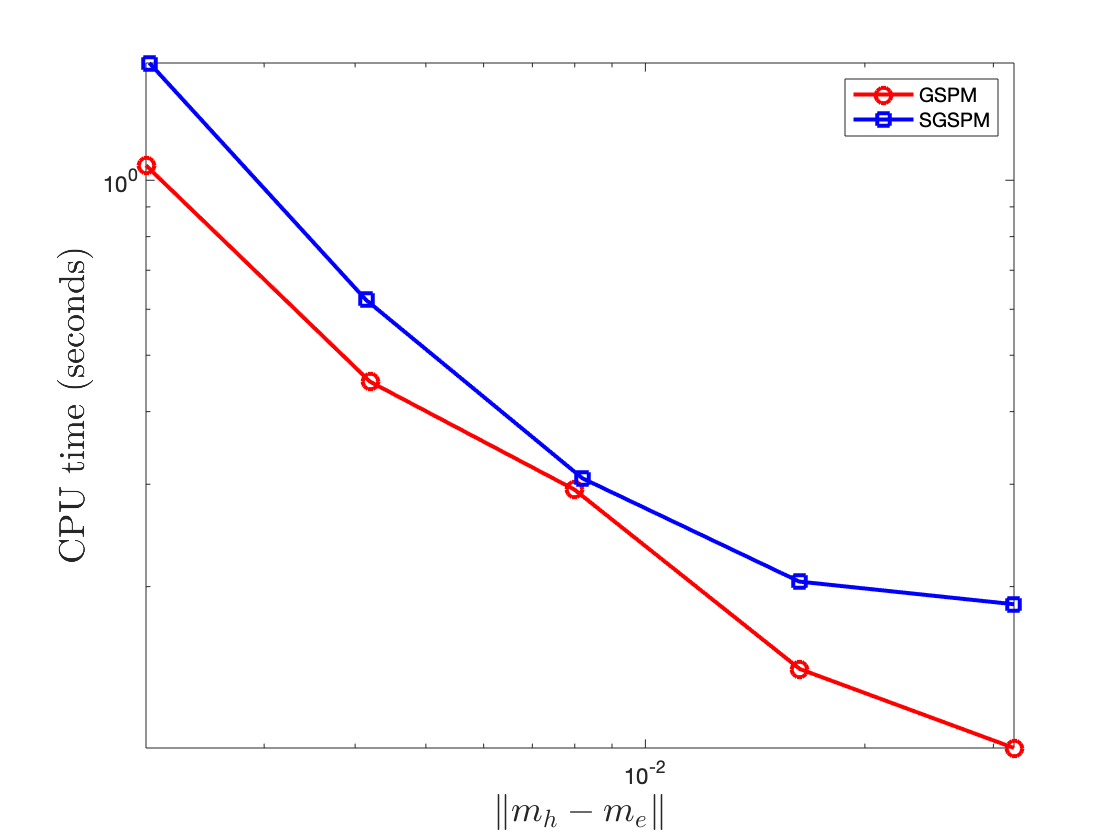}}
     \subfloat[1D, space]{\includegraphics[width=0.35\linewidth]{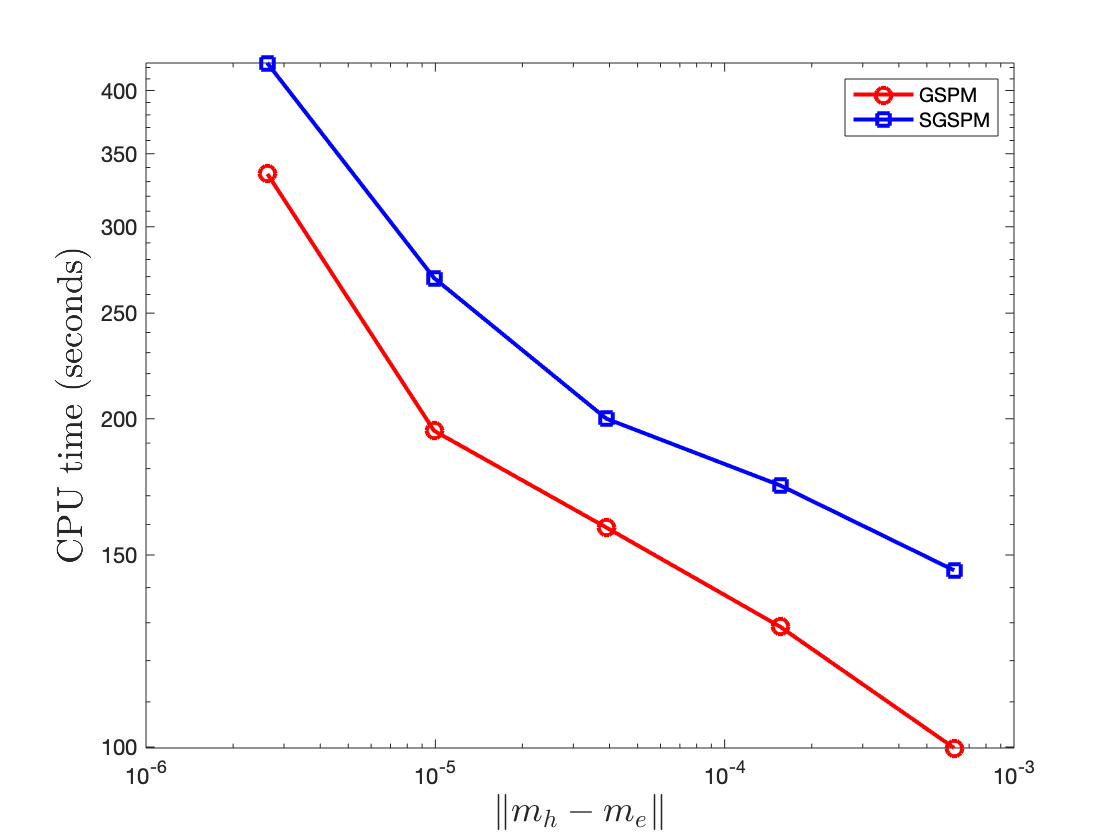}}
     \hspace{0.1in}
     \subfloat[3D, time]{\includegraphics[width=0.35\linewidth]{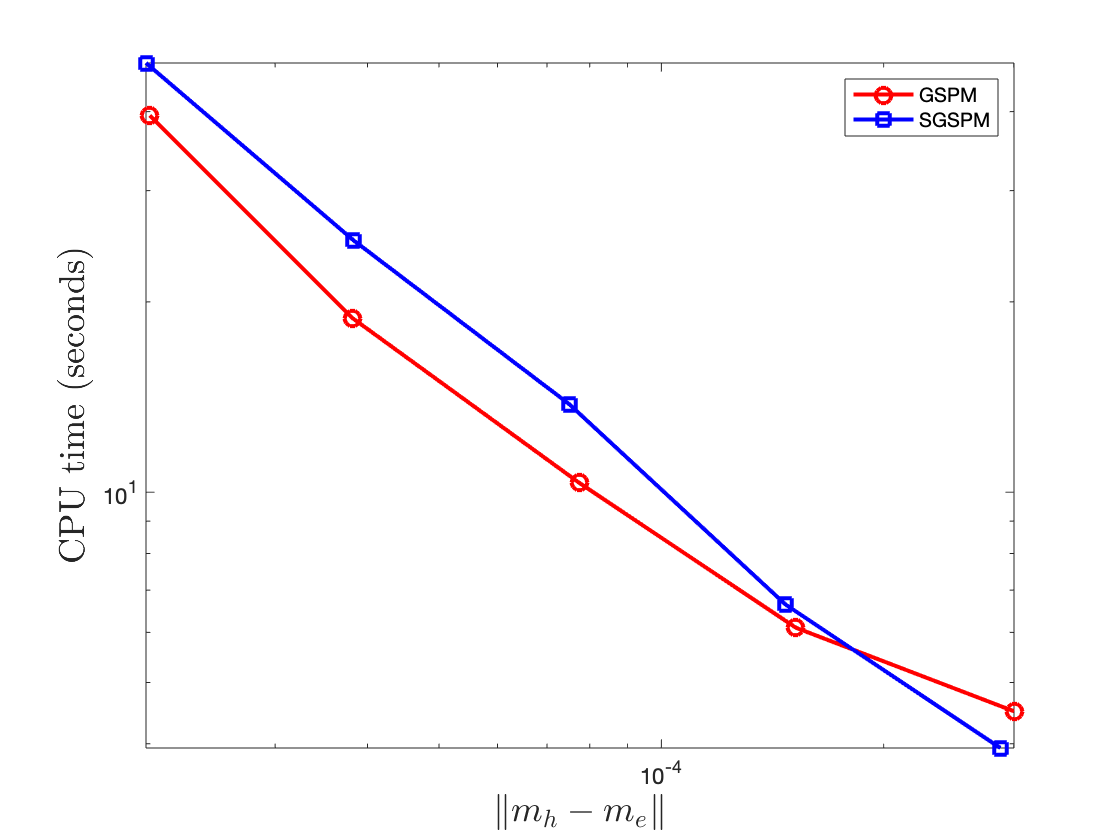}}
     \subfloat[3D, space]{\includegraphics[width=0.35\linewidth]{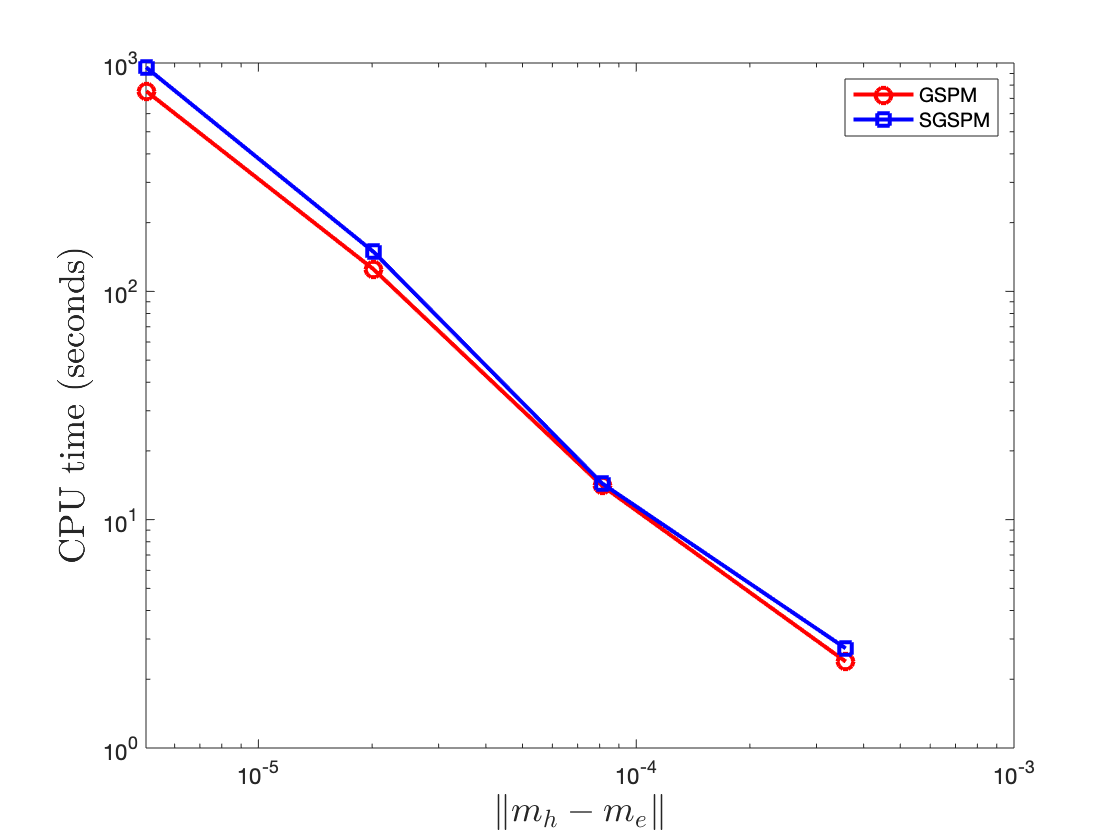}}
    \caption{The comparison of efficiency of GSPM and SGSPM for 1D and 3D. The GSPM is slightly faster than SGSPM.}
    \label{fig:rho-2}
\end{figure}

\subsection{The simulations with only exchange field}

\begin{figure}[htbp]
    \centering
    \subfloat[GS, $\alpha=0$]{\includegraphics[width=0.3\linewidth]{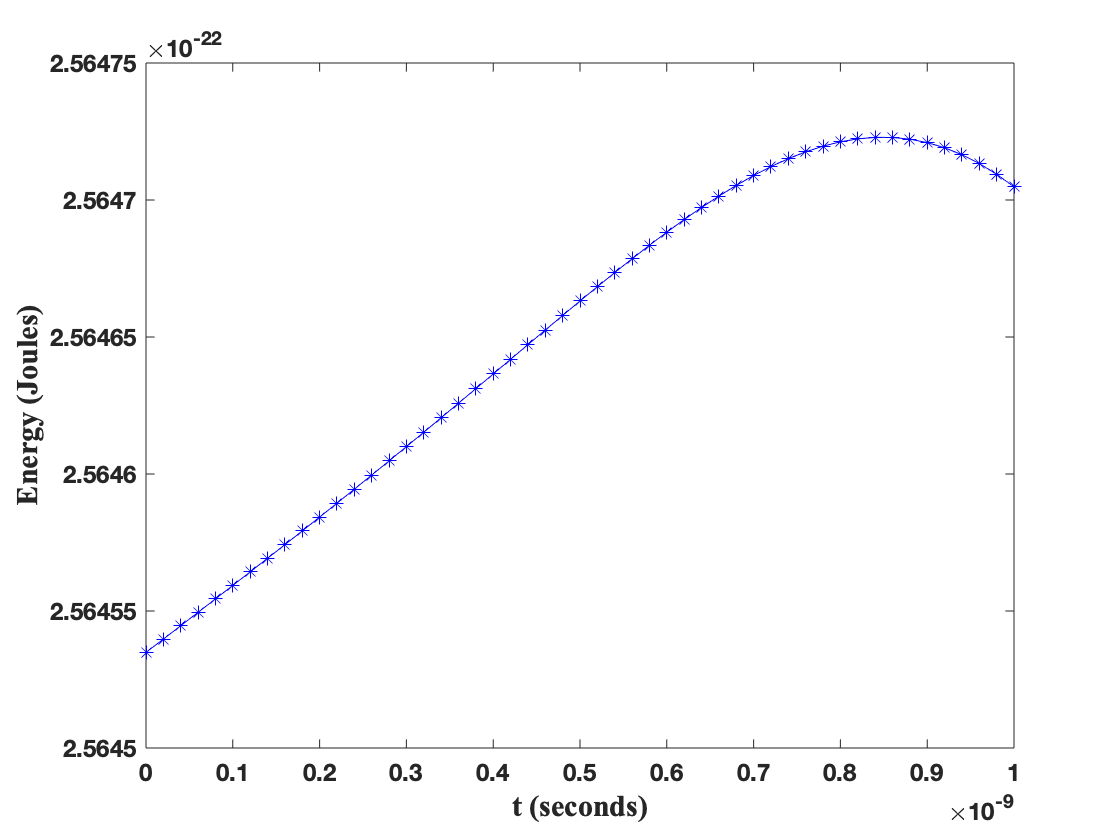}}
    \subfloat[GS, $\alpha=0.005$]{\includegraphics[width=0.3\linewidth]{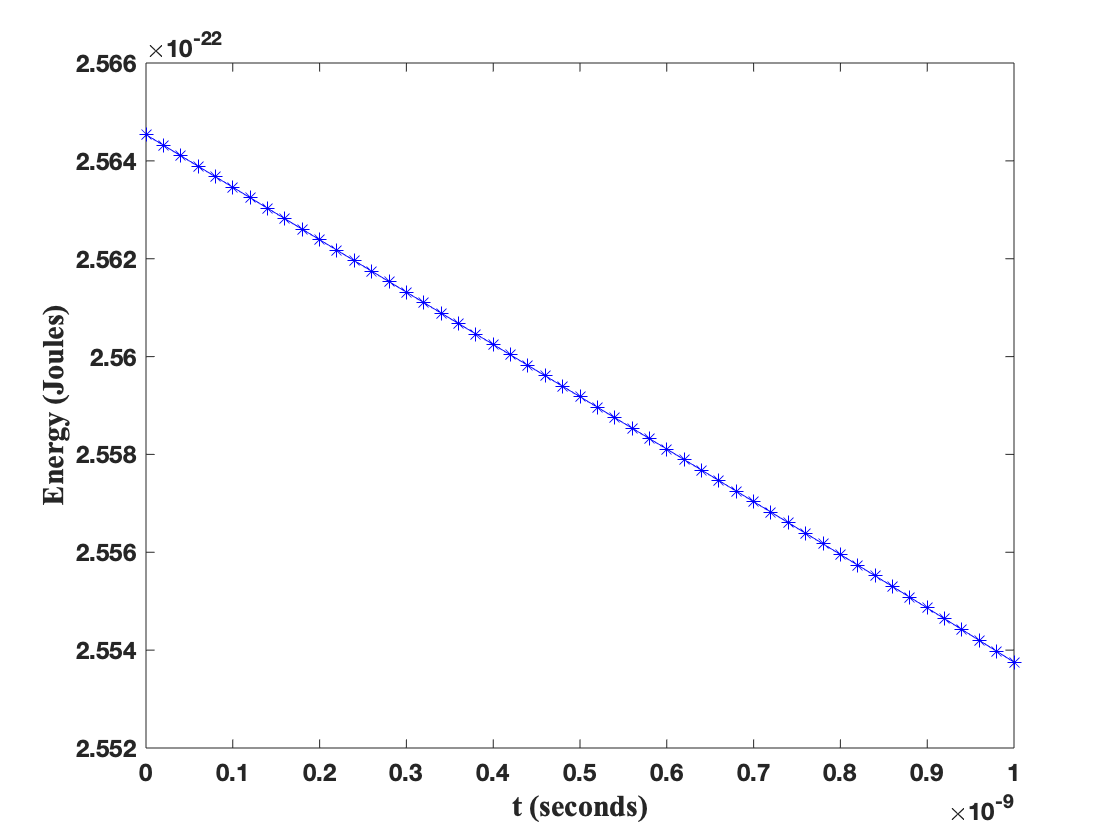}}
    \subfloat[GS, $\alpha=0.01$]{\includegraphics[width=0.3\linewidth]{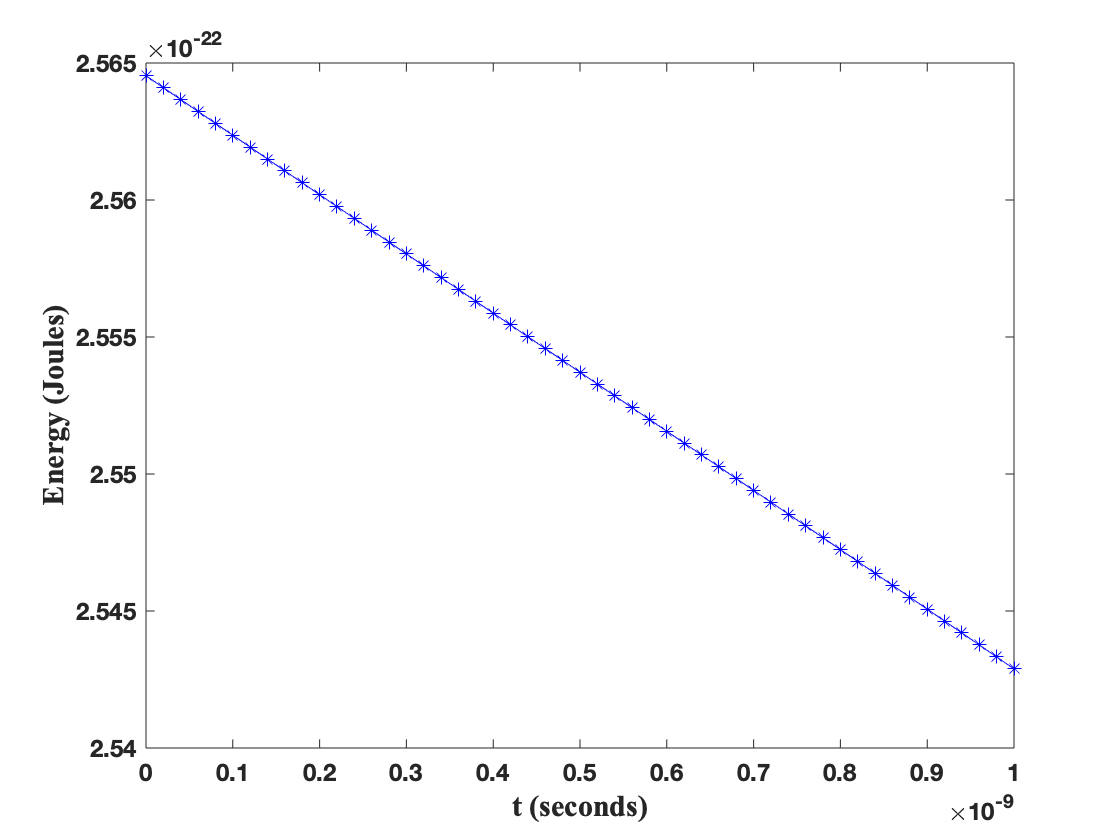}}
    \hspace{0.1in}
    \subfloat[SGS, $\alpha=0$]{\includegraphics[width=0.3\linewidth]{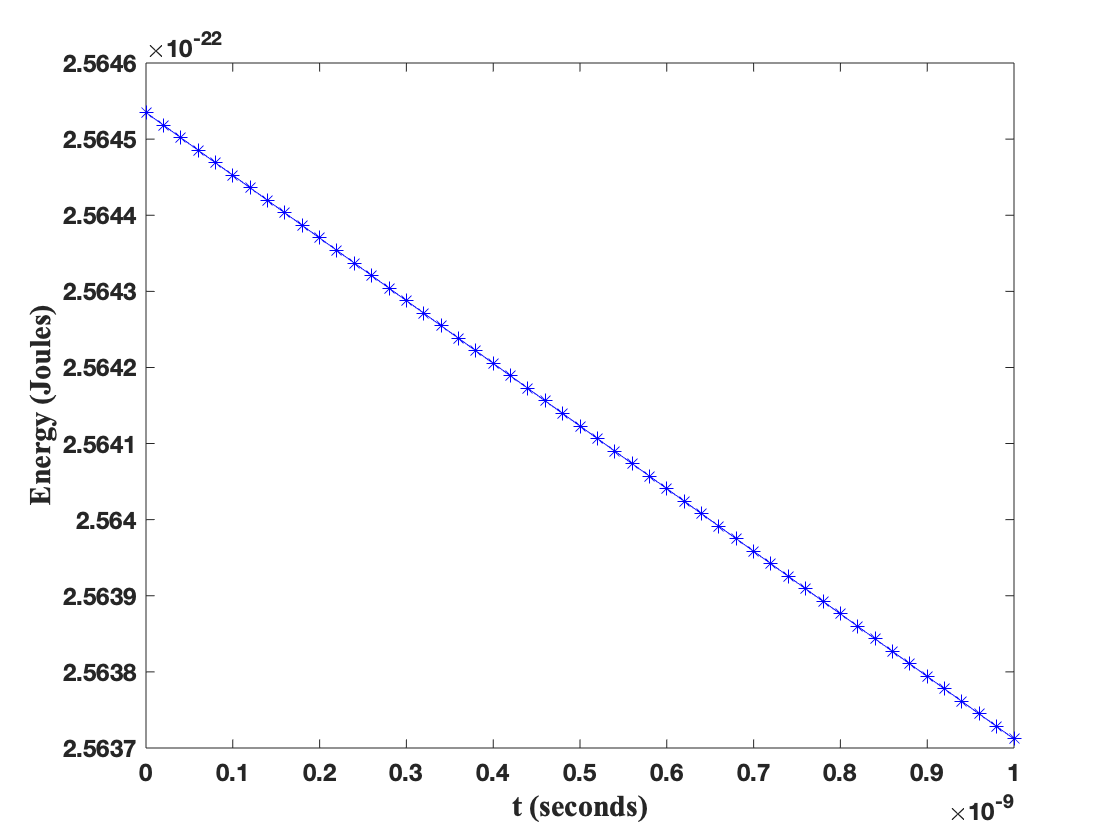}}
     \subfloat[SGS, $\alpha=0.005$]{\includegraphics[width=0.3\linewidth]{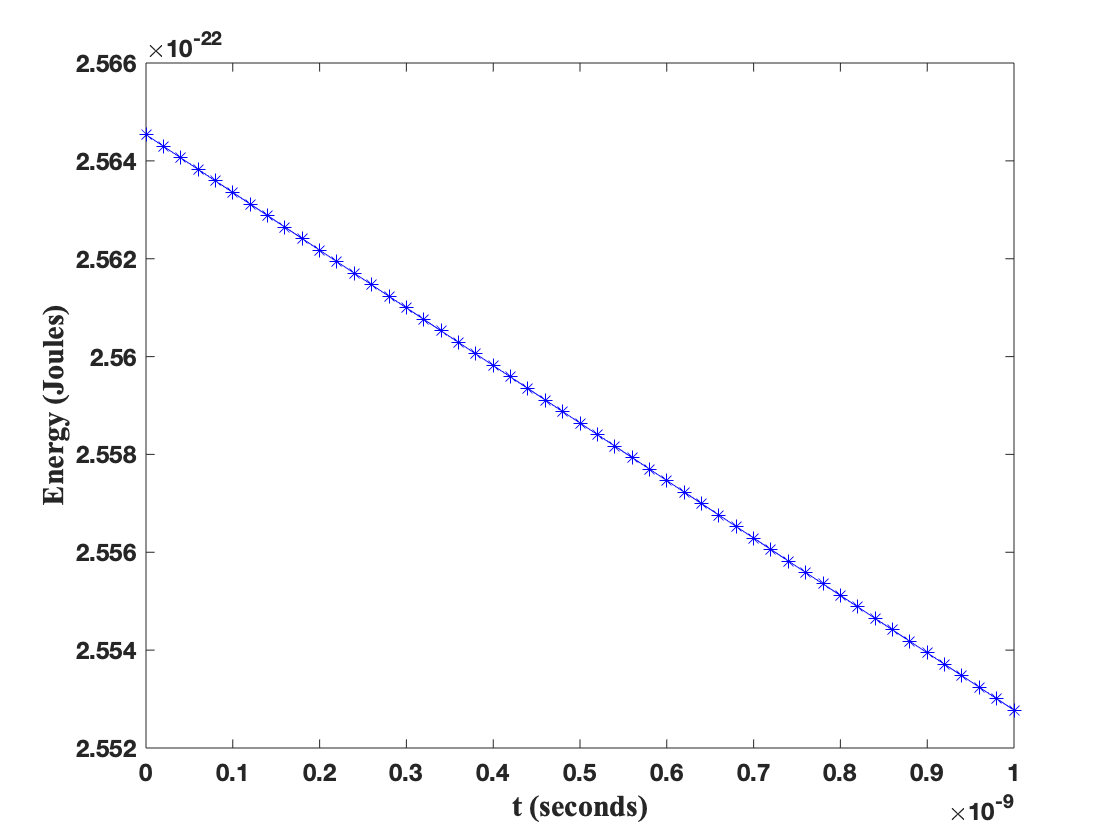}}
    \subfloat[SGS, $\alpha=0.01$]{\includegraphics[width=0.3\linewidth]{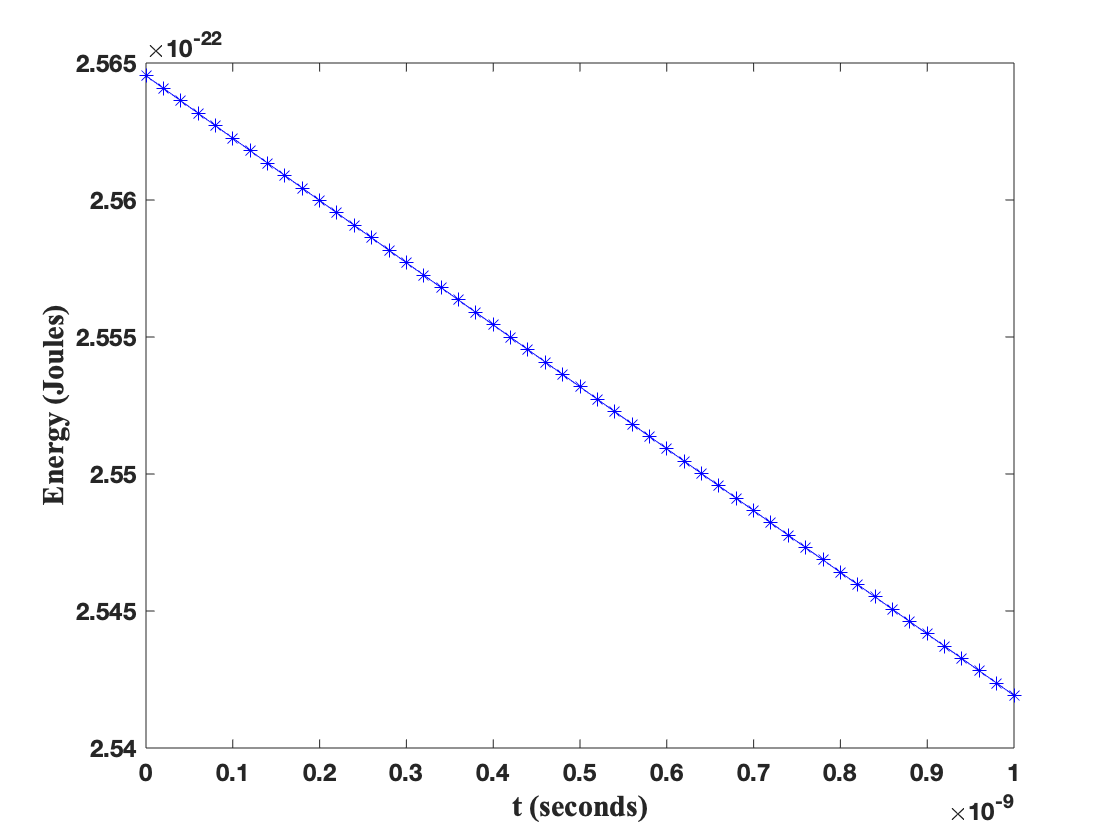}}
    \caption{The energy evolution for GS and SGS projection method with $\alpha=0,0.005,0.01$. The initial condition is specified by $\m_0=[\cos(\cos(\pi x))\sin(0.01),\sin(\cos(\pi x))\sin(0.01),\cos(0.01)]^T$.}
    \label{fig:energy-SGS-2}
\end{figure}

In this section, we take the field $\f=0$, say that the effective field only in the form of $\h_{\textbf{eff}}=\epsilon \Delta \m$. To see the performance of the stability, we kick off the initial condition as below,
\begin{align}\label{eq-m0-1}
    \m_0=[\cos(\cos(\pi x))\sin(0.01),\sin(\cos(\pi x))\sin(0.01),\cos(0.01)]^T.
\end{align}
The simulation model uses a ferromagnetic thin film with dimensions $480\times 480 \times 20\;nm^3$ discretized on a $100\times 100\times 4$ grid. The temporal step size we take is $\Delta t=1\;ps$. The result of energy evolution by GSPM method and SGSPM method for the equilibrium at final time $T=1\;ns$ is presented in \Cref{fig:energy-SGS-2}. The SGS projection scheme demonstrates superior energy stability compared to the standard GS method: SGS preserves the physical energy dissipation property for both zero and non-zero damping, while the standard GS method fails to guarantee energy decay when $\alpha=0$.


\begin{figure}[htbp]
    \centering
    \subfloat[GS, Initial S state, arrow profile]{\includegraphics[width=0.25\linewidth]{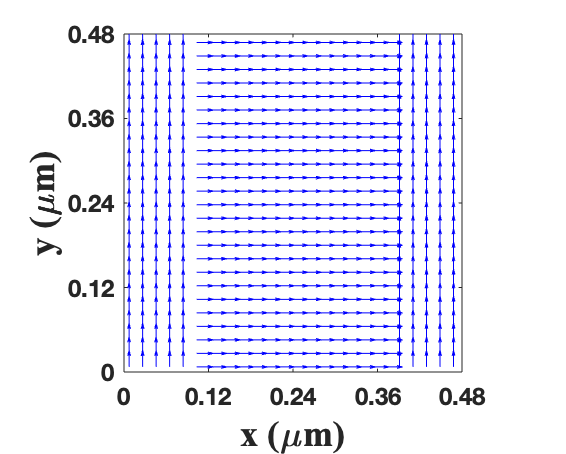}}
    \subfloat[GS, arrow profile, $\alpha=0$]{\includegraphics[width=0.25\linewidth]{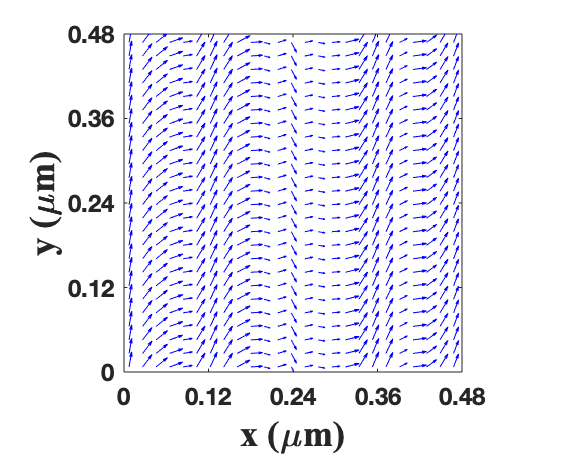}}
    \subfloat[GS, arrow profile, $\alpha=0.01$]{\includegraphics[width=0.25\linewidth]{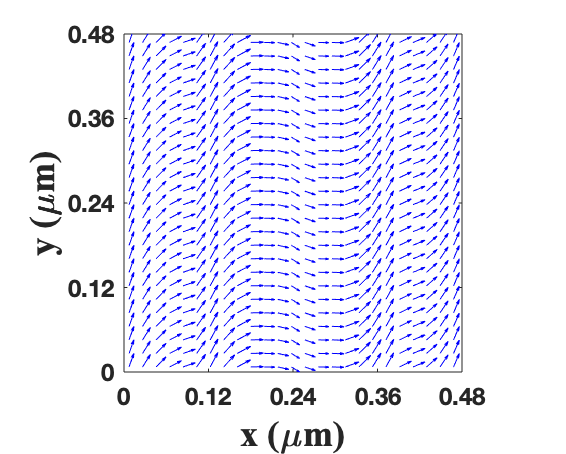}}
    \subfloat[GS, arrow profile, $\alpha=0.1$]{\includegraphics[width=0.25\linewidth]{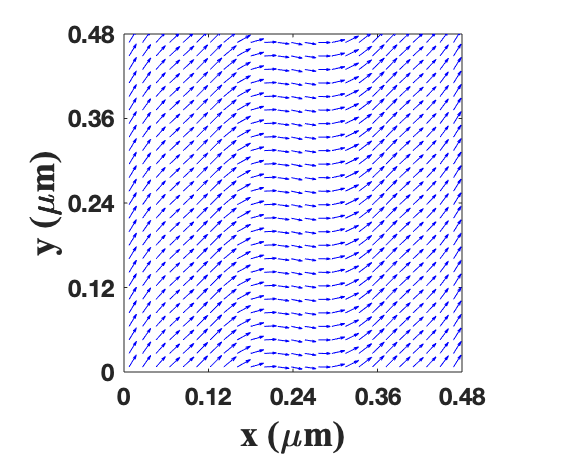}}
    \hspace{0.1in}
    \subfloat[SGS, Initial S state, arrow profile]{\includegraphics[width=0.25\linewidth]{arrow_initial_S_v1.png}}
    \subfloat[SGS, arrow profile, $\alpha=0$]{\includegraphics[width=0.25\linewidth]{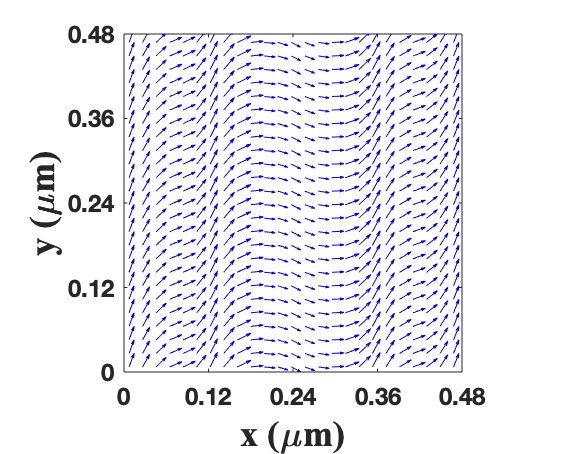}}
    \subfloat[SGS, arrow profile, $\alpha=0.01$]{\includegraphics[width=0.25\linewidth]{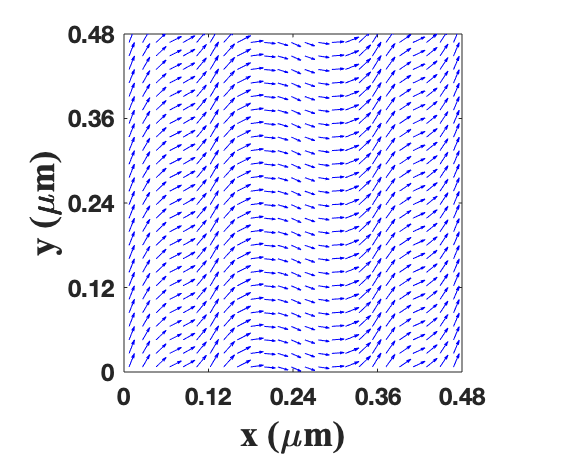}}
    \subfloat[SGS, arrow profile, $\alpha=0.1$]{\includegraphics[width=0.25\linewidth]{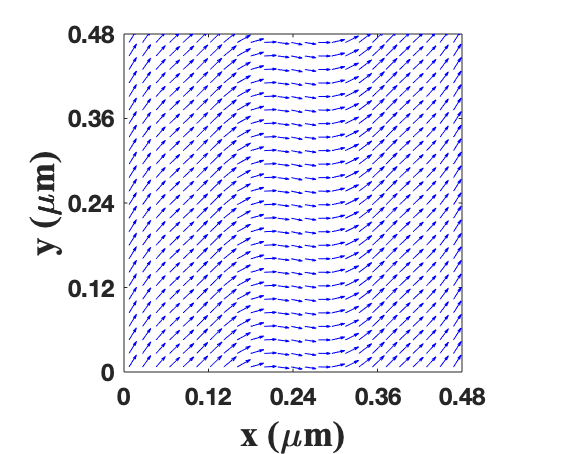}}
    \caption{GSPM and SGSPM for the simulation only for exchange field with initial condition $\m_0=[\cos(\cos(\pi x))\sin(0.01),\sin(\cos(\pi x))\sin(0.01),\cos(0.01)]^T$. The final time $T=1\;ns$, the time step $\Delta t=1\;ps$. The top row is GS with $\alpha=0,0.01,0.1$, the bottom row is SGS with $\alpha=0,0.01,0.1$.}
    \label{fig:GS-exchange-1}
\end{figure}

\begin{figure}[htbp]
    \centering
    \subfloat[GS]{\includegraphics[width=0.4\linewidth]{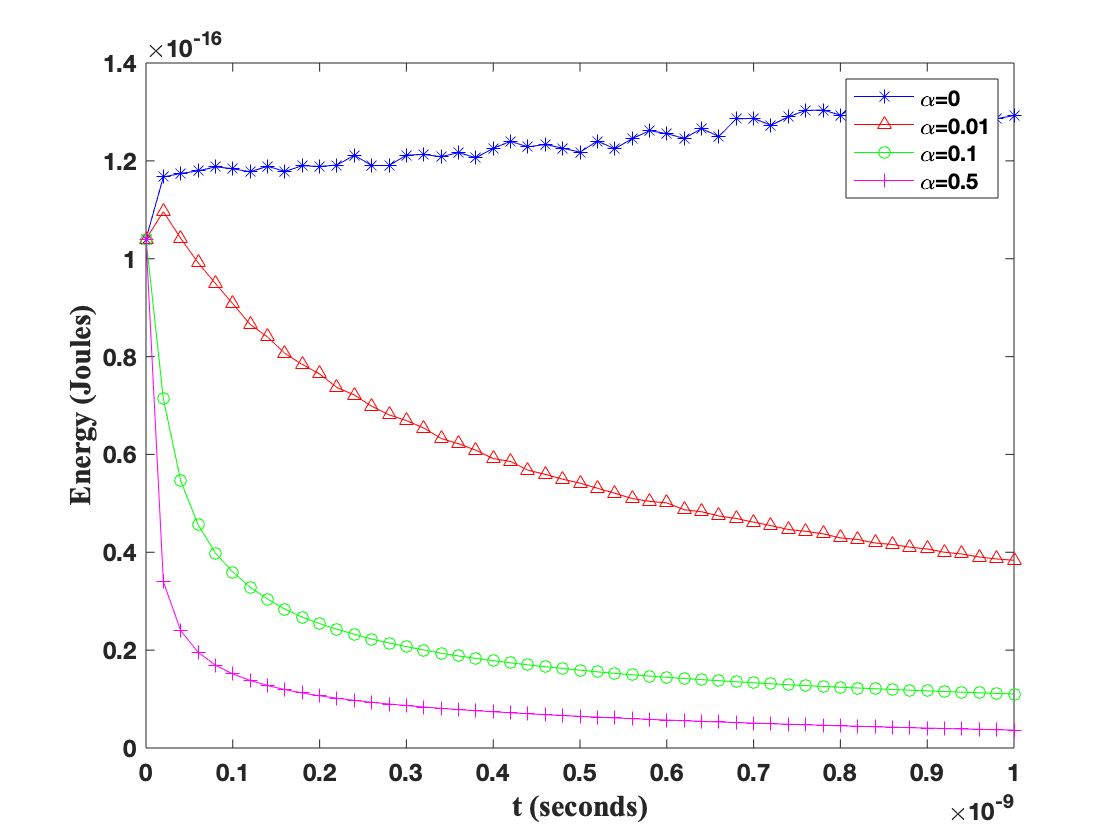}}
    \hspace{0.1in}
      \subfloat[SGS]{\includegraphics[width=0.4\linewidth]{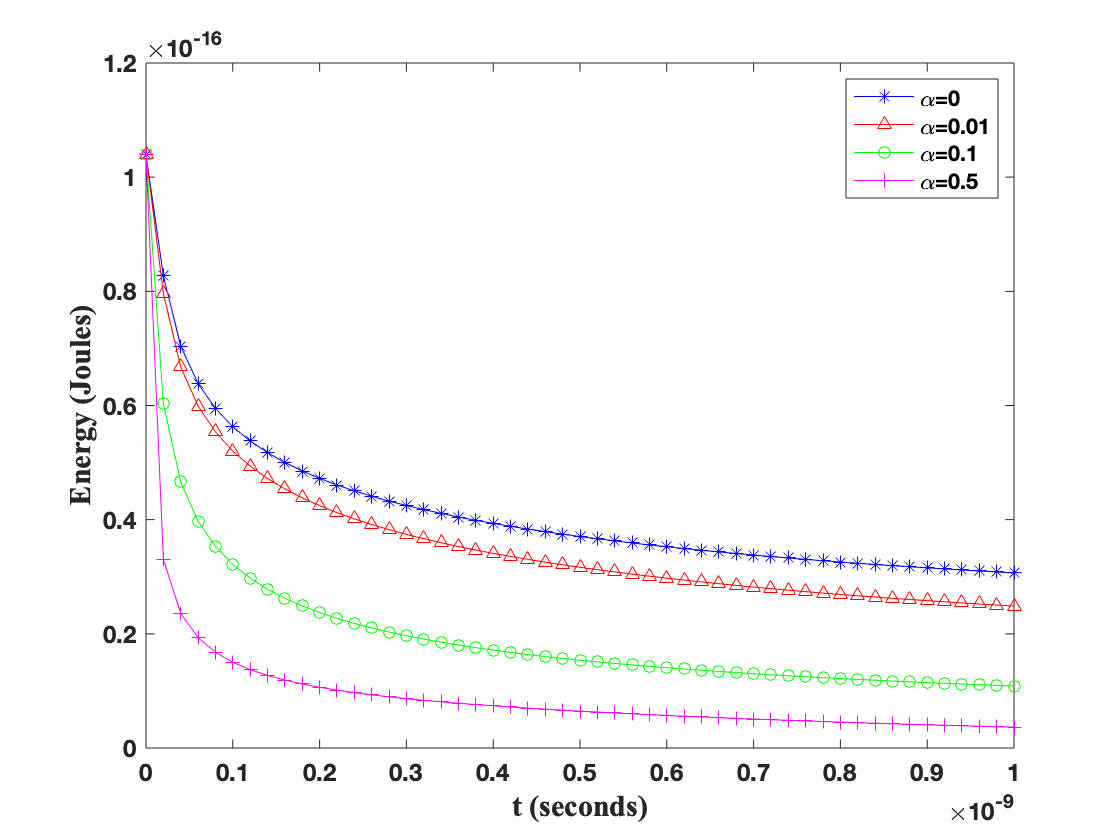}}
    \caption{The energy evolution for GS and SGS projection method with $\alpha=0,0.01,0.1,0.5$. The initial condition is specified as a S state.}
    \label{fig:energy-SGS-1}
\end{figure}

For robustness, we change the initial condition to be a S state specified in the left panel in \Cref{fig:GS-exchange-1}. The setup of magnetization simulation is the same as before.
The corresponding magnetization profile using GSPM and SGSPM with different damping parameters $\alpha=0,0.01,0.1$ is shown in \Cref{fig:GS-exchange-1}. The energy evolution using GSPM and SGSPM is presented in \Cref{fig:energy-SGS-1}.
The SGS projection method exhibits robust energy stability regardless of the damping parameter, while the conventional GS scheme fails to preserve energy dissipation when $\alpha=0$. For both algorithms, stronger Gilbert damping (\(\alpha\)) speeds up the decay of magnetic energy.

\begin{figure}[htbp]
    \centering
    \subfloat[$\alpha=0$]{\includegraphics[width=0.4\linewidth]{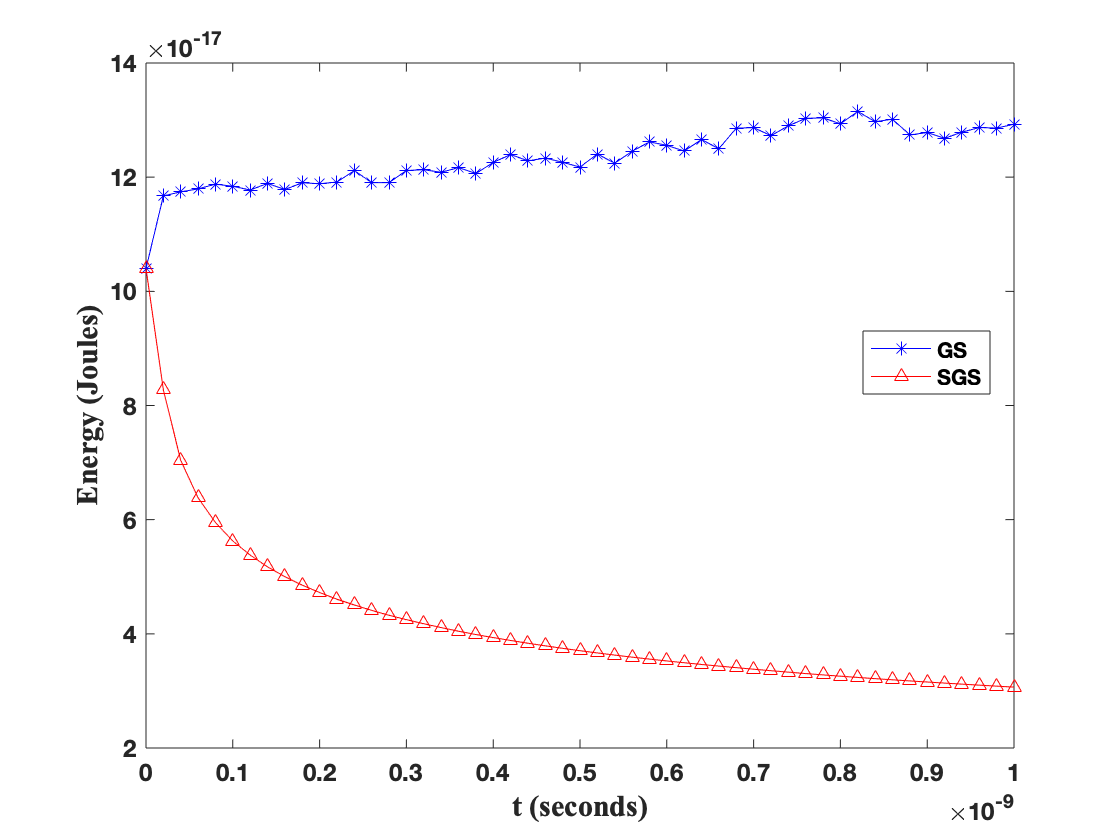}}
    \subfloat[$\alpha=0.01$]{\includegraphics[width=0.4\linewidth]{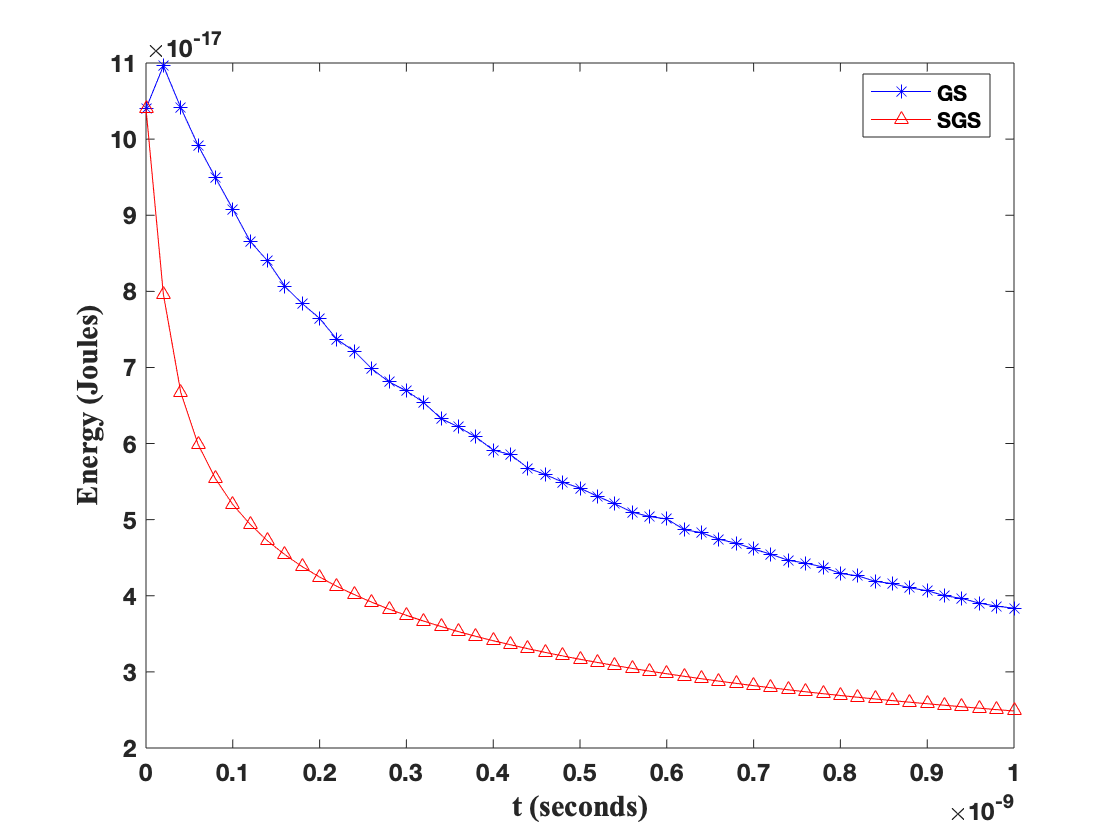}}
    \hspace{0.1in}
     \subfloat[$\alpha=0.1$]{\includegraphics[width=0.4\linewidth]{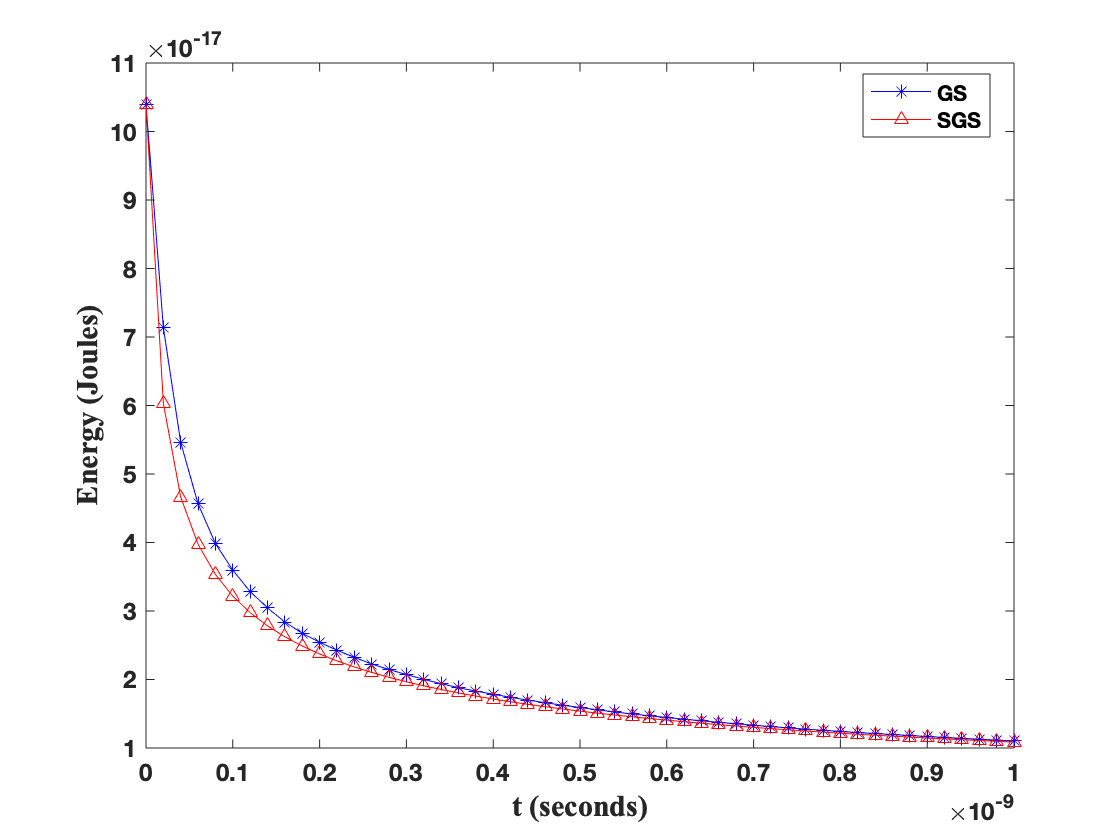}}
      \subfloat[$\alpha=0.5$]{\includegraphics[width=0.4\linewidth]{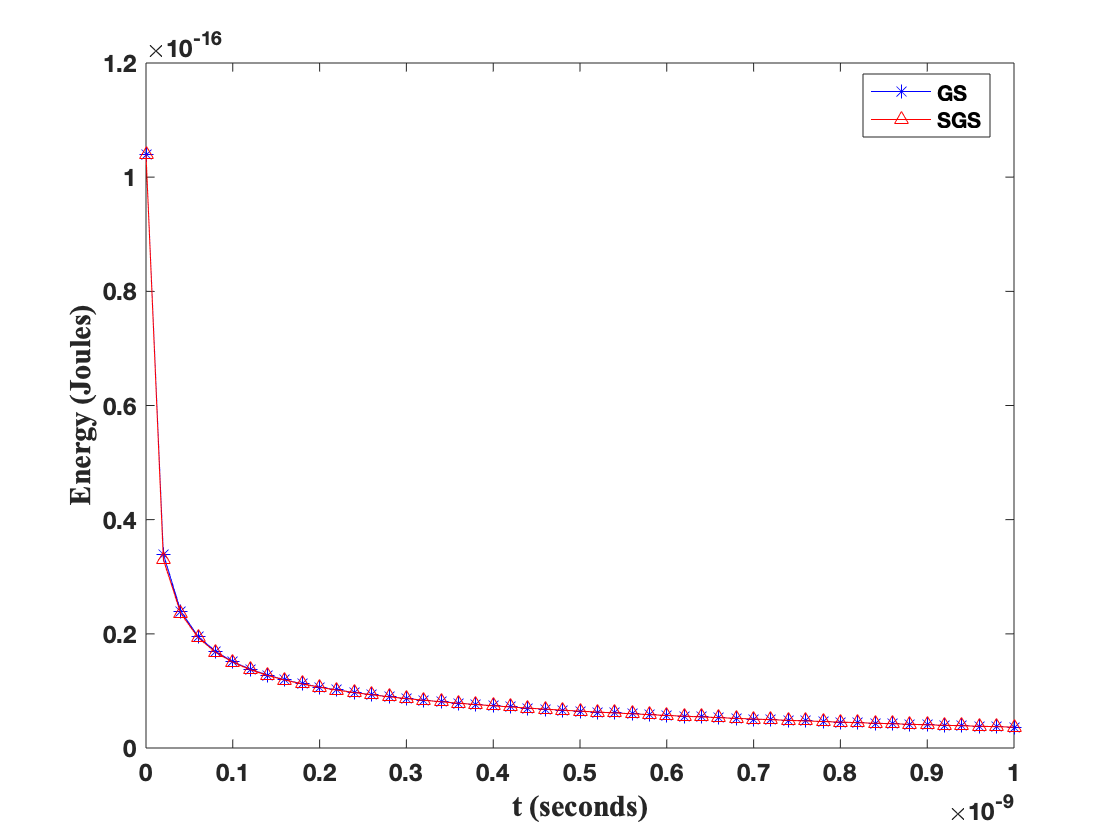}}
    \caption{Comparison of the energy evolution between GS and SGS projection method with $\alpha=0,0.01,0.1,0.5$ for the only exchange case given the initial S state.}
    \label{fig:energy-GS-SGS-exchange-v1}
\end{figure}

To see the performance with different damping parameters of $\alpha=0,0.01,0.1,0.5$ in details, we get the results in \Cref{fig:energy-GS-SGS-exchange-v1}. Across all tested damping values \(\alpha = 0, 0.01, 0.1, 0.5\) under the exchange-only model with an S-state initial condition, the SGS projection method outperforms the standard GS scheme in physical energy stability. When \(\alpha=0\), GS produces unphysical rising and fluctuating magnetic energy that violates energy dissipation, whereas SGS maintains strict monotonic energy decay. For all positive \(\alpha\), both GS and SGS yield dissipative energy curves, yet SGS always achieves a faster energy decay rate and lower residual energy than GS at the same damping coefficient. As the damping parameter \(\alpha\) increases, the gap in energy dissipation speed between GS and SGS gradually narrows, but SGS consistently retains superior energy-damping performance throughout the entire simulation time window.

\subsection{The simulations with only exchange field, anisotropy field, stray field}

\begin{figure}[htbp]
    \centering
    \subfloat[Initial S state, arrow profile]{\includegraphics[width=0.25\linewidth]{arrow_initial_S_v1.png}}
    \subfloat[arrow profile, $\alpha=0$]{\includegraphics[width=0.25\linewidth]{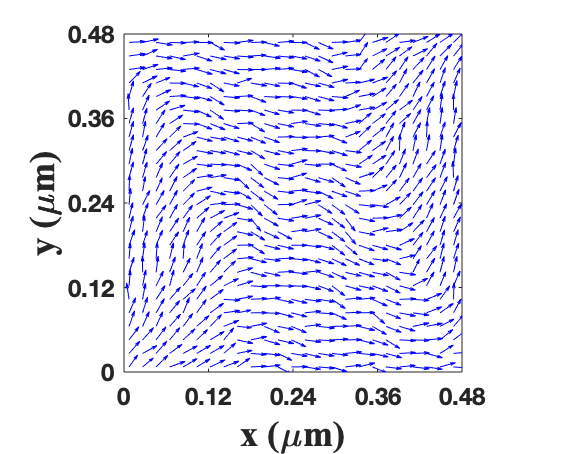}}
    \subfloat[arrow profile, $\alpha=0.005$]{\includegraphics[width=0.25\linewidth]{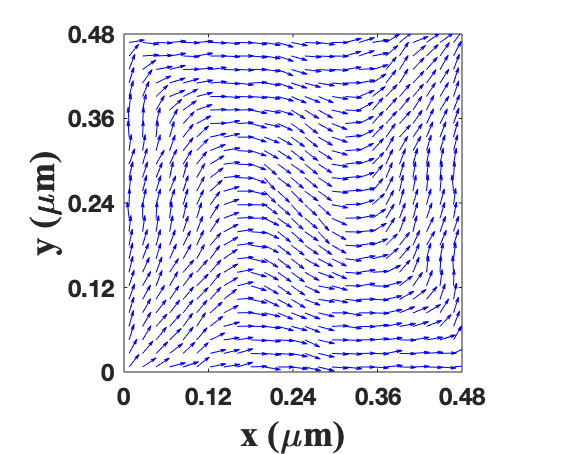}}
    \subfloat[arrow profile, $\alpha=0.01$]{\includegraphics[width=0.25\linewidth]{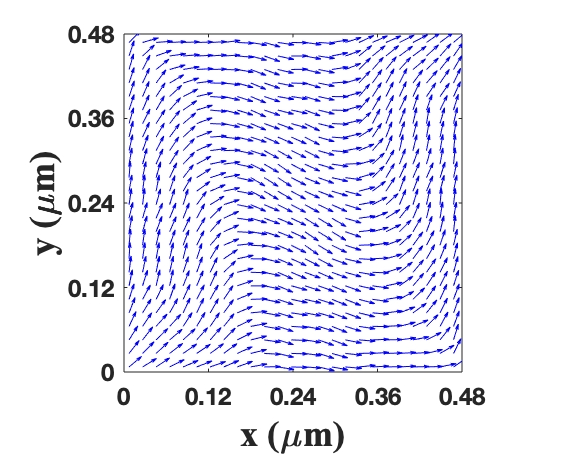}}
    \hspace{0.1in}
     \subfloat[Initial S state, color profile]{\includegraphics[width=0.25\linewidth]{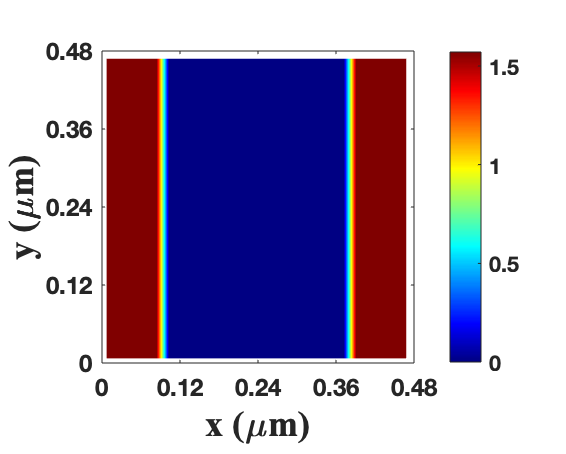}}
     \subfloat[color profile, $\alpha=0$]{\includegraphics[width=0.25\linewidth]{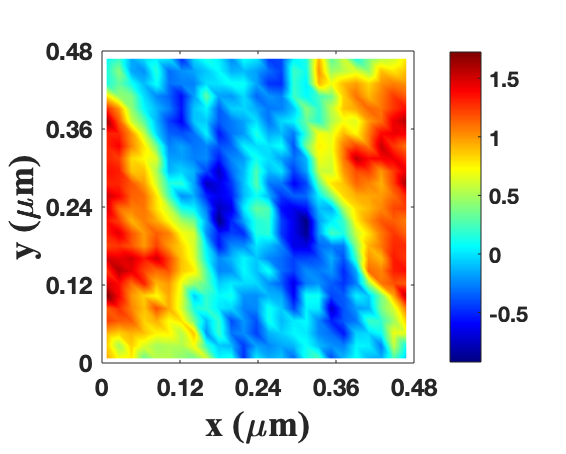}}
     \subfloat[arrow profile, $\alpha=0.005$]{\includegraphics[width=0.25\linewidth]{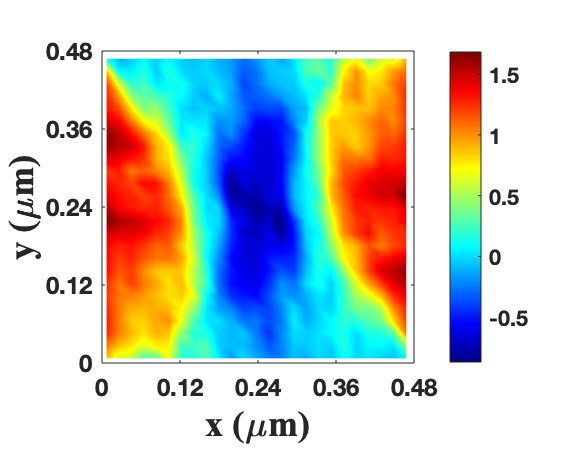}}
      \subfloat[arrow profile, $\alpha=0.01$]{\includegraphics[width=0.25\linewidth]{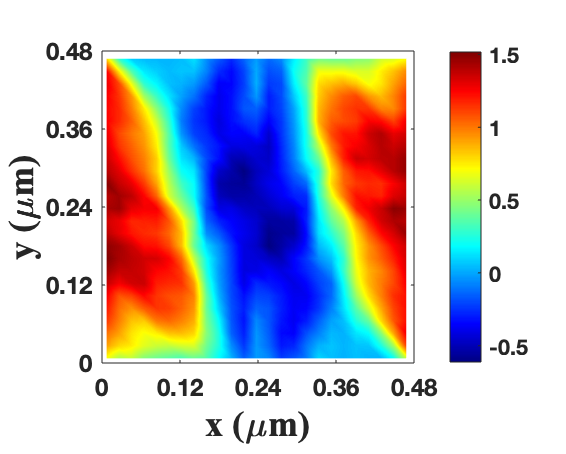}}
    \caption{GSPM for the simulation of the full model with exchange field, stray field, anisotropy field  up to the final time $T=1\;ns$ with the time step $\Delta t=1\;ps$.}
    \label{fig:GS-1}
\end{figure}

\begin{figure}[htbp]
    \centering
    \subfloat[Initial S state, arrow profile]{\includegraphics[width=0.25\linewidth]{arrow_initial_S_v1.png}}
    \subfloat[arrow profile, $\alpha=0$]{\includegraphics[width=0.25\linewidth]{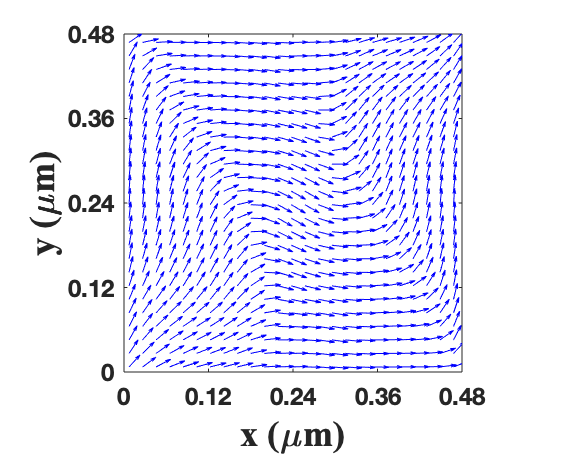}}
    \subfloat[arrow profile, $\alpha=0.005$]{\includegraphics[width=0.25\linewidth]{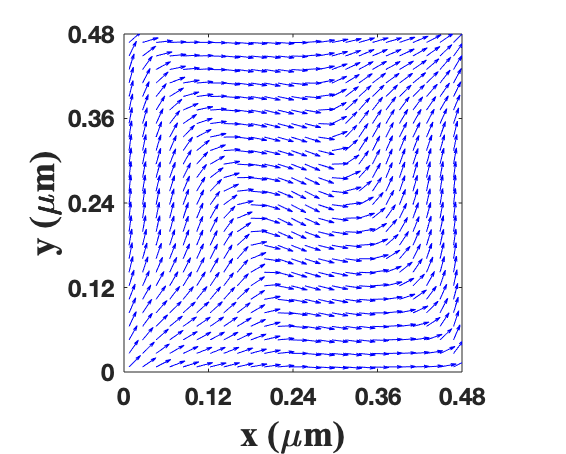}}
    \subfloat[arrow profile, $\alpha=0.01$]{\includegraphics[width=0.25\linewidth]{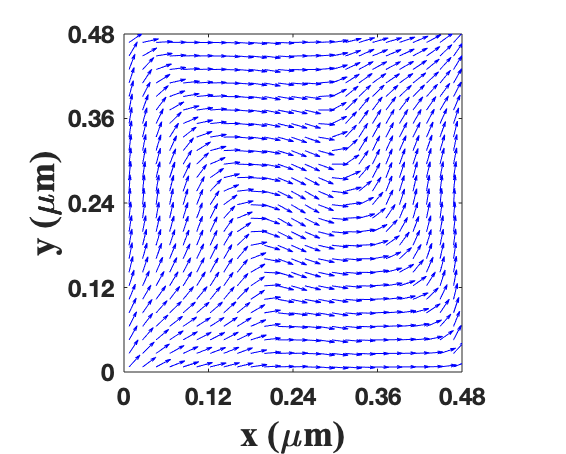}}
    \hspace{0.1in}
     \subfloat[Initial S state, color profile]{\includegraphics[width=0.25\linewidth]{color_initial_S_v1.png}}
     \subfloat[color profile, $\alpha=0$]{\includegraphics[width=0.25\linewidth]{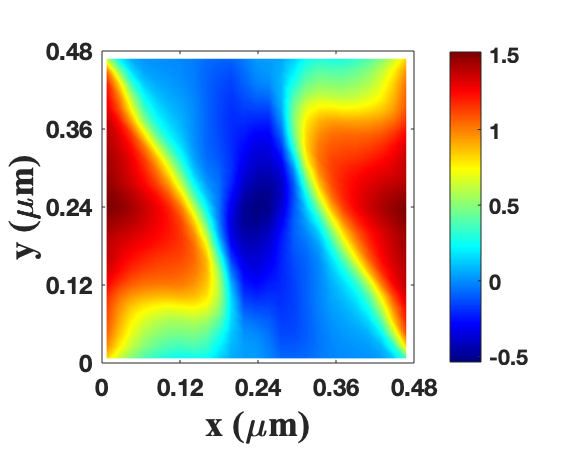}}
     \subfloat[arrow profile, $\alpha=0.005$]{\includegraphics[width=0.25\linewidth]{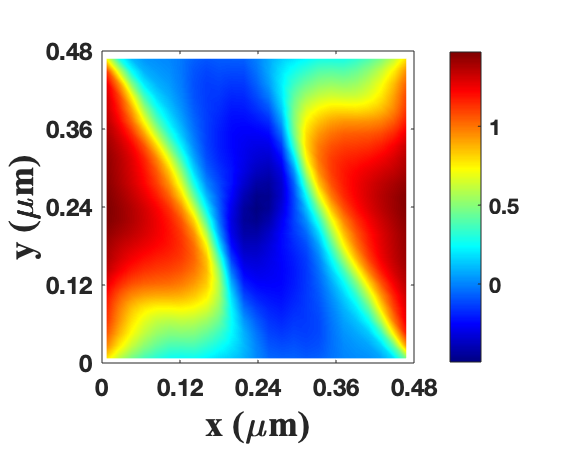}}
      \subfloat[arrow profile, $\alpha=0.01$]{\includegraphics[width=0.25\linewidth]{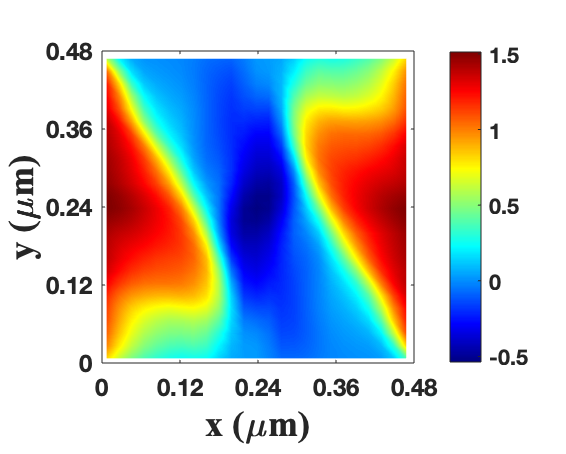}}
    \caption{SGSPM for the simulation of the full model with exchange field, stray field, anisotropy field  up to the final time $T=1\;ns$ with the time step $\Delta t=1\;ps$.}
    \label{fig:SGS-1}
\end{figure}

In this section, we test the full Landau-Lifshitz equation with the effective field below,
\begin{align*}
    \h_{\textbf{eff}}=\epsilon\Delta \m-Q(m_2\e_2+m_3\e_3)+\h_s+\h_e.
\end{align*}
The simulation model uses a ferromagnetic thin film with dimensions $480\times 480 \times 20\;nm^3$ discretized on a $100\times 100\times 4$ grid. The temporal step size we take is $\Delta t=1\;ps$. Given initial S state, the magnetization profile using GSPM with $\alpha=0,0.005,0.01$ is presented in \Cref{fig:GS-1}, while the profile using SGSPM with $\alpha=0,0.005,0.01$ is presented in \Cref{fig:SGS-1}. It suggests that for the full micromagnetic model incorporating exchange, stray and anisotropy fields, both GSPM and SGSPM yield visually similar steady-state magnetization arrow distributions and color contour maps at the final time \(T=1\ \text{ns}\) under \(\alpha=0, 0.005, 0.01\), which demonstrates that the two methods can produce consistent static magnetization patterns after long-time evolution. However, SGSPM maintains unconditional monotonic energy decay for all damping coefficients including \(\alpha=0\) during the whole simulation process, complying with the physical energy dissipation law strictly. In contrast, GSPM generates unphysical energy growth when \(\alpha=0\), lacking inherent energy stability. Overall, SGSPM achieves identical steady magnetization configurations as GSPM while providing physically valid energy evolution at all tested damping parameters, making it a more robust numerical scheme for micromagnetic simulations.

\begin{figure}[htbp]
    \centering
    \subfloat[GS, $\alpha=0$]{\includegraphics[width=0.4\linewidth]{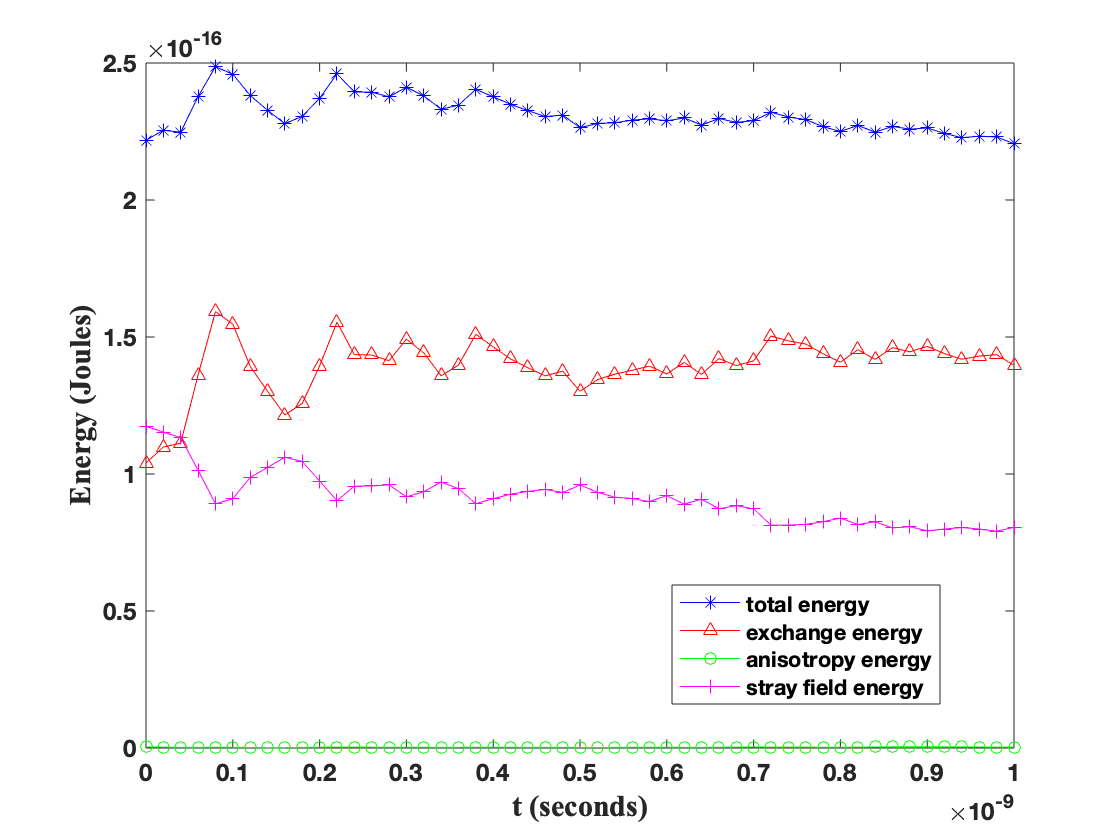}}
    \subfloat[GS, $\alpha=0.01$]{\includegraphics[width=0.4\linewidth]{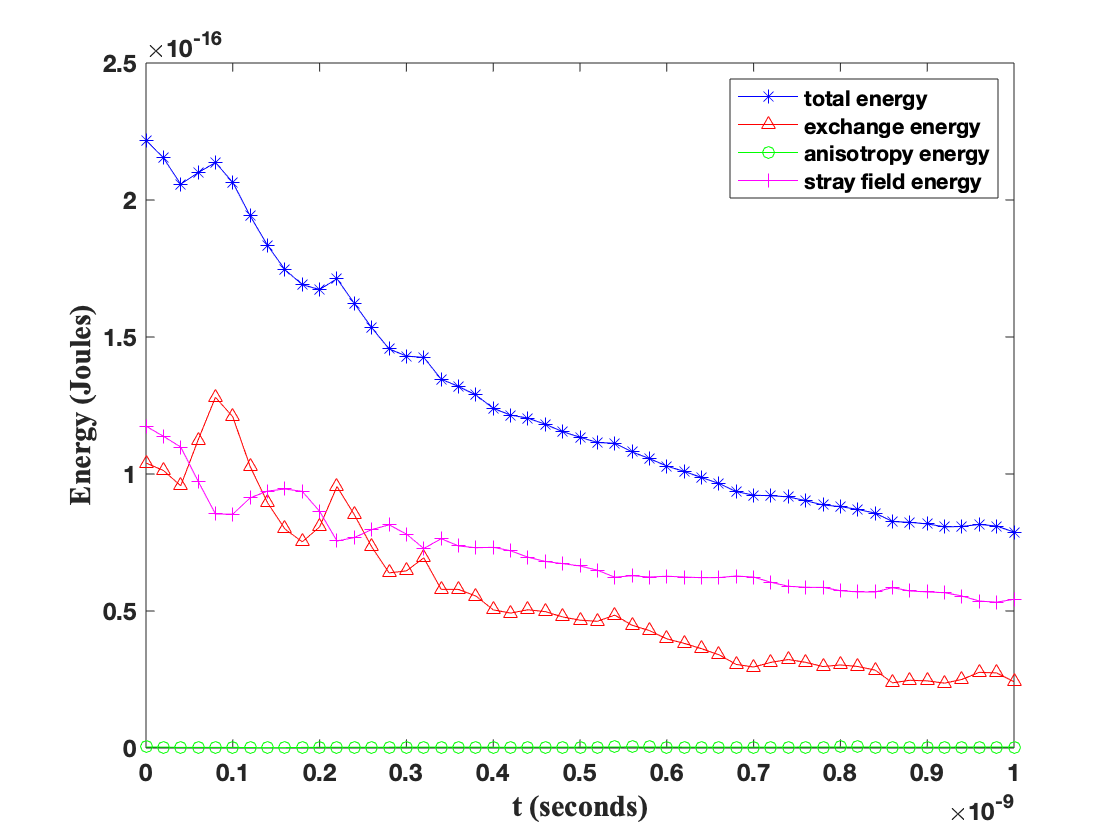}}
    \hspace{0.1in}
     \subfloat[GS, $\alpha=0.1$]{\includegraphics[width=0.4\linewidth]{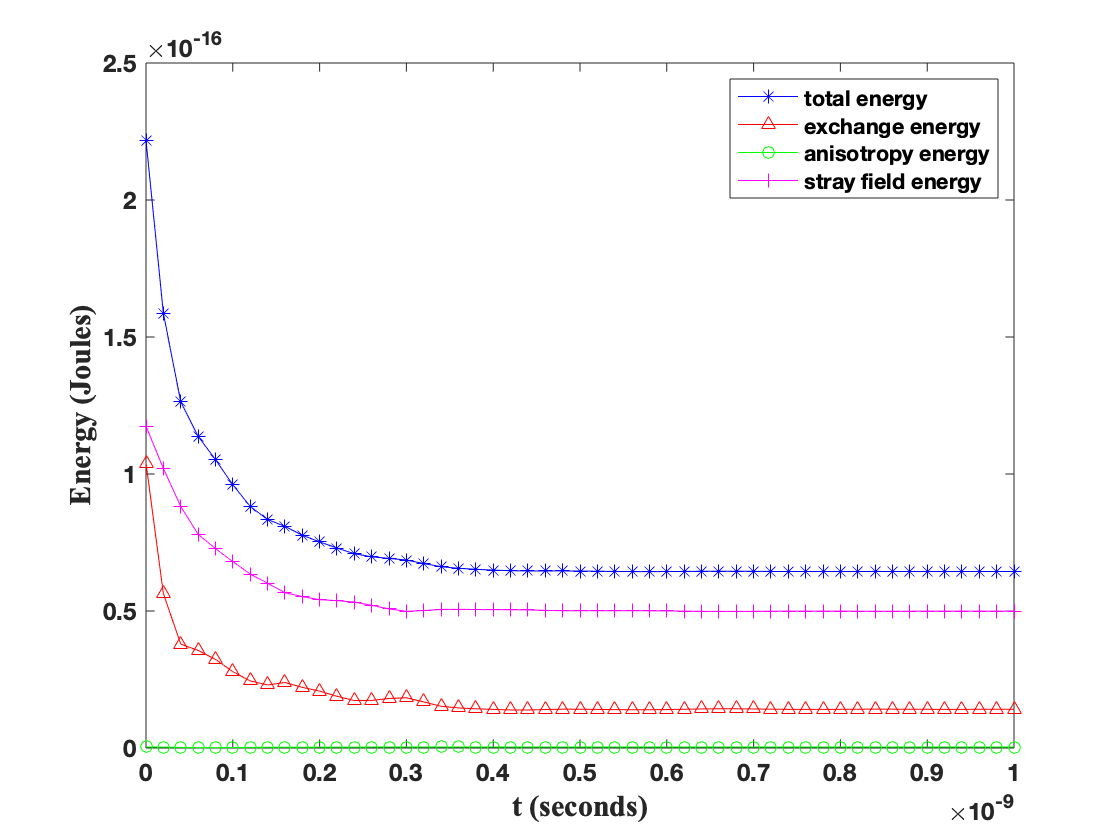}}
      \subfloat[GS, $\alpha=0.5$]{\includegraphics[width=0.4\linewidth]{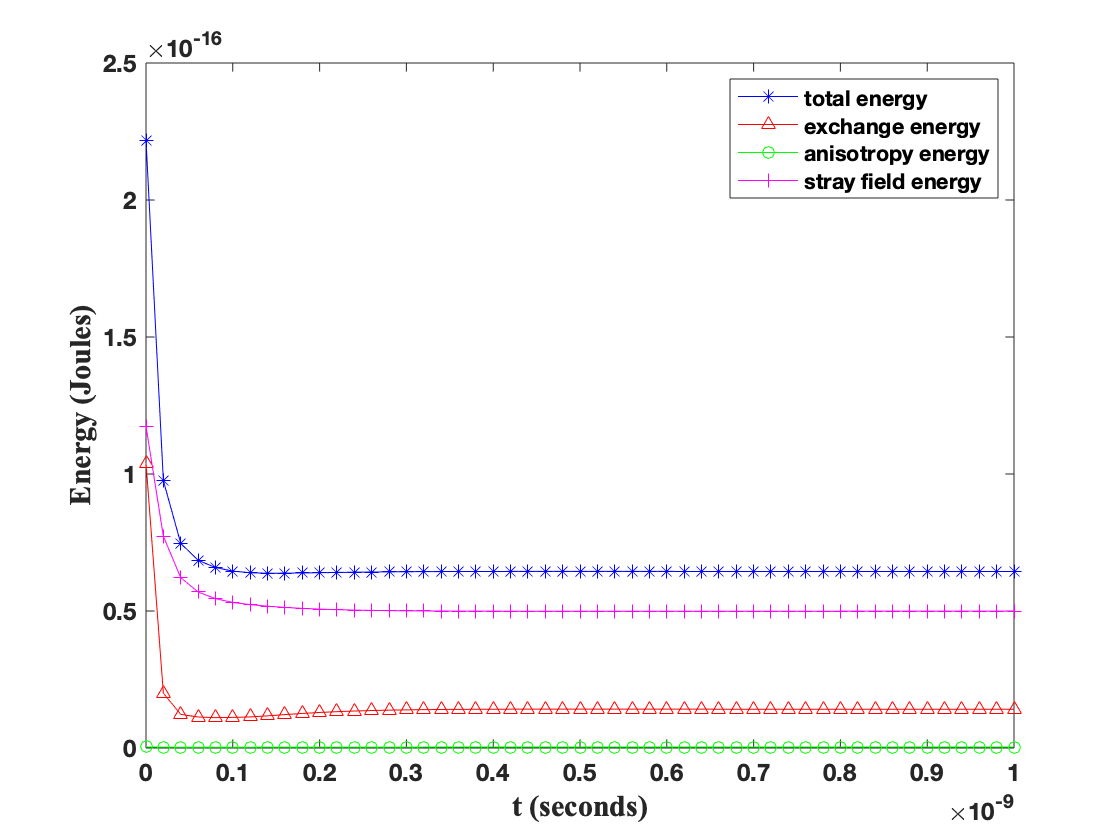}}
    \caption{The energy evolution for the total energy, the exchange energy, the anisotropy energy and stray field energy by GS projection method with $\alpha=0,0.01,0.1,0.5$.}
    \label{fig:energy-GS-v1}
\end{figure}

\begin{figure}[htbp]
    \centering
    \subfloat[SGS, $\alpha=0$]{\includegraphics[width=0.4\linewidth]{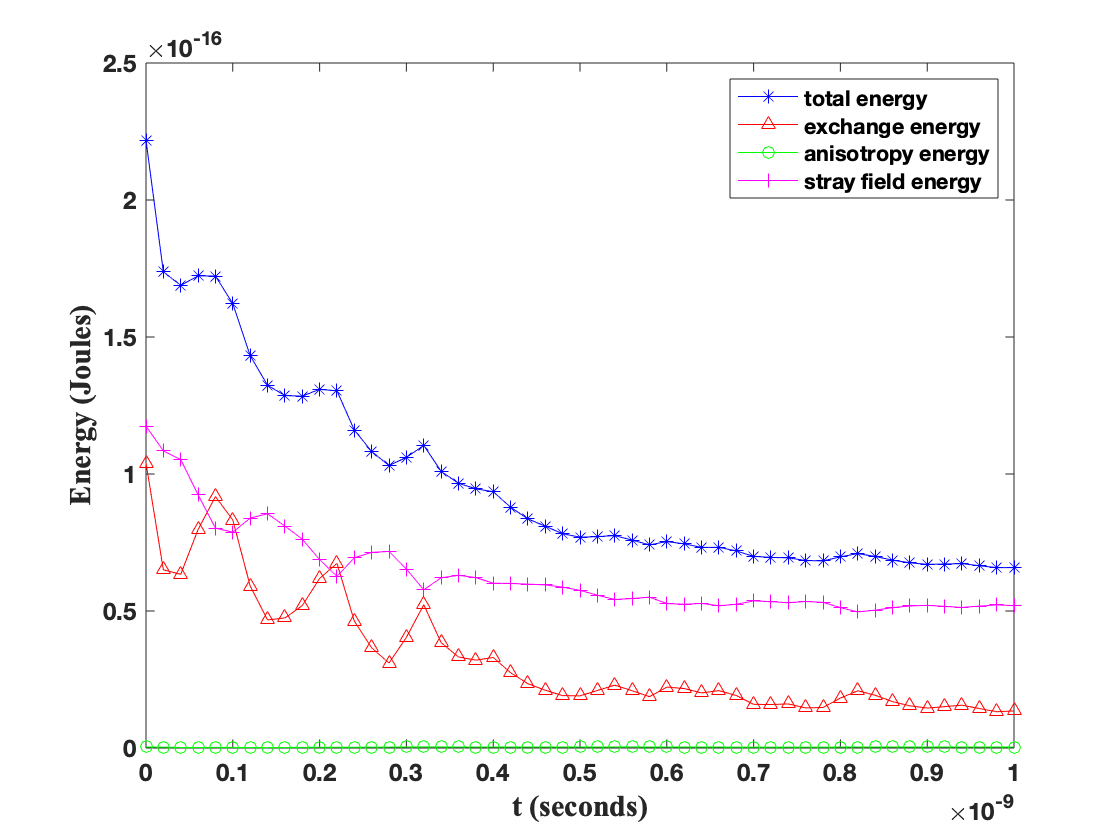}}
    \subfloat[SGS, $\alpha=0.01$]{\includegraphics[width=0.4\linewidth]{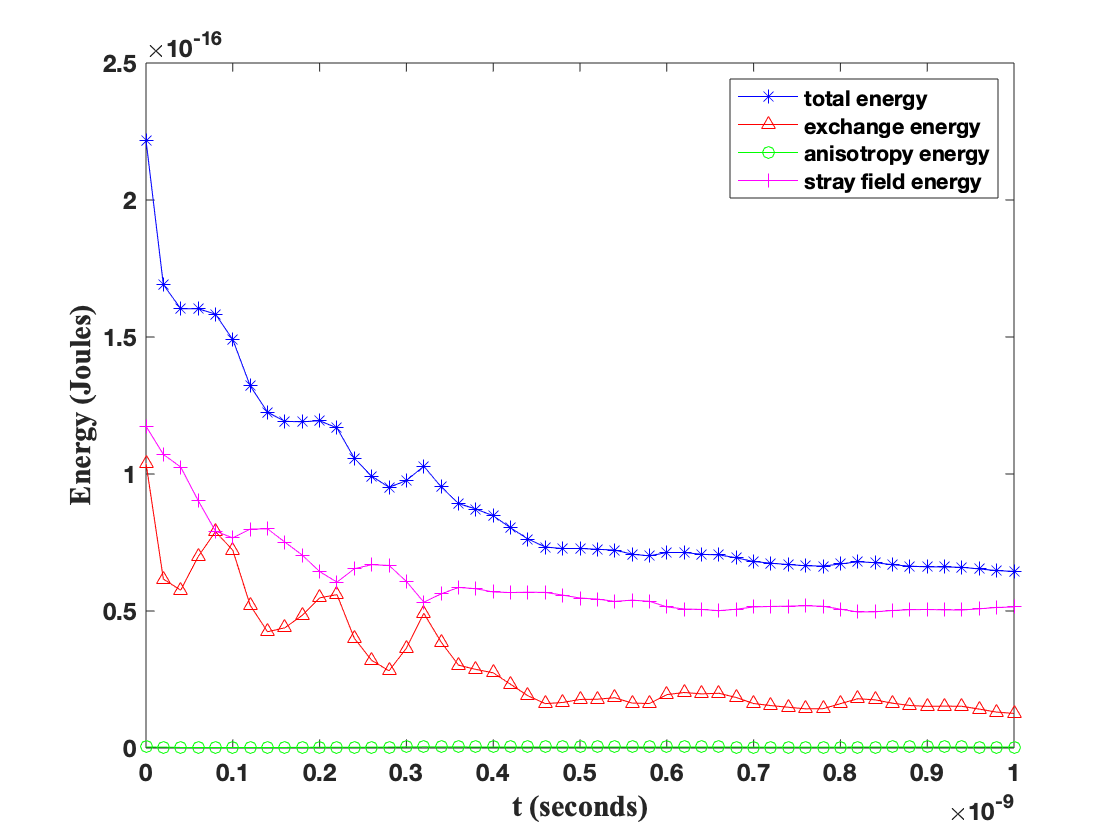}}
    \hspace{0.1in}
     \subfloat[SGS, $\alpha=0.1$]{\includegraphics[width=0.4\linewidth]{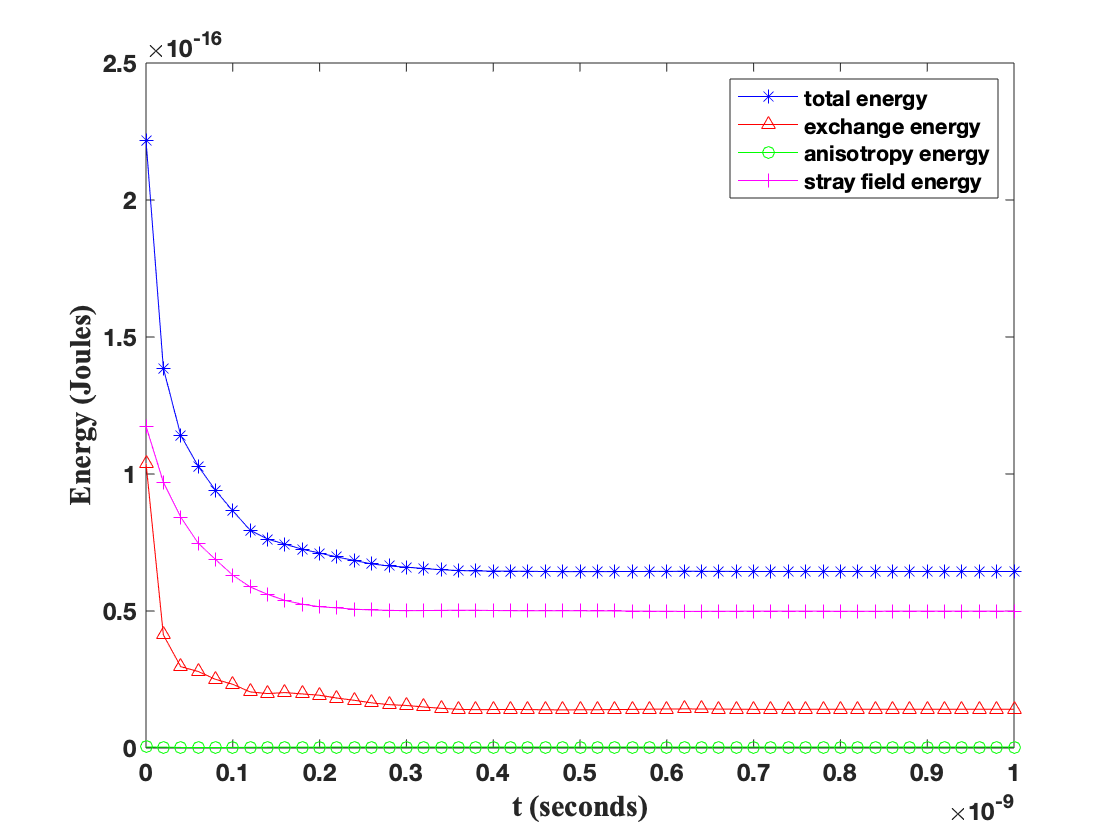}}
      \subfloat[SGS, $\alpha=0.5$]{\includegraphics[width=0.4\linewidth]{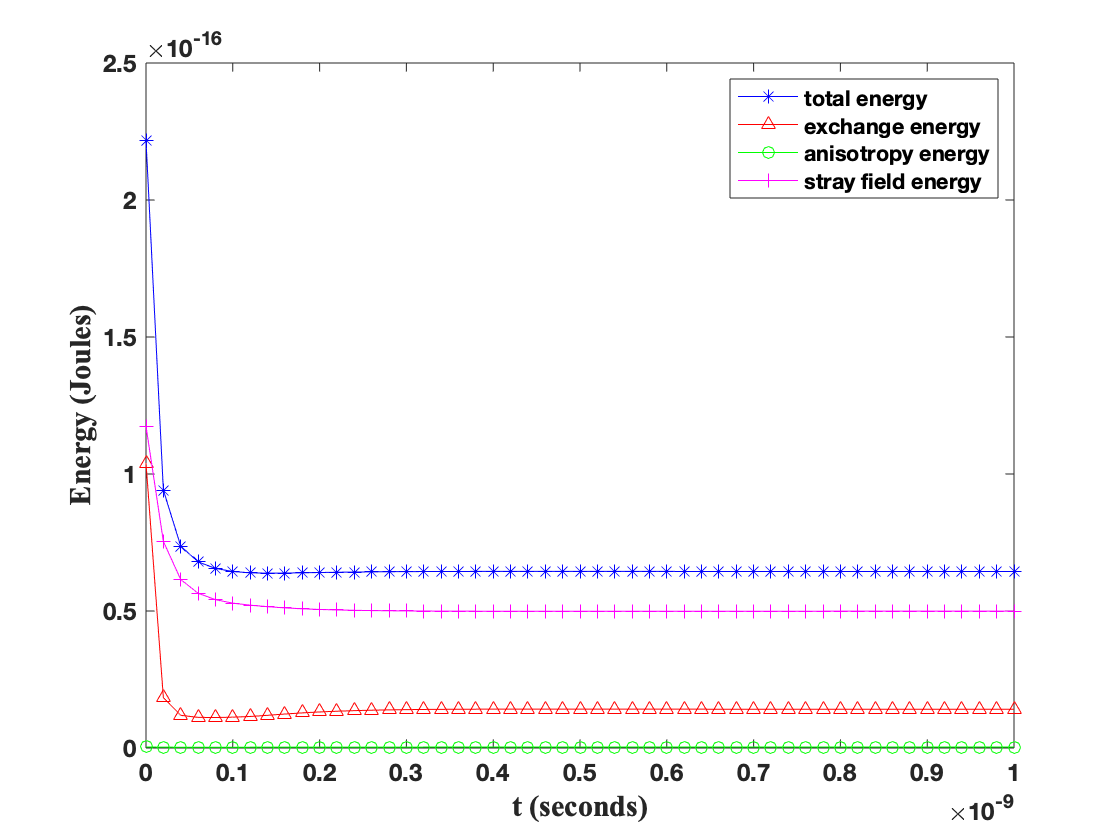}}
    \caption{The energy evolution for the total energy, the exchange energy, the anisotropy energy and stray field energy by SGS projection method with $\alpha=0,0.01,0.1,0.5$.}
    \label{fig:energy-SGS-v1}
\end{figure}

The results of energy evolution by GSPM method and SGSPM method for the equilibrium at final time $T=1\;ns$ are presented in \Cref{fig:energy-GS-v1} and \Cref{fig:energy-SGS-v1}, respectively. It follows from \Cref{fig:energy-GS-v1} that for the GS projection method, the decomposition curves of total, exchange, anisotropy and stray field energies show distinct behaviours dependent on the damping coefficient \(\alpha\). At \(\alpha=0\), the total energy fluctuates and rises above its initial value, which violates the physical energy dissipation rule; both exchange and stray field energies also present unphysical oscillatory growth, while anisotropy energy remains nearly zero throughout the simulation. For all positive damping values (\(\alpha=0.01,0.1,0.5\)), the total energy decays monotonically, accompanied by continuous decreases in exchange and stray field energies, and anisotropy energy stays negligible. Larger \(\alpha\) accelerates the decay rate of total, exchange and stray field energies, and yields lower saturated residual energy at the end of the time window. The unphysical energy growth observed under \(\alpha=0\) demonstrates that the GS scheme cannot guarantee the dissipation of each energy component or the overall total energy without Gilbert damping, revealing its fundamental energy stability defect for undamped micromagnetic simulations. It follows from \Cref{fig:energy-SGS-v1} that for the SGS projection method, the total energy maintains strictly monotonic decay for all damping coefficients \(\alpha = 0, 0.01, 0.1, 0.5\), fully satisfying the physical energy dissipation principle. The exchange energy and stray field energy gradually decline over the whole time span, while anisotropy energy stays close to zero in all cases.
When \(\alpha=0\), exchange and stray field energies show mild oscillations in the early stage but still follow an overall downward trend without abnormal energy rise; the total energy decreases steadily, which is physically valid, unlike the unphysical energy growth of GS under zero damping.
For positive \(\alpha\) values (\(0.01, 0.1, 0.5\)), larger damping parameters produce faster decay rates and lower residual levels for total, exchange and stray field energies, with nearly no visible oscillations in energy curves.
Overall, SGS guarantees physically reasonable dissipative evolution of all decomposed energy components in both undamped and damped full micromagnetic simulations, possessing intrinsic energy stability absent from the standard GS scheme.

\begin{figure}[htbp]
    \centering
    \subfloat[$\alpha=0$]{\includegraphics[width=0.4\linewidth]{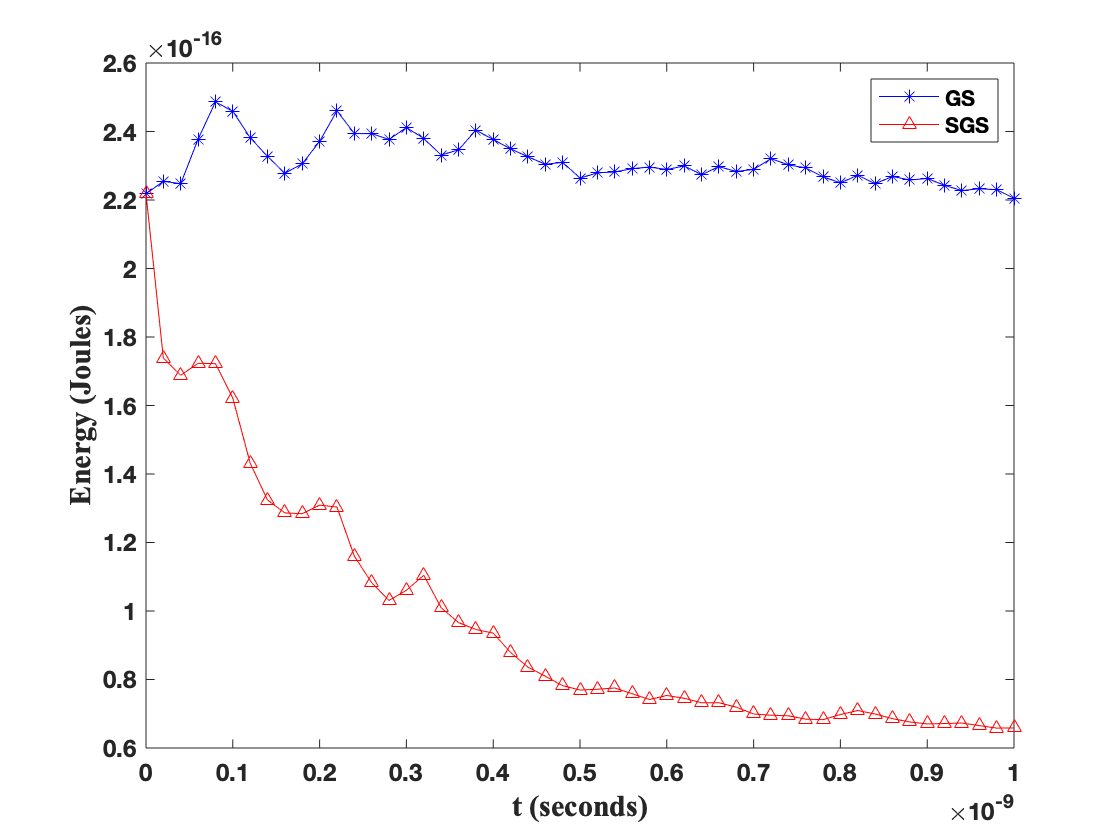}}
    \subfloat[$\alpha=0.01$]{\includegraphics[width=0.4\linewidth]{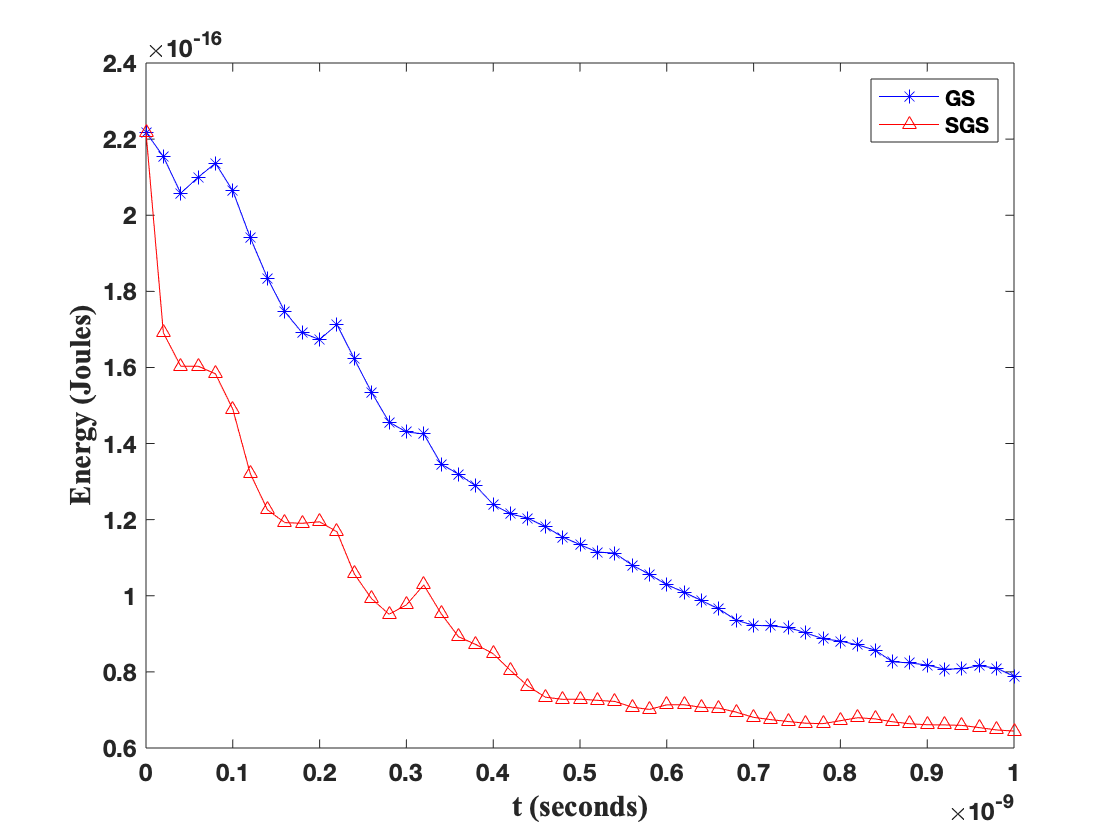}}
    \hspace{0.1in}
     \subfloat[$\alpha=0.1$]{\includegraphics[width=0.4\linewidth]{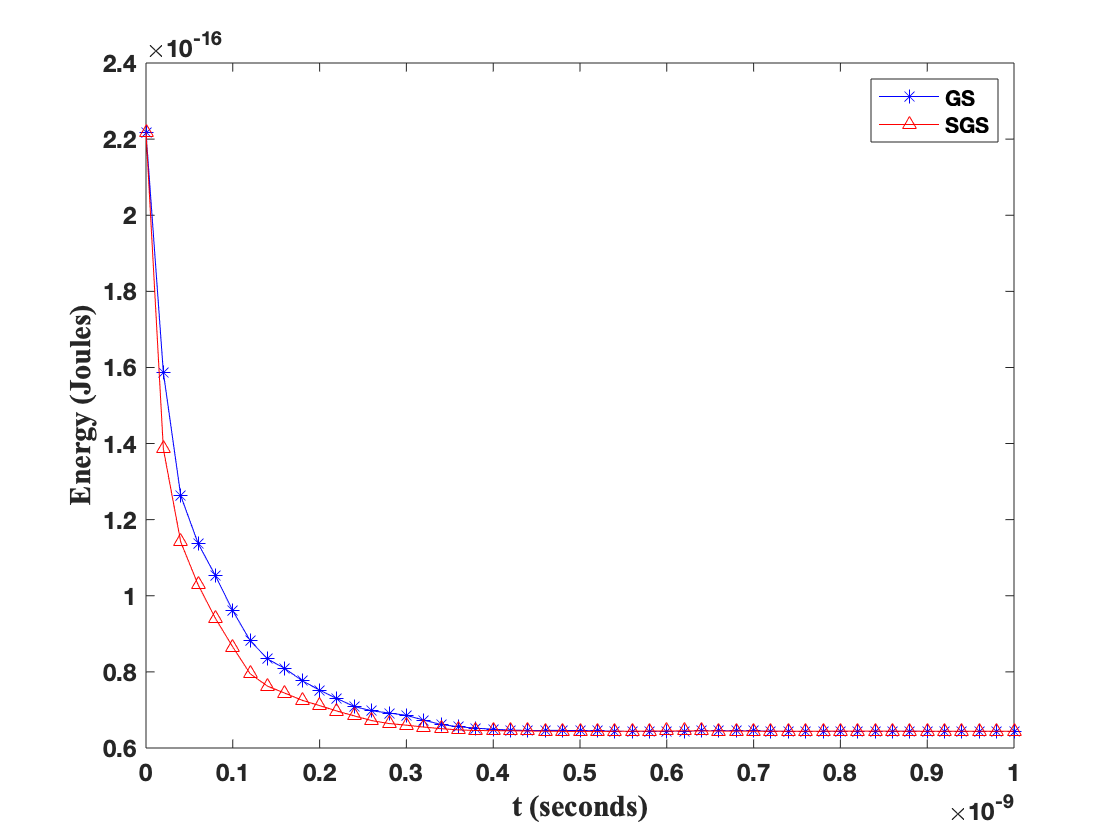}}
      \subfloat[$\alpha=0.5$]{\includegraphics[width=0.4\linewidth]{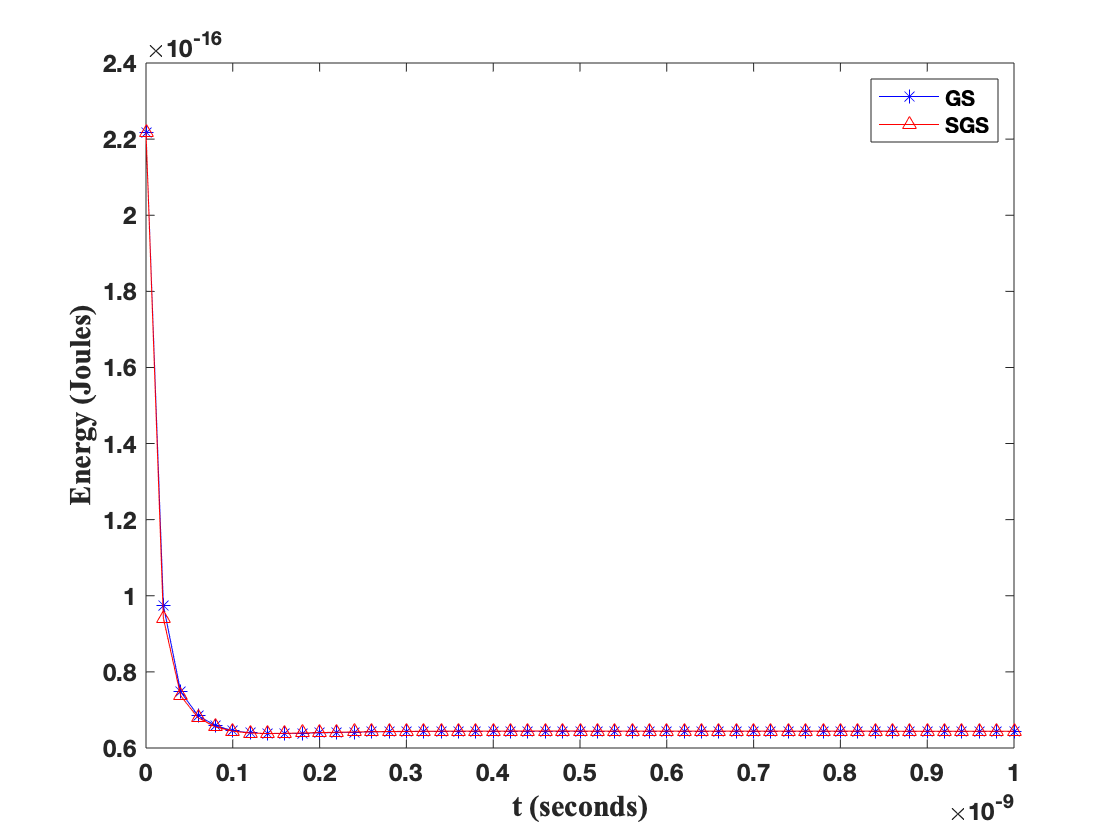}}
    \caption{Comparison of the energy evolution between GS and SGS projection method with $\alpha=0,0.01,0.1,0.5$ up to the final time $T=1\;ns$ with the time step size $\Delta t=1\;ps$.}
    \label{fig:energy-GS-SGS-v1}
\end{figure}

For more details, we get the results presented in \Cref{fig:energy-GS-SGS-v1}. Across all tested damping coefficients \(\alpha = 0, 0.01, 0.1, 0.5\) under the full micromagnetic model with \(T=1\ \text{ns}\) and \(\Delta t=1\ \text{ps}\), the SGS projection method consistently outperforms the GS scheme in energy stability and physical consistency. At \(\alpha=0\), GS exhibits unphysical fluctuating total energy that sits above its initial value and violates energy dissipation, while SGS delivers steady monotonic energy decay. For positive \(\alpha\) values, both schemes produce dissipative energy curves, yet SGS always achieves a faster energy decay rate and lower residual energy than GS at identical damping strength. As \(\alpha\) increases, the discrepancy in energy decay speed between GS and SGS diminishes, but SGS still maintains superior dissipative performance throughout the entire simulation interval. Unlike GS, SGS unconditionally satisfies the energy dissipation law for both undamped and damped cases, confirming its robust inherent energy stability for full micromagnetic simulations.

\begin{figure}[htbp]
    \centering
    \subfloat[GS, arrow profile]{\includegraphics[width=0.2\linewidth]{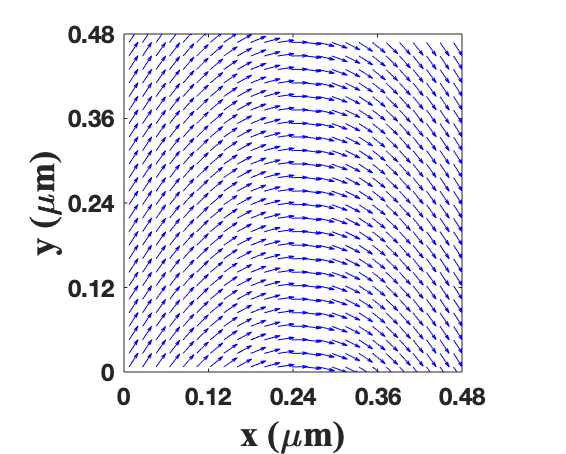}}
    \subfloat[GS, arrow profile, $\alpha=0$]{\includegraphics[width=0.2\linewidth]{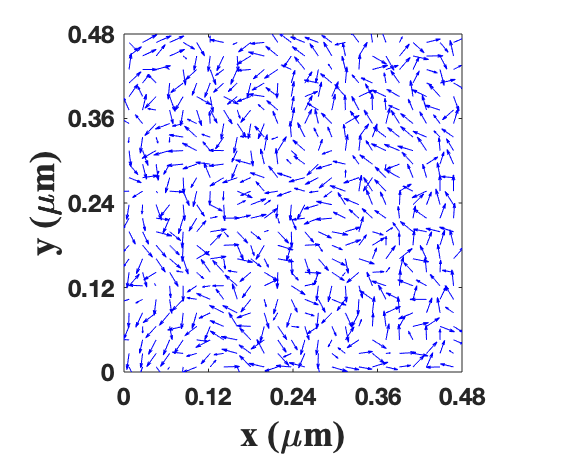}}
    \subfloat[GS, arrow profile, $\alpha=0.01$]{\includegraphics[width=0.2\linewidth]{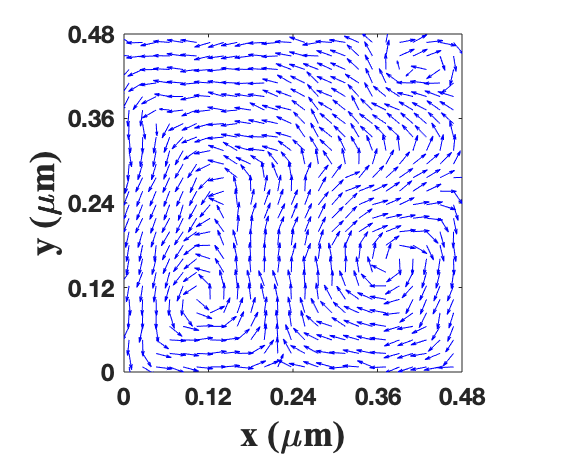}}
    \subfloat[GS, arrow profile, $\alpha=0.1$]{\includegraphics[width=0.2\linewidth]{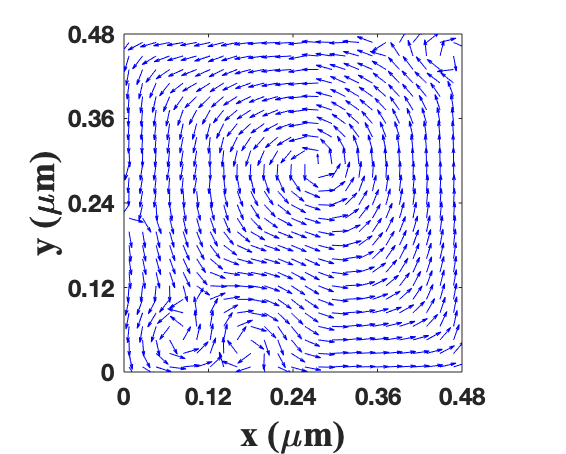}}
    \subfloat[GS, arrow profile, $\alpha=0.5$]{\includegraphics[width=0.2\linewidth]{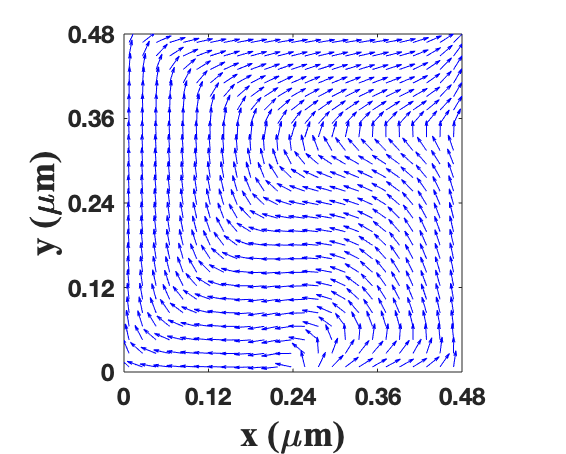}}
    \hspace{0.1in}
    \subfloat[SGS, arrow profile]{\includegraphics[width=0.2\linewidth]{arrow_SGS_all_alpha_0_initial_v1.png}}
    \subfloat[SGS, arrow profile, $\alpha=0$]{\includegraphics[width=0.2\linewidth]{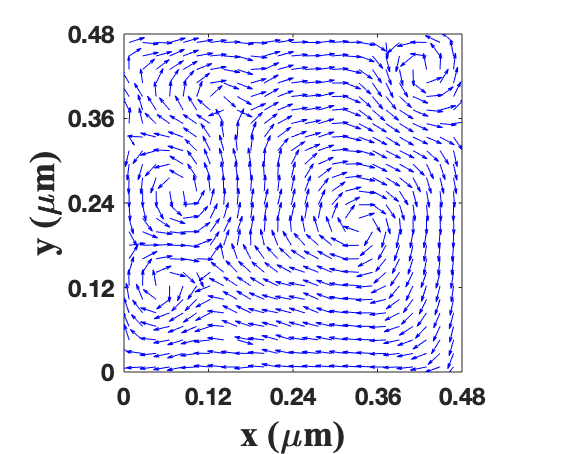}}
    \subfloat[SGS, arrow profile, $\alpha=0.01$]{\includegraphics[width=0.2\linewidth]{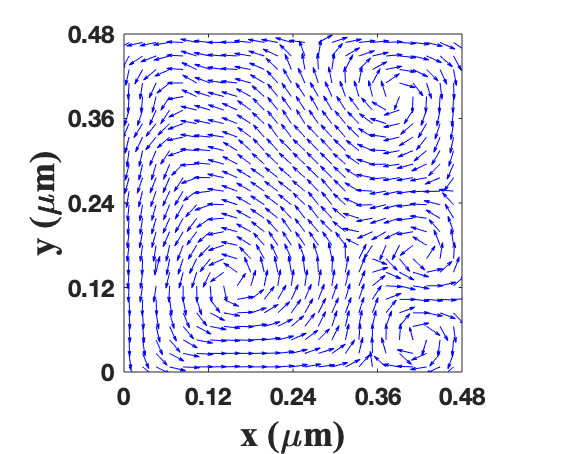}}
    \subfloat[SGS, arrow profile, $\alpha=0.1$]{\includegraphics[width=0.2\linewidth]{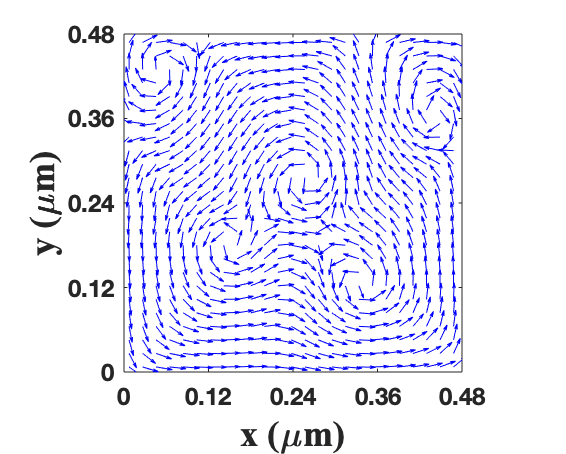}}
    \subfloat[SGS, arrow profile, $\alpha=0.5$]{\includegraphics[width=0.2\linewidth]{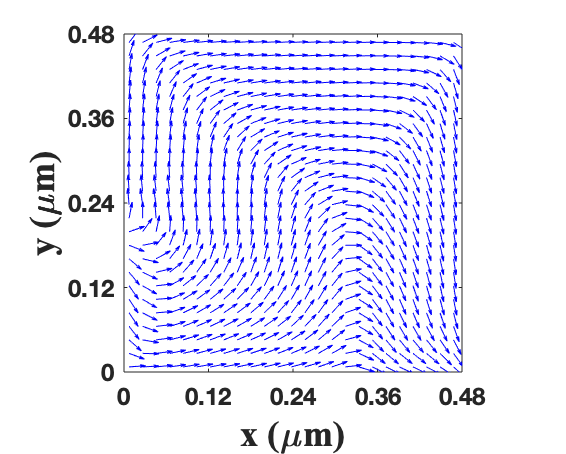}}
    \caption{GSPM and SGSPM for the simulation o for the full model with exchange field, stray field, anisotropy field with initial condition $\m_0=[\cos(\cos(\pi x))\sin(0.01),\sin(\cos(\pi x))\sin(0.01),\cos(0.01)]^T$. The final time $T=1\;ns$, the time step $\Delta t=1\;ps$. The top row is GS with $\alpha=0,0.01,0.1,0.5$, the bottom row is SGS with $\alpha=0,0.01,0.1,0.5$.}
    \label{fig:GS-SGS-full-1}
\end{figure}

To see the consistency of the full micromagnetic simulations, we test the case with initial condition as above as in \eqref{eq-m0-1}, we get the micromagnetic profile presented in \Cref{fig:GS-SGS-full-1}.

\begin{figure}[htbp]
    \centering
    \subfloat[GS, $\alpha=0$]{\includegraphics[width=0.25\linewidth]{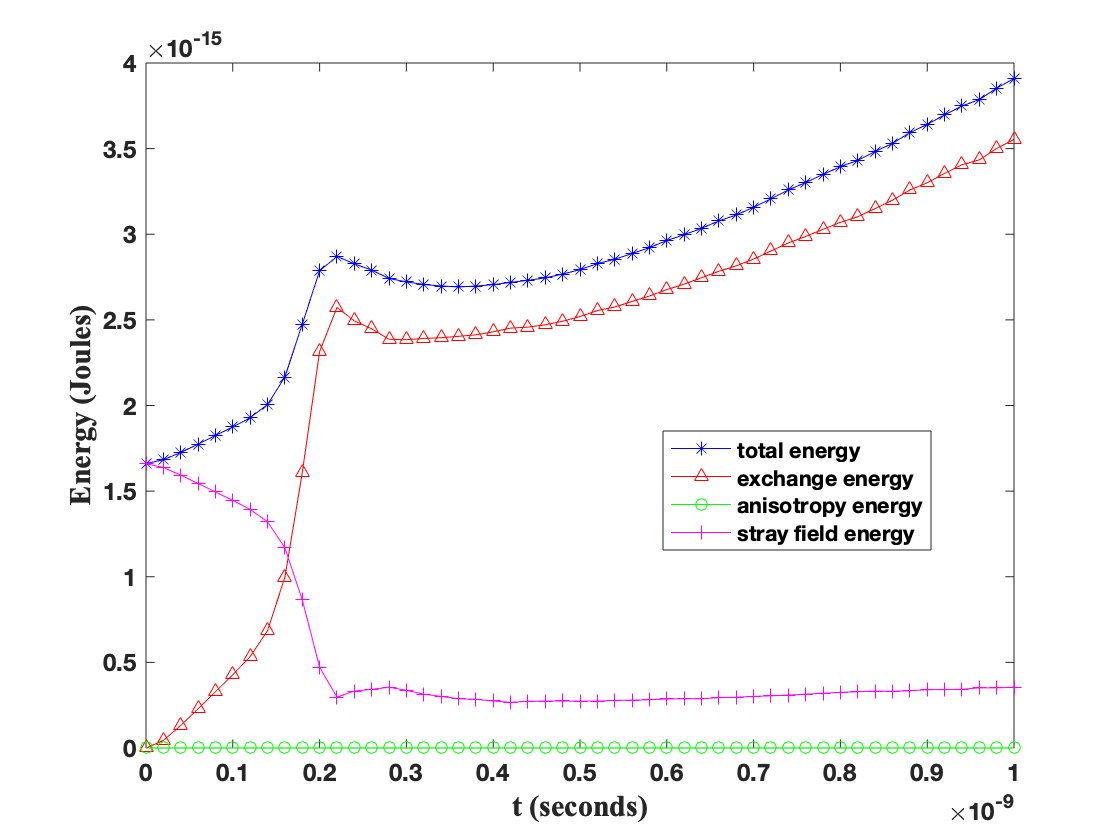}}
    \subfloat[GS, $\alpha=0.01$]{\includegraphics[width=0.25\linewidth]{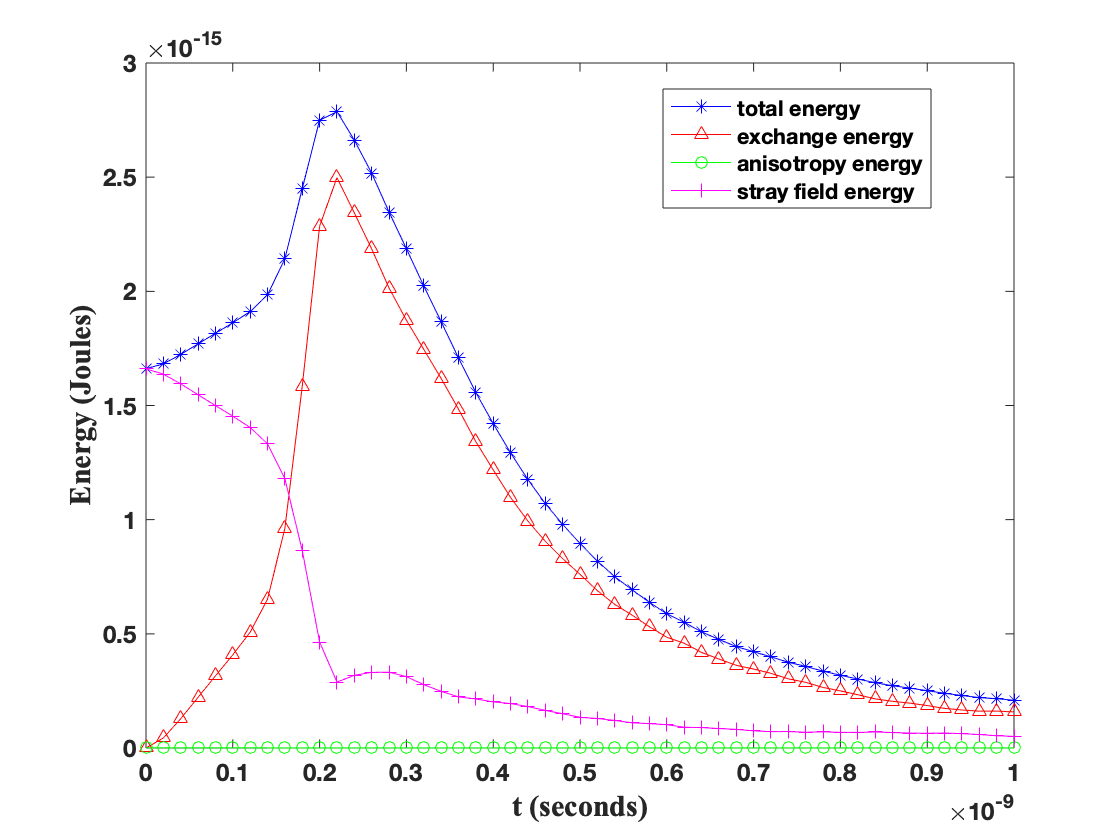}}
     \subfloat[GS, $\alpha=0.1$]{\includegraphics[width=0.25\linewidth]{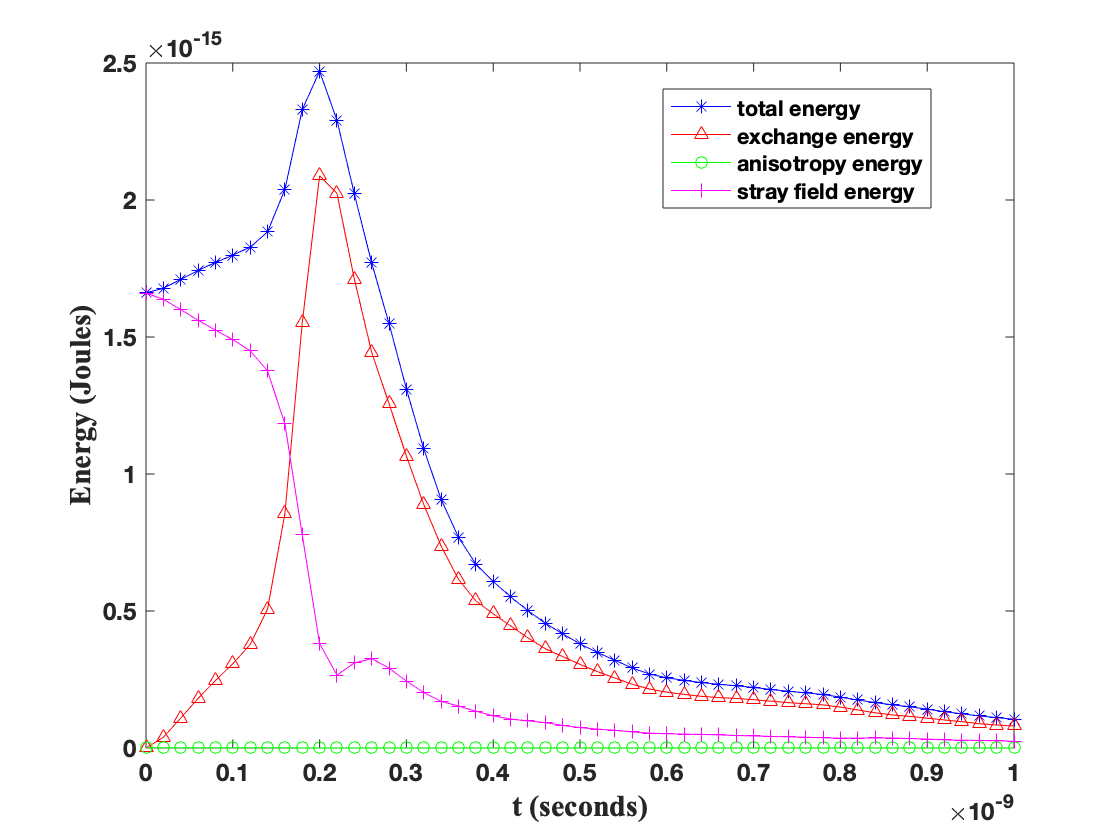}}
      \subfloat[GS, $\alpha=0.5$]{\includegraphics[width=0.25\linewidth]{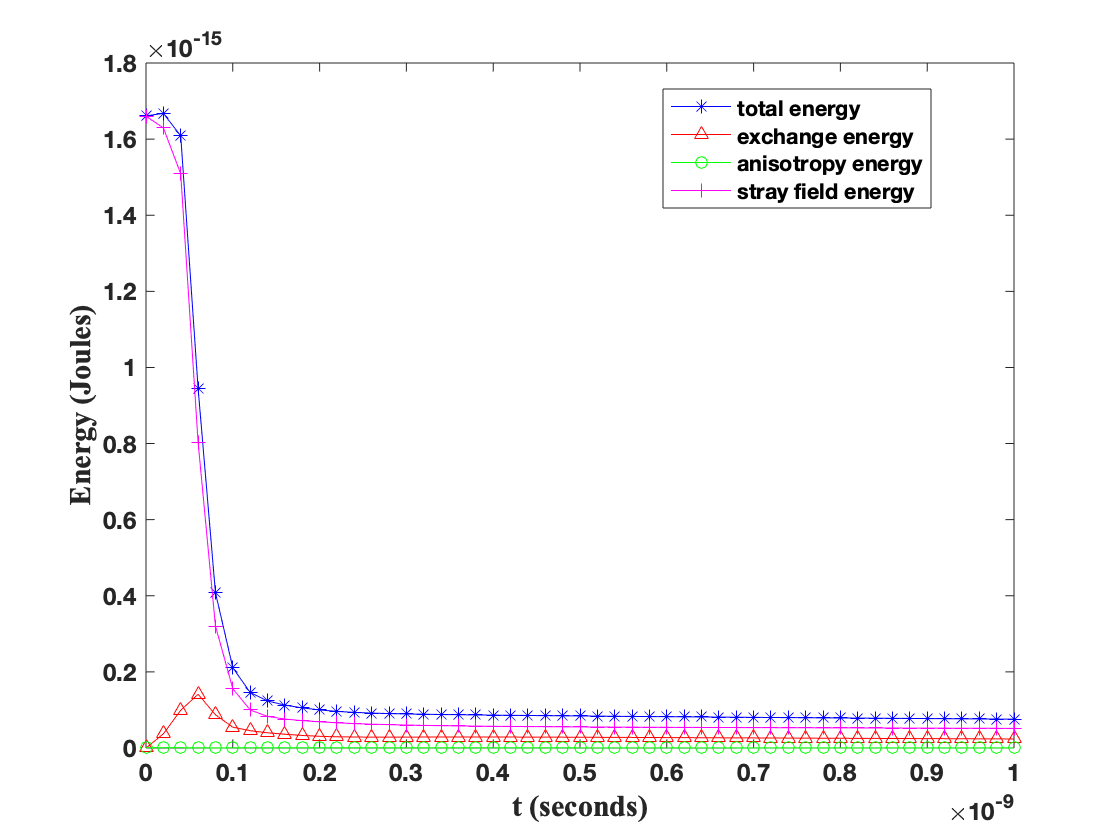}}
      \hspace{0.1in}
      \subfloat[SGS, $\alpha=0$]{\includegraphics[width=0.25\linewidth]{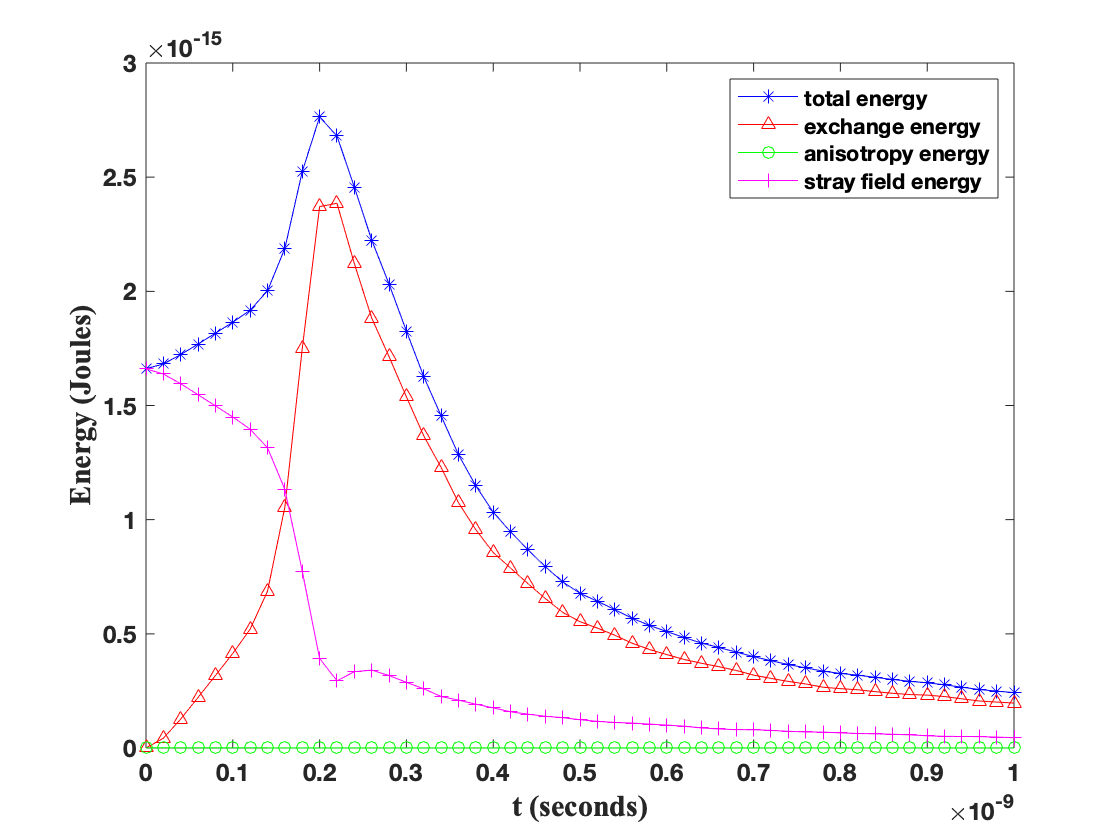}}
    \subfloat[SGS, $\alpha=0.01$]{\includegraphics[width=0.25\linewidth]{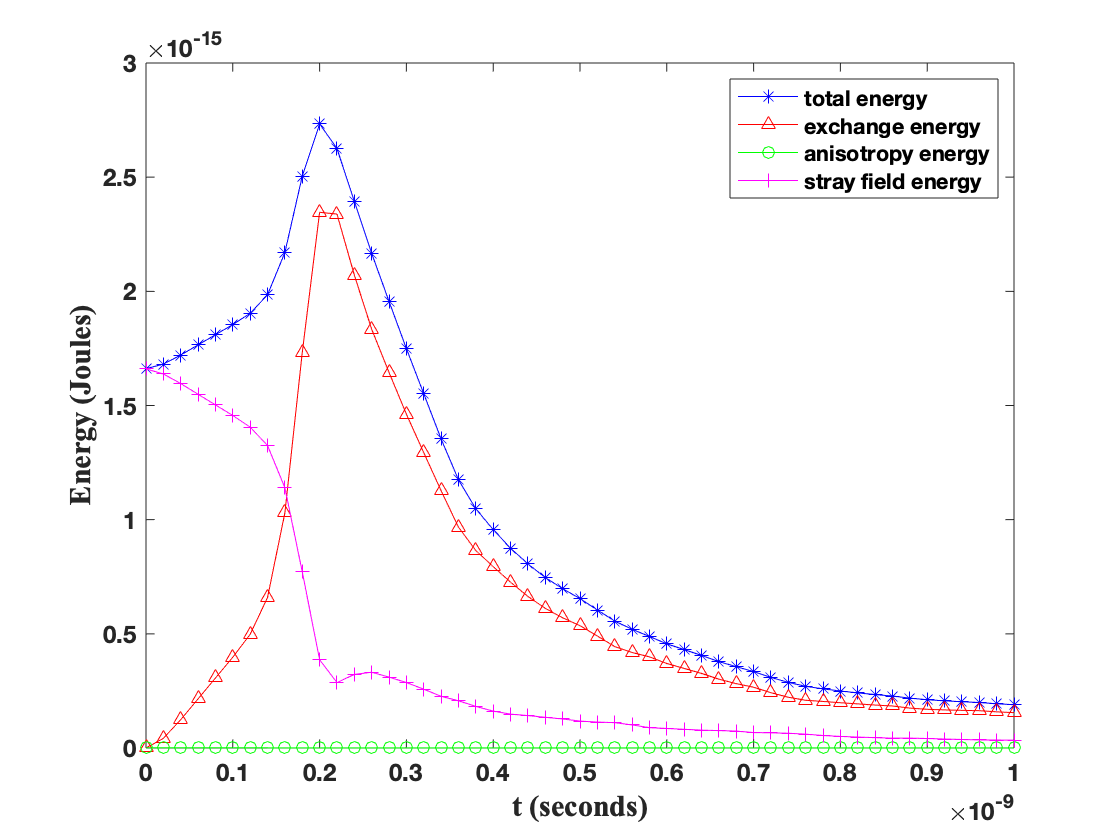}}
     \subfloat[SGS, $\alpha=0.1$]{\includegraphics[width=0.25\linewidth]{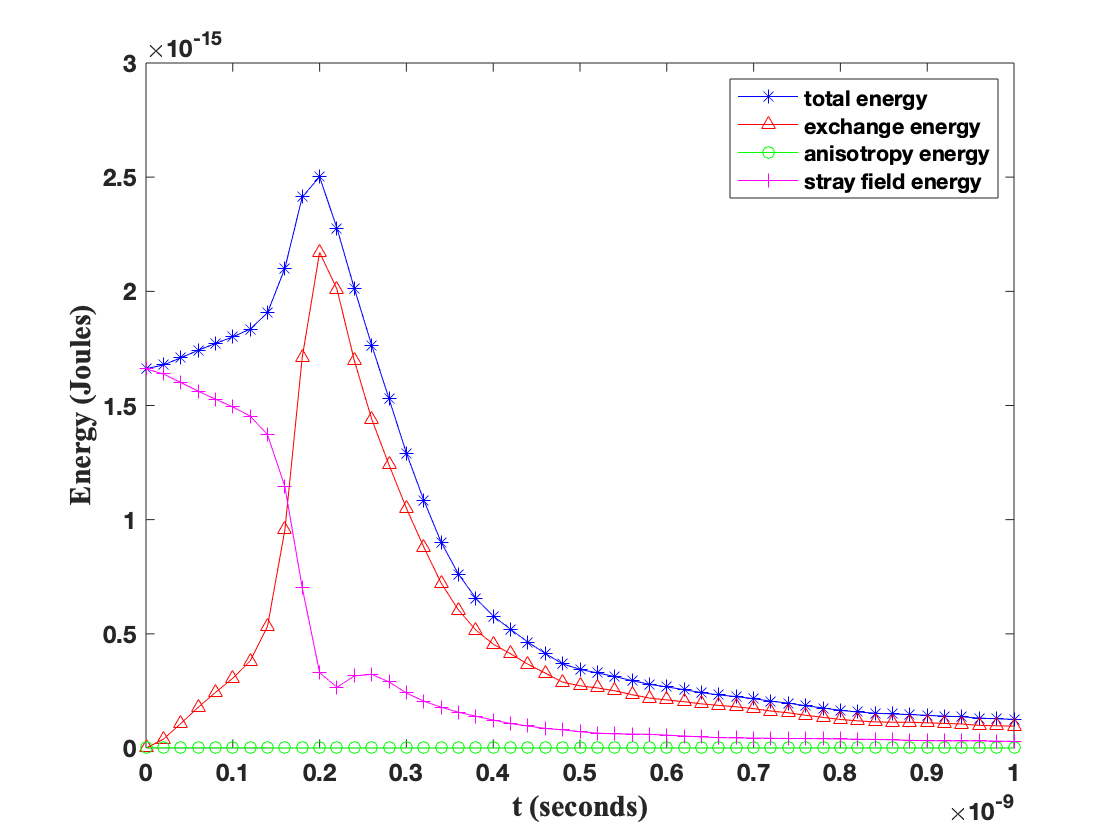}}
      \subfloat[SGS, $\alpha=0.5$]{\includegraphics[width=0.25\linewidth]{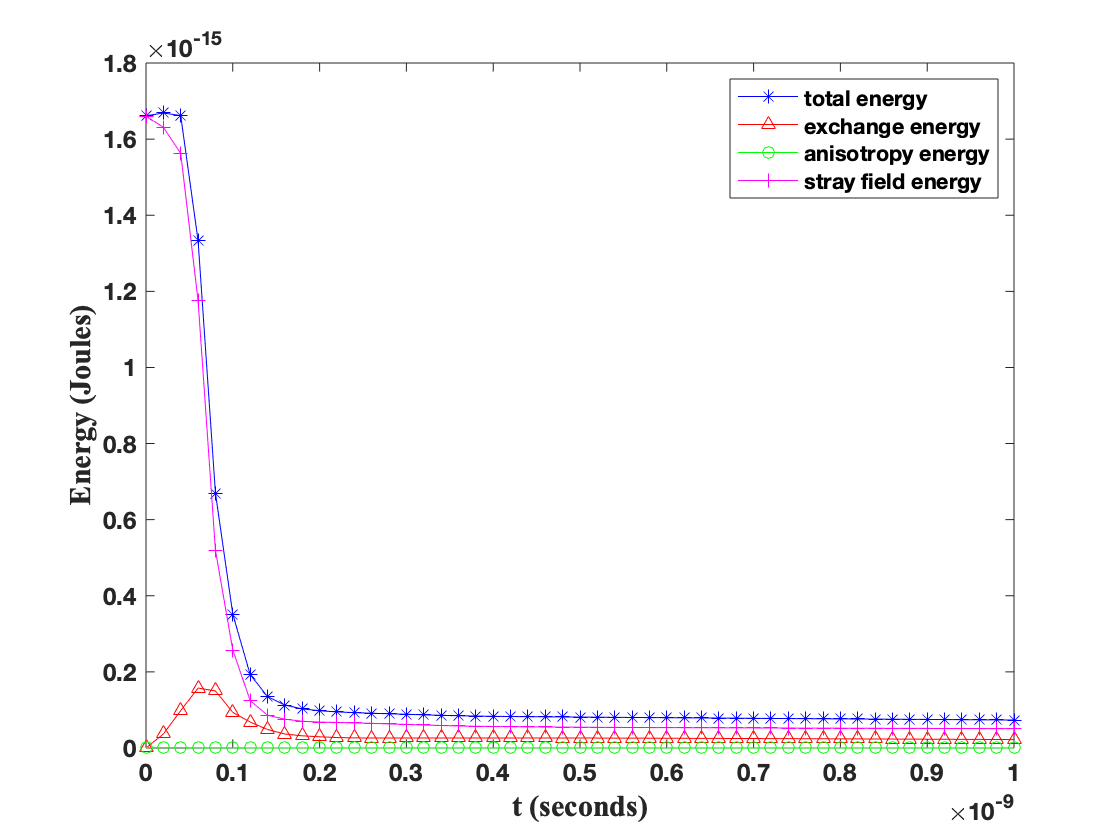}}
    \caption{Given inital condition $\m_0=[\cos(\cos(\pi x))\sin(0.01),\sin(\cos(\pi x))\sin(0.01),\cos(0.01)]^T$. The energy evolution for the total energy, the exchange energy, the anisotropy energy and stray field energy by GS and SGS projection method with $\alpha=0,0.01,0.1,0.5$.}
    \label{fig:energy-GS-SGS-v1-v1}
\end{figure}

\begin{figure}[htbp]
    \centering
    \subfloat[$\alpha=0$]{\includegraphics[width=0.4\linewidth]{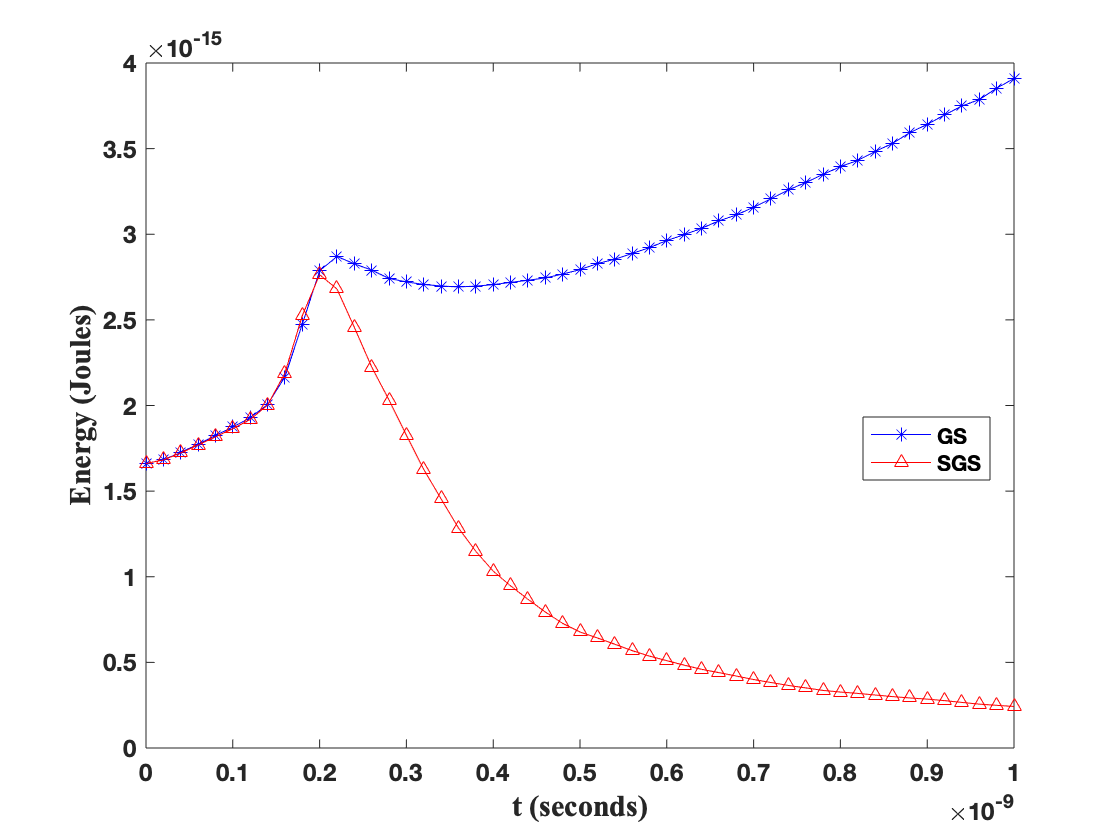}}
    \subfloat[$\alpha=0.01$]{\includegraphics[width=0.4\linewidth]{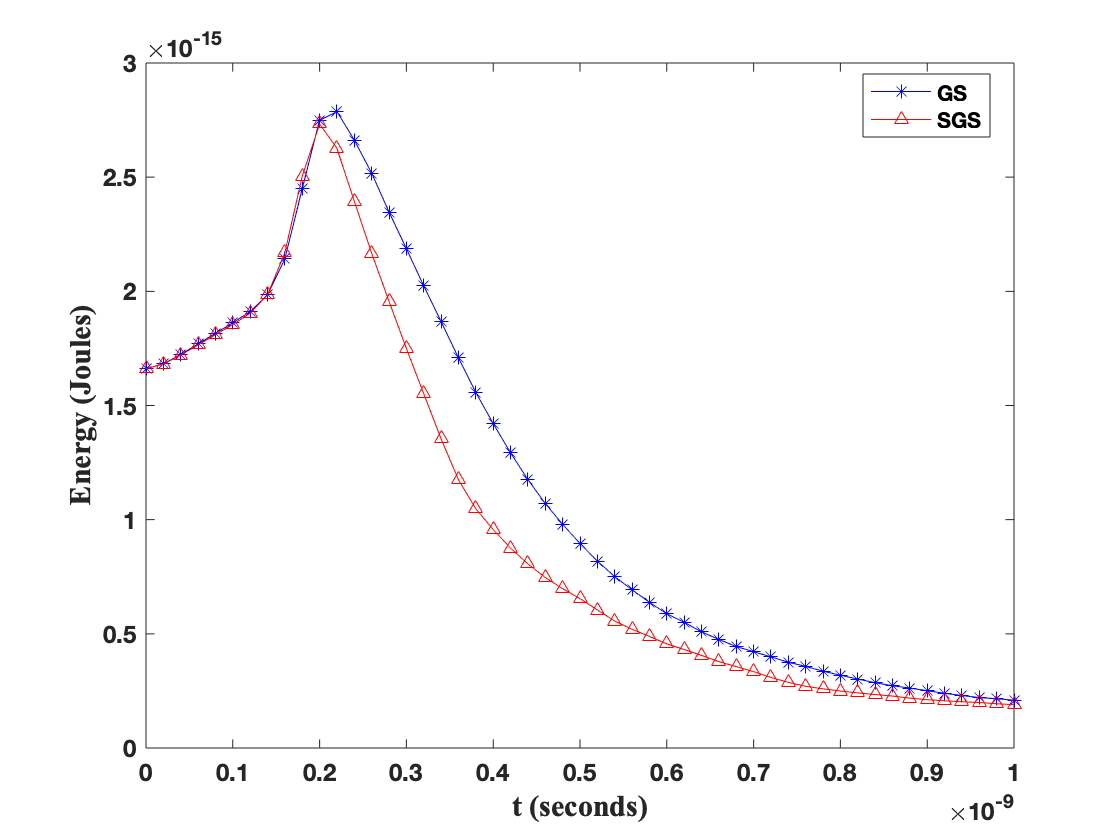}}
    \hspace{0.1in}
     \subfloat[$\alpha=0.1$]{\includegraphics[width=0.4\linewidth]{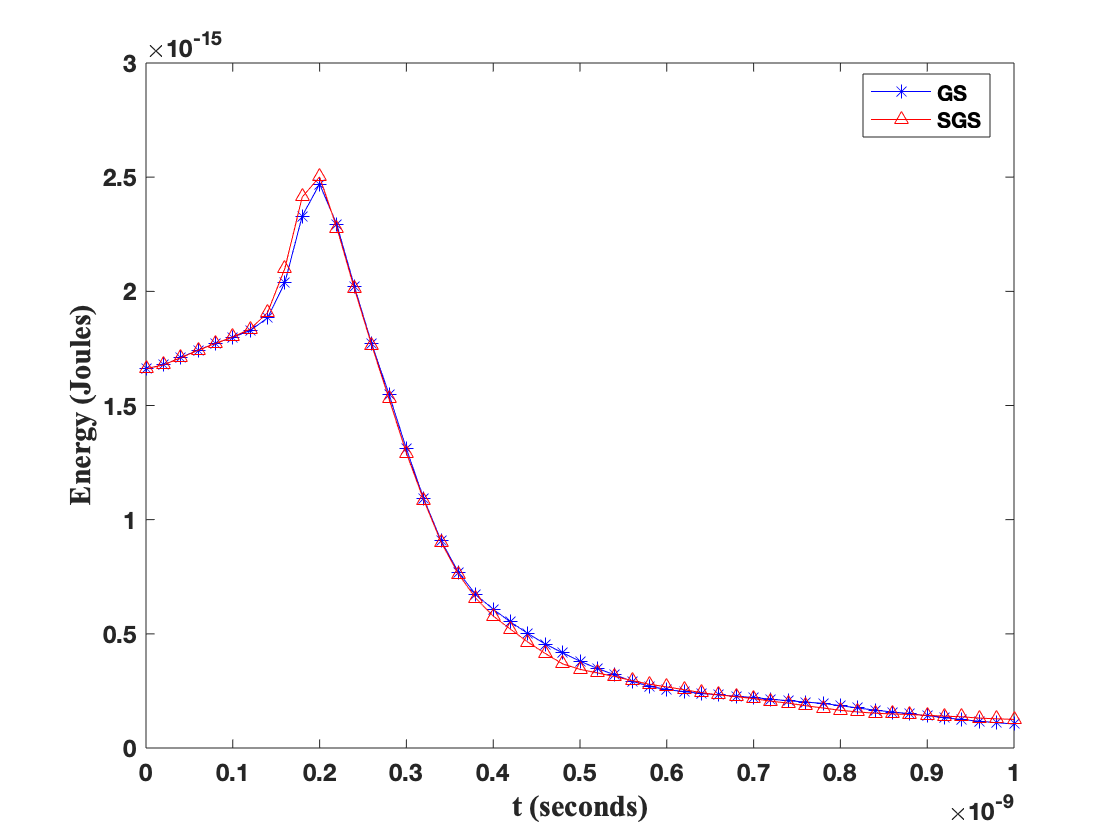}}
      \subfloat[$\alpha=0.5$]{\includegraphics[width=0.4\linewidth]{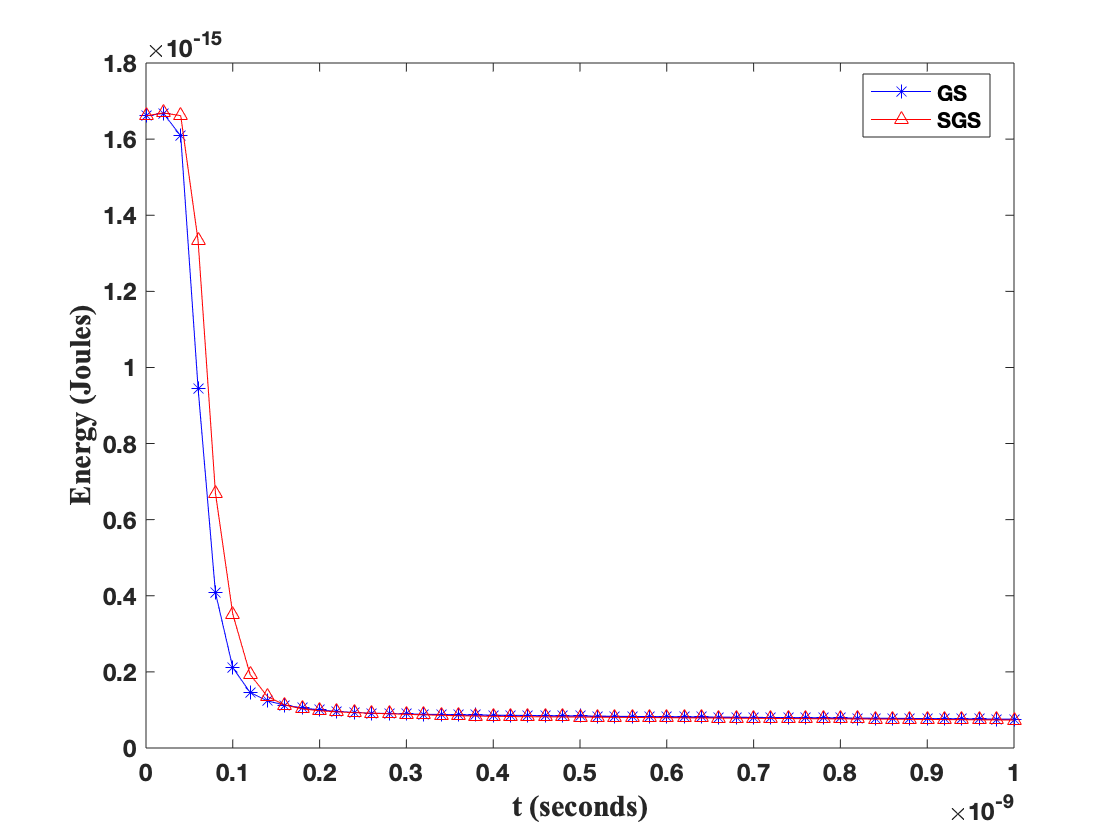}}
    \caption{Given inital condition $\m_0=[\cos(\cos(\pi x))\sin(0.01),\sin(\cos(\pi x))\sin(0.01),\cos(0.01)]^T$, the comparison of the energy evolution between GS and SGS projection method with $\alpha=0,0.01,0.1,0.5$ up to the final time $T=1\;ns$ with the time step size $\Delta t=1\;ps$.}
    \label{fig:energy-GS-SGS-v1-v2}
\end{figure}

The energy evolution given the initial condition \eqref{eq-m0-1} is presented in \Cref{fig:energy-GS-SGS-v1-v1} and \Cref{fig:energy-GS-SGS-v1-v2} for the comparison of GSPM and SGSPM. \Cref{fig:energy-GS-SGS-v1-v1} compares the evolution of total, exchange, anisotropy and stray-field energies computed by the standard GS (top row) and symmetric SGS (bottom row) projection schemes under four damping values \(\alpha=0,0.01,0.1,0.5\). For both methods, all energy components first rise to a peak then decay; larger \(\alpha\) accelerates energy dissipation and suppresses peak amplitudes. GS and SGS produce qualitatively identical energy profiles across tested positive damping coefficients, while different as $\alpha=0$. \Cref{fig:energy-GS-SGS-v1-v2} compares total energy evolutions of GS and SGS projection schemes under \(\alpha=0,0.01,0.1,0.5\). For all damping values, both methods share matching peak energy timing, yet SGS yields consistently lower energy afterward. Larger \(\alpha\) speeds up energy decay and narrows the energy gap between the two schemes.

\subsection{The simulations with only exchange field, anisotropy field, stray field, external field}

\begin{figure}[htbp]
    \centering
    \subfloat[Initial Ne\'el state]{\includegraphics[width=0.5\linewidth]{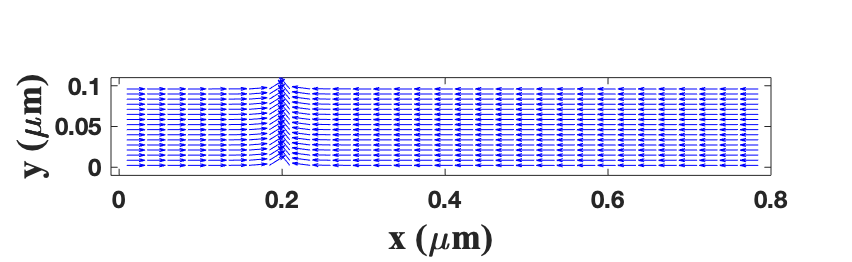}}
    \hspace{0.1in}
    \subfloat[GS, $T=0.2\;ns$]{\includegraphics[width=0.5\linewidth]{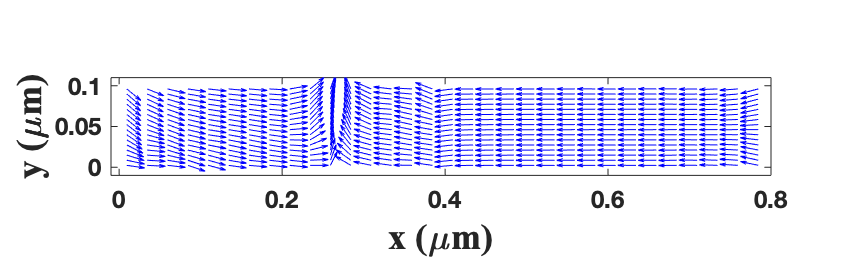}}
    \subfloat[SGS, $T=0.2\;ns$]{\includegraphics[width=0.5\linewidth]{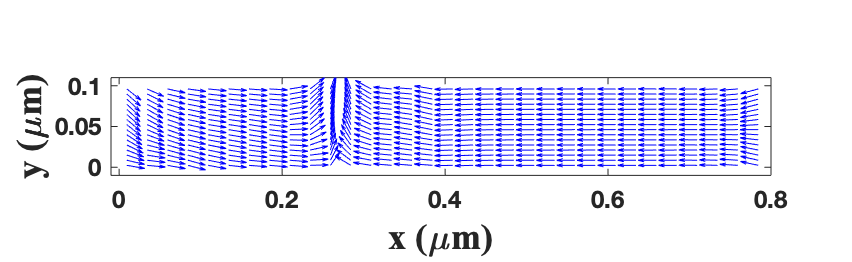}}
     \hspace{0.1in}
    \subfloat[GS, $T=0.4\;ns$]{\includegraphics[width=0.5\linewidth]{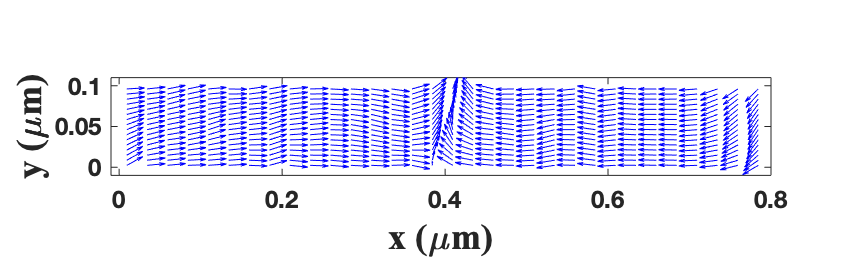}}
    \subfloat[SGS, $T=0.4\;ns$]{\includegraphics[width=0.5\linewidth]{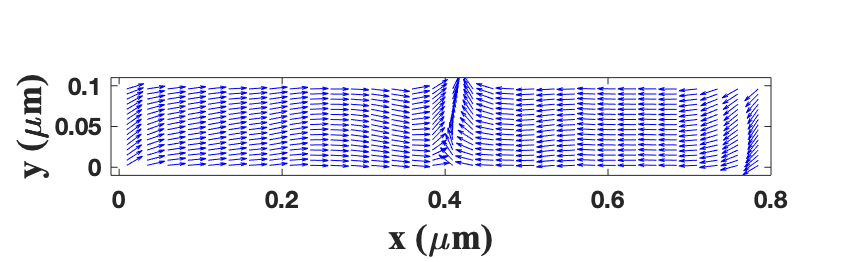}}
    \hspace{0.1in}
    \subfloat[GS, $T=0.6\;ns$]{\includegraphics[width=0.5\linewidth]{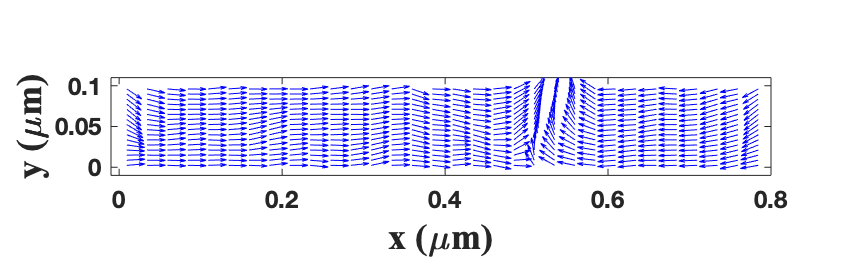}}
    \subfloat[SGS, $T=0.6\;ns$]{\includegraphics[width=0.5\linewidth]{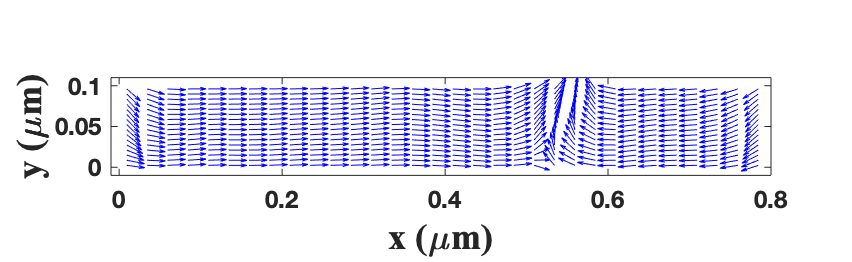}}
    \hspace{0.1in}
    \subfloat[GS, $T=0.8\;ns$]{\includegraphics[width=0.5\linewidth]{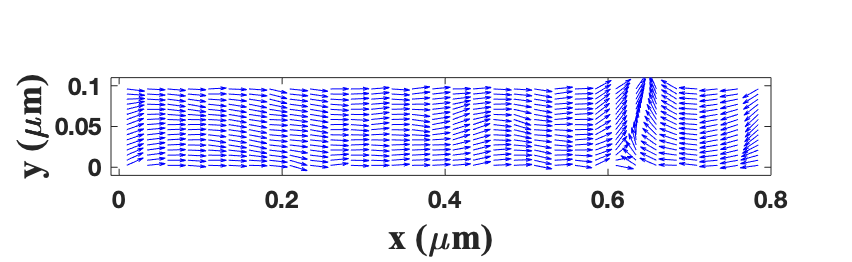}}
    \subfloat[SGS, $T=0.8\;ns$]{\includegraphics[width=0.5\linewidth]{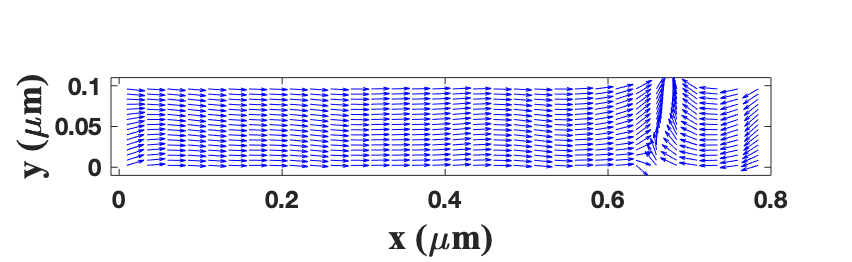}}
 \caption{The domain wall motion using GSPM and SGSPM with $\alpha=0$. Left panel with GSPM; Right panel with SGSPM.}
    \label{fig:DW-GS-SGS-v1-v2}
\end{figure}

\begin{figure}[htbp]
    \centering
    \subfloat[Initial Ne\'el state]{\includegraphics[width=0.5\linewidth]{Neel_GS_alpha_0dot01_he_5mT_initial_v1.png}}
    \hspace{0.1in}
    \subfloat[GS, $T=0.2\;ns$]{\includegraphics[width=0.5\linewidth]{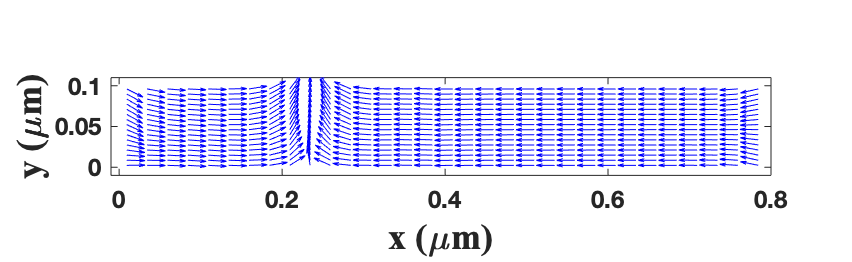}}
    \subfloat[SGS, $T=0.2\;ns$]{\includegraphics[width=0.5\linewidth]{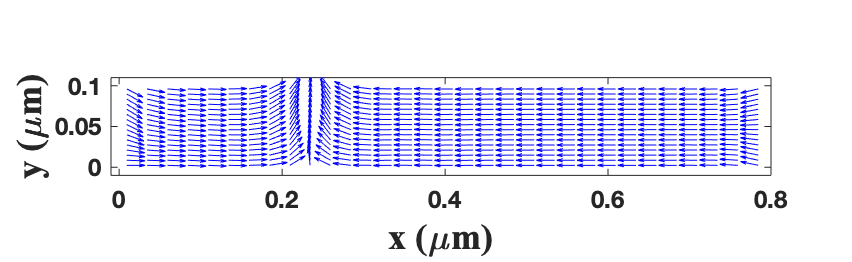}}
     \hspace{0.1in}
    \subfloat[GS, $T=0.4\;ns$]{\includegraphics[width=0.5\linewidth]{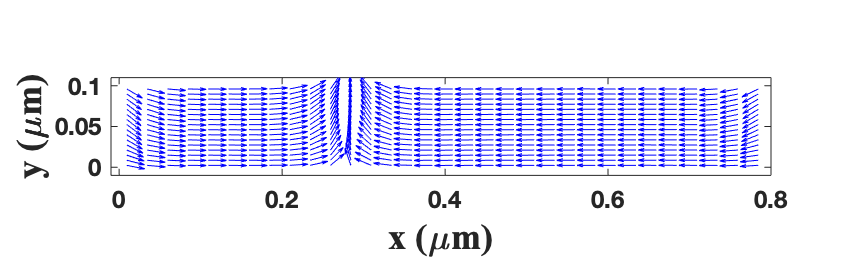}}
    \subfloat[SGS, $T=0.4\;ns$]{\includegraphics[width=0.5\linewidth]{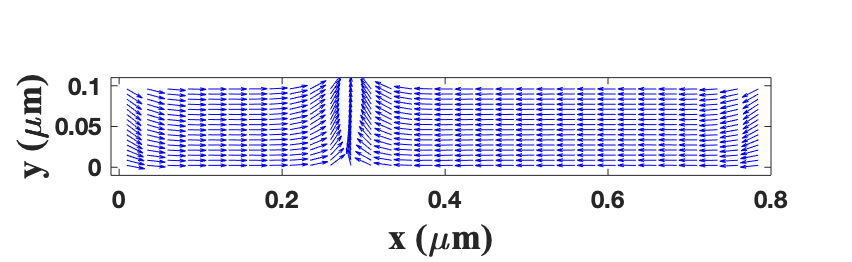}}
    \hspace{0.1in}
    \subfloat[GS, $T=0.6\;ns$]{\includegraphics[width=0.5\linewidth]{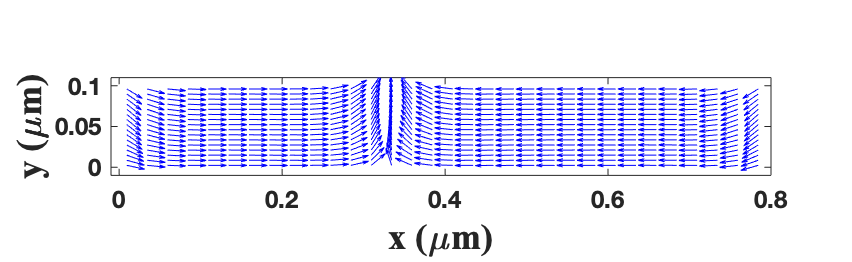}}
    \subfloat[SGS, $T=0.6\;ns$]{\includegraphics[width=0.5\linewidth]{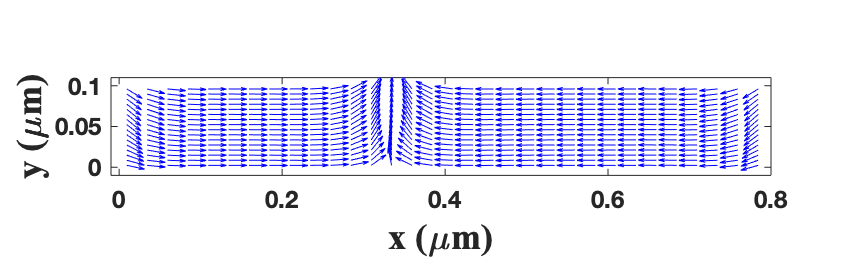}}
    \hspace{0.1in}
    \subfloat[GS, $T=0.8\;ns$]{\includegraphics[width=0.5\linewidth]{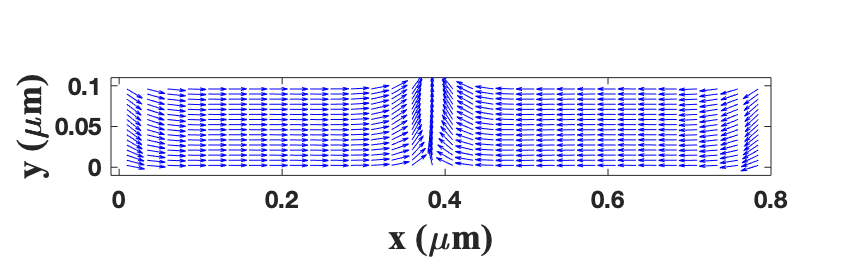}}
    \subfloat[SGS, $T=0.8\;ns$]{\includegraphics[width=0.5\linewidth]{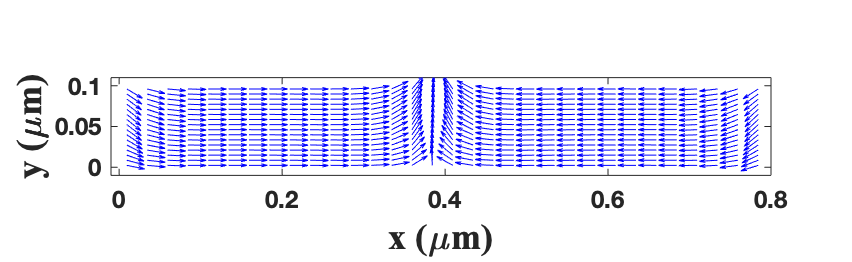}}
 \caption{The domain wall motion using GSPM and SGSPM with $\alpha=0.1$. Left panel with GSPM; Right panel with SGSPM.}
    \label{fig:DW-GS-SGS-v1-v3}
\end{figure}

In this section, we study the domain wall motion. A Ne\'el domain wall was initialized as the initial magnetic state within a ferromagnetic nanostrip of size $800 \times 100 \times 4\,\mathrm{nm}^3$. To ensure sufficient spatial resolution for capturing domain wall features while maintaining computational feasibility, the nanostrip was discretized using a structured grid of $128 \times 64 \times 4$ nodes. Following initialization, an external magnetic field of magnitude $\boldsymbol{h}_e=5\,\mathrm{mT}$ was applied along the positive $x$-direction to drive domain wall motion. Micromagnetic simulations of domain wall dynamics were performed over a time interval of up to $1\,\mathrm{ns}$. The results for GSPM and SGSPM given $\alpha=0$ and $\alpha=0.1$ are presented in \Cref{fig:DW-GS-SGS-v1-v2} and \Cref{fig:DW-GS-SGS-v1-v3}. The domain wall moves to right driven by the external field. For $\alpha=0$, SGSPM gives more stable results than that of GSPM. For the damping $\alpha=0.1$ domain wall propagation problem, SGSPM produces numerically consistent, visually matching magnetization solutions compared with the classic GSPM method at all examined time points, confirming that the symmetric iteration preserves the core micromagnetic dynamics without compromising solution quality.

\section{Conclusions and discussions}
\label{sec:conclusions}


In this paper, we propose a new stable Gauss-Seidel projection method (GSPM), namely the symmetric Gauss-Seidel projection method (SGSPM), which can be solved by the fast Fourier transform. Such an efficient method is of first order accuracy in time, norm preserving. It is much more stable than the GSPM, especially for small damping parameters. The numerical tests among the exchange field, stray field, anisotropy field and external field indicate the effectiveness of the proposed method. Such a symetric Gauss-Seidel strategy can be extended to the second order methods proposed in \cite{li2026enhanced} and two improved GSPM proposed in \cite{Li2020TwoIG} (the efficiency improvement) directly to improve the stability under weak damping parameters. Current work indicates that the model with only exchange field using our proposed method gives the energy decay property, whereas the popular GSPM is unstable with zero-damping. The enriched experiments indicate that representative stray field affects the energy evolution. The improvement of the energy decay of the full model will be studied in the future work.

\section*{Data availability}
The data will be made available on reasonable request.

\section*{Conflict of Interest Statement}
The authors have no conflicts of interest to declare. 

\section*{Acknowledgments}
This work is partially supported by the Basic Research Program of Jiangsu Province under Grant BK20250468, and the Research and Development Fund of XJTLU under Grant RDF-24-01-015.

\vspace{1cm}

\bibliographystyle{elsarticle-num-names}
\bibliography{references.bib}

\end{document}